\documentclass{book}

\usepackage[T1]{fontenc}
\usepackage[utf8]{inputenc}
\usepackage[UKenglish]{babel}
\usepackage[UKenglish]{isodate}

\usepackage[margin=1.5in]{geometry}
\usepackage{graphicx}
\usepackage{enumitem}
\usepackage{csquotes}

\usepackage[
type={CC},
modifier={by-nc-nd},
version={4.0},
]{doclicense}

\usepackage{amsmath,amssymb,amsthm}

\usepackage[
natbib,
maxcitenames=3, 
giveninits, 
maxbibnames=10,  
sorting=nyt,
url=false,
doi=false,
sortcites,
backref,
backend=biber
]{biblatex}
\renewbibmacro{in:}{}

\makeatletter
\let\blx@rerun@biber\relax
\makeatother

\usepackage{hyperref}

\usepackage[textwidth=0.8\marginparwidth]{todonotes}
\setuptodonotes{inline}
\usepackage{url}

\usepackage[nameinlink, capitalise, noabbrev]{cleveref}

\usepackage{subcaption}
\usepackage{xfrac}
\usepackage{nicefrac}
\usepackage{tikz}
\usepackage{pgfplots}
\pgfplotsset{compat=1.17}
\usetikzlibrary{3d,calc}
\usepackage[overload]{empheq}

\newcommand{\Hspace}{\mathcal H}
\newcommand{\Bspace}{\mathcal X}

\newtheorem{theorem}{Theorem}
\newtheorem{proposition}[theorem]{Proposition}
\newtheorem{lemma}[theorem]{Lemma}
\newtheorem{corollary}[theorem]{Corollary}
\newtheorem{remark}[theorem]{Remark}

\theoremstyle{definition}
\newtheorem{example}[theorem]{Example}
\newtheorem{definition}[theorem]{Definition}

\numberwithin{theorem}{chapter}
\numberwithin{equation}{chapter}

\newcommand{\N}{\mathbb{N}}

\newcommand{\R}{\mathbb{R}}
\newcommand{\E}{\mathcal{E}}
\newcommand{\dom}{{\mathrm{dom}}}
\newcommand{\closure}[1]{\overline{#1}}
\newcommand{\norm}[1]{\left\Vert #1 \right\Vert}
\newcommand{\abs}[1]{\left\vert #1 \right\vert}
\renewcommand{\H}{\mathcal{H}}

\renewcommand{\div}{\operatorname{div}}
\DeclareMathOperator{\argmin}{arg\,min}
\DeclareMathOperator{\argmax}{arg\,max}
\DeclareMathOperator{\esssup}{ess\,sup}
\newcommand{\st}{\,:\,}
\newcommand{\supp}{\operatorname{supp}}
\newcommand{\dx}{\,\mathrm{d}x}
\renewcommand{\d}{\,\mathrm{d}}

\newcommand{\sign}{\operatorname{sign}}

\newcommand{\calA}{\mathcal{A}}
\newcommand{\calL}{\mathcal{L}}
\newcommand{\calN}{\mathcal{N}}

\newcommand{\calJ}{\mathcal{J}}
\newcommand{\calC}{\mathcal{C}}

\newcommand{\calR}{\mathcal{R}}
\newcommand{\calD}{\mathcal{D}}
\newcommand{\eps}{\varepsilon}

\newcommand{\bv}{\mathrm{BV}}
\DeclareMathOperator{\tv}{TV}

\newcommand{\dist}{\operatorname{dist}}

\newcommand{\subdiff}{\partial}
\newcommand{\grad}{\nabla}
\newcommand{\defeq}{:=}
\newcommand{\tex}{T_\mathrm{ex}}

\makeatletter
\newsavebox{\@brx}
\newcommand{\llangle}[1][]{\savebox{\@brx}{\(\m@th{#1\langle}\)}%
  \mathopen{\copy\@brx\kern-0.5\wd\@brx\usebox{\@brx}}}
\newcommand{\rrangle}[1][]{\savebox{\@brx}{\(\m@th{#1\rangle}\)}%
  \mathclose{\copy\@brx\kern-0.5\wd\@brx\usebox{\@brx}}}
\makeatother

\let\sp\relax

\newcommand{\wsto}{\rightharpoonup^\ast}
\newcommand{\weakto}{\rightharpoonup}

\newcommand{\tauto}{\xrightarrow{\tau}}
\newcommand{\lsc}{lower-semi\-continuous}

\newcommand{\ef}{eigenfunction}
\newcommand{\ev}{eigenvalue}
\newcommand{\dJ}{\subdiff \calJ}
\newcommand{\sp}[1]{\langle #1 \rangle}
\newcommand{\dualball}{K_{\calJ}}
\newcommand{\nullspace}{{\calN}}
\newcommand{\laplace}{\mathop{}\!\mathbin\bigtriangleup}
\DeclareMathOperator{\prox}{prox}
\newcommand{\RI}{{\R \cup \{+\infty\}}}
\newcommand{\charf}{\chi}
\newcommand{\compemb}{\subset\subset}
\newcommand{\contemb}{\subset}

\newcommand{\innerp}[2]{\left\langle #1, #2 \right\rangle}
\newcommand{\scalprod}{\innerp}
\newcommand{\FunA}{E}
\newcommand{\FunB}{F}
\newcommand{\FunFromTo}[3]{#1\colon#2\rightarrow#3}
\newcommand{\dual}[1]{#1^\ast}
\newcommand{\conj}[1]{#1^\ast}
\newcommand{\biconj}[1]{#1^{\ast\ast}}

\newcommand{\imageA}{u}
\newcommand{\imageB}{v}

\newcommand{\subgradA}{\zeta}
\newcommand{\subgradB}{\eta}
\newcommand{\subgradC}{\xi}

\newcommand{\spaceA}{\mathcal H}

\graphicspath{{figures/}}

\begin{document}

\begin{titlepage}
\newgeometry{top=2in,bottom=1cm,right=1in,left=1in}
   \begin{center}
       \vspace*{1cm}
 
       \textbf{\huge{Introduction to  Nonlinear Spectral Analysis}}
 
       \vspace{1cm}
        These notes are based on the lectures taught by the authors \\
        at the universities of Bonn and Cambridge in 2022
 
       \vspace{1.5cm}
 
       \textbf{Leon Bungert and Yury Korolev}
 
      \vspace{9cm}
 Latest update: \today
 \vspace{3cm}
    \doclicenseThis
   \end{center}  
\end{titlepage}

\tableofcontents

\chapter*{Introduction}

Eigenfunctions of linear operators, i.e. nontrivial solutions of $\lambda u = Au$, and corresponding spectral decompositions have found numerous applications in partial differential equations (PDEs), signal processing and many other fields. Perhaps the most common example is the Fourier transform
which can be used to decompose a signal into a sum of trigonometric functions---eigenfunctions of the
Laplace operator on the unit cube.

While this is useful for some types of signals, e.g. audio signals, the trigonometric basis is not well suited for signals with discontinuities, such as images. In this case, it turns out advantageous to use an appropriate generalization of the concept of an \ef{} to nonlinear operators such as the 1-Laplacian~\cite{gilboa2018nonlinear}. Under an \ef{} of a nonlinear---or even set-valued---operator $\calA$ we understand a solution of the inclusion $\lambda u \in \calA(u)$. A~prototypical example of such an operator is the subdifferential of a convex function $\calA(u) = \partial \calJ(u)$. In the case of the 1-Laplacian, this functional is the Total Variation~\cite{rudin1992nonlinear}. 

Unlike the linear case, the explicit form of nonlinear \ef{s} is rarely known and the spectral decomposition of a given signal $f$ is usually computed by solving a \mbox{(sub-)gradient} flow $u_t \in -\calA(u)$ with the initial condition $u(0) = f$. For images this was first proposed in~\cite{gilboa2013spectral, gilboa2014total}, based on the Total Variation flow~\cite{andreu2002some}, and further developed in~\cite{gilboa2016nonlinear, burger2016spectral, bungert2019nonlinear, bungert2019computing}. 
While such decompositions very rarely are constituted of actual eigenfunctions of the nonlinear operator at hand (see \cite{bungert2019nonlinear} for a comprehensive characterization of these rare cases), they still can be successfully applied for certain analysis and synthesis tasks in signal and image processing, e.g., in image denoising~\cite{moeller2015learning}, image fusion~\cite{benning2017nonlinear} and shape processing~\cite{fumero2020nonlinear,brokman2021nonlinear,brokman2024spectral}.

Apart from spectral decompositions, the so-called ground states---suitably normalized \ef{s} with the smallest value \ev{} $\lambda$---are of independent interest. In the case when $\calA$ is the subdifferential of a convex absolutely one-homogeneous functional, ground states were first studied in the context of variational regularisation in~\cite{benning2013ground}. They play an important role in graph clustering~\cite{buhler2009spectral, elmoataz2015p} and have been linked to the asymptotic behavior of time-dependent equations~\cite{andreu2002some, varvaruca2004exact, aujol2018theoretical, feld2019rayleigh, bungert2019asymptotic}. Ground states of certain $L^\infty$-type functionals are connected to the infinity Laplacian and distance functions~\cite{bungert2020structural, bungert2022eigenvalue,bozorgnia2024infinity}.
Furthermore, also the nonlinear spectral theory on weighted graphs \cite{deidda2025nonlinear,bungert2020structural,mazon2020total} as well as duality for nonlinear eigenvalue problems \cite{tudisco2022nonlinear,bungert2022eigenvalue} have been investigated recently.

These notes are meant as an introduction to the theory of nonlinear \ef{s}. We will discuss the variational form of  eigenvalue problems and the corresponding non-linear Euler--Lagrange equations, as well  as connections with gradient flows. For the latter ones, we will give precise conditions for finite extinction and discuss convergence rates. We will use this theory to study asymptotic behaviour of nonlinear PDEs and present applications in $L^\infty$ variational problems. Finally we will discuss numerical methods for solving gradient flows and computing nonlinear eigenfunctions based on a nonlinear power method.

Our main tools are convex analysis and calculus of variations, necessary background on which will be provided. It is expected that the reader is familiar with Hilbert spaces; familiarity with Banach spaces is beneficial but not strictly necessary.

The notes are based on the lectures taught by the authors at the universities of Bonn and Cambridge in 2022. 
While here the focus is on the Hilbert space theory, the approach can be generalized to operators defined on a Banach space, see~\cite{bungert2021gradient, bungert2022eigenvalue}.

\paragraph{Acknowledgments}

Parts of these notes were written while we were in residence at the Institut Mittag-Leffler in Djursholm, Sweden during the semester on \textit{Geometric Aspects of Nonlinear Partial Differential Equations} in 2022, supported by the Swedish Research Council under grant no. 2016-06596.
Finally, we would like to thank the many people who shaped our views on the subject, including Martin Burger, Guy Gilboa, Carola Schönlieb, Martin Benning, Michael Möller, and José Mazón.

\chapter[Eigenfunctions of the Laplace operator]{Motivation and Recap: Eigenfunctions of the Laplace Operator}
\label{ch:lin-ef}

In this section we will recall some results from linear spectral theory, which we will generalise in the next chapter. For simplicity, we will work with one particular example, the Laplace operator $\laplace \,\cdot \defeq \div \, (\grad \,\cdot)$. 

We start by recalling some necessary concepts from functional analysis and calculus of variations.
In these notes we let $\Bspace$ denote a generic Banach space and $\Hspace$ a generic Hilbert space.

\section{Background}\label{sec:laplace-background}
We assume that the reader is familiar with Hilbert and Banach spaces. 
We also assume basic knowledge on the $L^p(\Omega)$ spaces on a domain $\Omega\subset\R^d$. 
$L^2(\Omega)$ is a Hilbert space equipped with the inner product $\langle u,v\rangle_{L^2}:=\int_\Omega uv\d x$.

\subsubsection{Sobolev spaces}
Let $\Omega \subset \R^d$ be an open bounded domain. We will denote $\N_0 \defeq \N \cup \{0\}$. For any multi-index $\alpha = (\alpha_1,...,\alpha_d) \subset \N_0^d$ of order $k \in \N$, i.e. such that $\abs{\alpha} \defeq \sum_{i=1}^d \alpha_i = k$, we define the following differential operator
\begin{equation*}
    D^\alpha \defeq \frac{\partial^{\abs{\alpha}}}{\partial^{\alpha_1}_{x_1} \dots \partial^{\alpha_d}_{x_d}}.
\end{equation*}

Let $C_c^\infty(\Omega)$ be space of smooth compactly supported functions on $\Omega$.
\begin{definition}
Let $\Omega \subset \R^d$ be an open bounded domain, $p\in[1,\infty]$, and $\alpha \subset \N_0^d$ a multi-index. A function $f \in L^p(\Omega)$ is said to have a weak derivative $D^\alpha f \in L^p(\Omega)$ or order $\alpha$ if the following integral exists for any function $\phi \in C_c^\infty(\Omega)$
\begin{equation*}
     (-1)^{\abs{\alpha}} \int_\Omega f \, D^\alpha \phi \d x \defeq \int_\Omega \phi \, D^\alpha f \d x, \quad \phi \in C_c^\infty(\Omega).
\end{equation*}
\end{definition}

Now we are ready to define \emph{Sobolev spaces}, the spaces of Lebesgue integrable functions whose weak derivatives are also Lebesgue integrable.
\begin{definition}
Let $\Omega \subset \R^d$ be an open bounded domain. The Sobolev space $W^{k,p}(\Omega)$ of order $k \in \N$ is defined as follows
\begin{equation*}
    W^{k,p}(\Omega) \defeq \{ f \in L^p(\Omega) \colon  D^\alpha f \in L^p(\Omega) \text{ for all multi-indices $\alpha$ s.t. $\abs{\alpha} \leq k$} \}.
\end{equation*}
The Sobolev spaces $H^k(\Omega):=W^{k,2}(\Omega)$ are Hilbert spaces with the following inner product
\begin{equation*}
    \sp{f,g}_{H^k} \defeq \sum_{\abs{\alpha} \leq k} \sp{\partial_x^{\alpha} f,\partial_x^{\alpha} g}_{L^2} = \sum_{\abs{\alpha} \leq k} \int_\Omega \partial_x^{\alpha} f \; \partial_x^{\alpha} g \d x.
\end{equation*}
\end{definition}

It can be shown that Sobolev functions of order $k \in \N$ have well-defined values on subsets of $\Omega$ of codimension $k$.  In particular, functions in $W^{1,p}(\Omega)$ have well-defined values on submanifolds with dimension $d-1$. 

\begin{definition}
Let $\Omega \subset \R^d$ be an open bounded domain. The space $W^{k,p}_0(\Omega)$ is the space Sobolev functions satisfying zero boundary conditions defined as the following closed subspace of $W^{k,p}(\Omega)$:
\begin{equation*}
    W^{k,p}_0(\Omega) \defeq \closure{C^\infty_c(\Omega)}^{W^{k,p}(\Omega)}.
\end{equation*}
We also define $H^{k}_0(\Omega):=W^{k,2}_0(\Omega)$.
\end{definition}

The following result will be very useful to us.
\begin{theorem}[Rellich–Kondrachov embedding theorem]\label{thm:rellich-kondrachov}
    Let $\Omega \subset \R^d$ be an open bounded connected domain with Lipschitz boundary.
    Let $1\leq p \leq d$ and $p^*:= \frac{dp}{d-p}$ if $p<d$ and $p^*:=\infty$ if $p = d$.
    Then the embedding
    $W^{1,p}(\Omega) \compemb L^{q}(\Omega)$
    is compact for all $1\leq q < p^*$ (in particular, for $q:=p$).
\end{theorem}
\begin{corollary}\label{cor:rellich}
The embedding $W^{1,p}_0(\Omega)\compemb L^q(\Omega)$ is compact for all $1\leq q < p^*$ (in particular, for $q:=p$).
\end{corollary}

\subsubsection{Laplace operator}
Let $\Omega \subset \R^d$ be an open bounded domain and $L^2(\Omega)$ be the Lebesgue space of square-integrable functions on $\Omega$. If for a function $u\in L^2(\Omega)$ there exists a vector-valued function $g\in L^2(\Omega;\R^d)$ such that
\begin{align}\label{eq:grad_weak}
 \int_\Omega u \, \grad\phi \d x = -\int_\Omega \phi \, g \d x,\quad\forall\phi\in C^\infty_c(\Omega),   
\end{align}
we write $\grad u := g$. 

Similarly, if for a vector-valued function $f\in L^2(\Omega;\R^d)$ there exists a scalar-valued function $v\in L^2(\Omega)$ such that
\begin{align}\label{eq:div_weak}
 \int_\Omega f(x) \cdot \grad\phi(x) \d x = -\int_\Omega  \phi v \d x,\quad\forall\phi\in C^\infty_c(\Omega),   
\end{align}
we write $\div f := v$.

We will also consider the Laplace operator $\laplace \,\cdot \defeq \div \, (\grad \,\cdot)$. 
If for a function $u\in L^2(\Omega)$ there exists a function $g\in L^2(\Omega)$ satisfying
\begin{align}\label{eq:laplace_weak}
 \int_\Omega u\laplace\phi\d x = \int_\Omega g\,\phi\d x,\quad\forall\phi\in C^\infty_c(\Omega),   
\end{align}
we write $\laplace u := g$ and note such $g$ is uniquely determined if it exists.
In this case we write $\laplace u\in L^2(\Omega)$.
Hence, the Laplace operator acts as follows:
\begin{equation*}
    \laplace \colon H^2(\Omega)\cap H^1_0(\Omega) \to L^2(\Omega), \quad \laplace u \defeq g
    \text{ where $g$ satisfies \eqref{eq:laplace_weak}.}
\end{equation*}

It can be easily checked that the negative Laplace operator $-\laplace$ is self-adjoint and positive definite. 
Its \ef{s} are solutions to
\begin{equation}\label{eq:helmholtz}
    -\laplace u = \lambda u,\quad u\in H^2(\Omega)\cap H^1_0(\Omega),\;\lambda\in\R.
\end{equation}

\begin{theorem}\label{thm:laplace-efs}
Let $\Omega \subset \R^d$ be an open bounded domain. 
Then there exists an orthonormal basis of $L^2(\Omega)$ consisting of solutions $\{u_i\}_{i\in\N}$ of~\eqref{eq:helmholtz} with real and positive \ev{s} $\{\lambda_i\}_{i\in\N} \subset \R_+$ satisfying $0<\lambda_1\leq\lambda_2\leq\dots$ and $\lim_{i\to\infty}\lambda_i = \infty$.
\end{theorem}

\subsubsection{The extended real line}

We will find it useful to work with the \emph{extended real line} $\RI \defeq \R \cup \{\pm\infty\}$.
\begin{definition} The extended real line is defined as $\RI \defeq \R \cup \{\pm\infty\}$ with the following rules that hold for any $a \in \R$ and $\alpha > 0$:
\begin{align*}
    & a \pm \infty \defeq \pm\infty + a \defeq \pm\infty , \\
    & \alpha \cdot (\pm\infty) \defeq \pm\infty \cdot \alpha \defeq \pm\infty \\
    & a / (\pm\infty) \defeq 0 \\ 
    & \pm\infty \pm \infty \defeq \pm\infty \\
    & \pm\infty \cdot \pm\infty \defeq \infty\\
    & \pm\infty \cdot \mp\infty \defeq -\infty.
\end{align*}
\end{definition}
\noindent Some operations are \emph{not defined}, e.g.,
\begin{equation*}
    +\infty - \infty, \quad \frac{\pm\infty}{\pm\infty}.
\end{equation*}

\begin{definition}
    Let $\calC$ be a subset of $\Bspace$. The \emph{characteristic function} $\charf_{\calC} \colon \Bspace \to \RI$ of $\calC$ is defined as follows
\begin{equation*}
    \charf_{\calC}(u) \defeq
    \begin{cases}
        0, \quad & u \in \calC, \\
        +\infty \quad &\text{otherwise}.
    \end{cases}
\end{equation*}
\end{definition}

We would like to exclude pathological cases when $E$ is equal to either $+\infty$ or $-\infty$ everywhere.
\begin{definition}
A functional $E \colon \Bspace \to \RI$ is called \emph{proper} if $\exists u \in \Bspace$ such that $E(u) \neq \pm \infty$. The set $\dom(E) \defeq \{u \in \Bspace \colon E(u) \neq \pm \infty\}$ is called the effective domain of $E$.
\end{definition}

\paragraph{Dual spaces}

\begin{definition}
    The dual space of a Banach space $\Bspace$ is given by the collection of all bounded linear functionals on $\Bspace$ and denoted by
    \begin{align*}
        \dual{\Bspace} := \big\{\phi:\Bspace\to\R \st \phi\text{ is linear and $\exists \,C\geq 0$ such that }\abs{\phi[u]}\leq C\norm{u}_\Bspace\,\forall u\in\Bspace\big\}.
    \end{align*}
\end{definition}
The dual space can be equipped with the dual norm 
\begin{align*}
    \norm{\phi}_{\dual{\Bspace}} := \sup_{\norm{u}_\Bspace=1} \phi[u].
\end{align*}
For a Hilbert space $\Hspace$ by the Riesz representation theorem $\Hspace^*$ is a Hilbert space isomorphic (but not equal!) to $\Hspace$.

If no confusion is possible, we will replace $\norm{\cdot}_\Bspace$ by $\norm{\cdot}$ and $\norm{\cdot}_{\dual{\Bspace}}$ by $\norm{\cdot}_\ast$. 

\subsubsection{First variation}

We define a generalisation of the concept of a directional derivative known as the first variation. It is related, but not equivalent to the concept of a Gateaux derivative.
\begin{definition}[first variation]\label{def:first_var}
Let $E \colon \Bspace \to \RI$ a function. 
Let $u,v \in \Bspace$ with $E(u)\neq\pm\infty$.
Then the first variation of $E$ at $u$ in the direction of $v$ is defined as the limit
\begin{align*}
\delta E(u;v) := \lim_{\eps \to 0}\frac{E(u + \eps v) - E(u)}{\eps}
\end{align*}
if the limit exists in $\RI$.
\end{definition}

\begin{proposition}[chain rule]\label{prop:chain_rule_first_var}
Consider two functions $E \colon \Bspace \to \RI$ and $\varphi \colon \RI \to \RI$. 
Let $u,v \in \Bspace$ with $E(u)\neq\pm\infty$.
If $\delta E(u;v)$ exists in $\R$ and $\varphi$ is differentiable at $E(u)$ then the first variation of the composition $\varphi \circ E \colon \Bspace \to \R$ exists and is given by
\begin{equation*}
    \delta(\varphi \circ E)(u;v) = \varphi^\prime(E(u)) \, \delta E(u;v).
\end{equation*}
\end{proposition}
\begin{proof}
Define an auxiliary function as
\begin{align*}
    \psi(s) :=
    \begin{cases}
        \frac{\varphi(s) - \varphi(E(u))}{s - E(u)},\quad&\text{if }s\neq E(u), \\
        \varphi'(E(u)),\quad&\text{else},
    \end{cases}
    \qquad s\in\R.
\end{align*}
Since $E(u)\neq\pm\infty$ and $\varphi$ is differentiable at $E(u)$, $\psi$ is continuous at $E(u)$.
Furthermore, since $\delta E(u;v)$ exists in $\R$, we know that $\eps\mapsto E(u+\eps v)$ is in particular continuous in $\eps=0$.
As a composition of continuous maps, the map $\eps\mapsto\psi(E(u+\eps v))$ is continuous at $\eps=0$.

Using the definition of $\psi$, for all $\eps\neq 0$ it holds
\begin{align*}
    \frac{\varphi(E(u+\eps v)) - \varphi(E(u))}{\eps} = 
    \psi(E(u+\eps v))\frac{E(u+\eps v)-E(u)}{\eps}.
\end{align*}
Hence, using that $\delta E(u;v)$ exists, that $\eps\mapsto \psi(E(u+\eps v))$ is continuous at $\eps=0$, and that $\psi$ is continuous at $E(u)$ we obtain
\begin{align*}
    \delta(\varphi\circ E)(u;v) 
    &= 
    \lim_{\eps\to 0}
    \frac{\varphi(E(u+\eps v)) - \varphi(E(u))}{\eps}
    =
    \lim_{\eps\to 0}
    \psi(E(u+\eps v))\frac{E(u+\eps v)-E(u)}{\eps}
    \\
    &=
    \psi(E(u))\delta E(u;v)
    =
    \varphi'(E(u))\delta E(u;v).
\end{align*}
\end{proof}

\begin{proposition}[product rule]\label{prop:prod_rule_first_var}
Let $E,F \colon \Bspace \to \RI$ and $u,v \in \Bspace$ with $E(u),F(u)\neq\pm\infty$.  
If the first variations of $E$ and $F$ at $u$ exist, then the first variation of the product $E\cdot F:\Bspace\to\RI$ exists at $u$ and is given by
\begin{equation*}
    \delta(E\cdot F)(u;v) = \delta E(u;v)\,F(u) + E(u)\,\delta F(u;v).
\end{equation*}
\end{proposition}

\subsubsection{Fr\'{e}chet derivative}

We will also need the following stronger concept of a derivative known as the Fr\'{e}chet derivative.

\begin{definition}[Fr\'echet derivative]
Let $E \colon \Bspace \to \R$ a mapping and $u \in \Bspace$. If there exists a bounded linear operator $A\in\dual{\Bspace}$ such that
\begin{align*}
\lim_{\norm{h} \to 0} \frac{\abs{E(u + h) - E(u) - A[h]}}{\norm{h}} = 0,
\end{align*}
holds true, then $E$ is called \emph{Fr\'{e}chet differentiable} at $u$ and $E^\prime(u):=A$ its \emph{Fr\'{e}chet derivative}. If the Fr\'{e}chet derivative exists for all $u \in \Bspace$, $E$ is called Fr\'{e}chet differentiable.
\end{definition}

\begin{proposition}[chain rule] \label{prop:chain-rule}
Consider two functions $E \colon \Bspace \to \R$ and $\varphi \colon \R \to \R$. Let $u \in \Bspace$.  If $E$ is Fr\'{e}chet differentiable at $u$ and $\varphi$ is differentiable at $E(u)$ then the composition $\varphi \circ E \colon \Bspace \to \R$ is Fr\'{e}chet differentiable at $u$ and its Fr\'{e}chet derivative is given by
\begin{equation*}
    (\varphi \circ E)^\prime(u) = \varphi^\prime(E(u)) E^\prime(u).
\end{equation*}
\end{proposition}

\begin{proposition}[product rule]
Let $E,F \colon \Bspace \to \R$ and $u \in \Bspace$.  
If $E$ and $F$ are Fr\'{e}chet differentiable at $u$ then the product $E\cdot F:\Bspace\to\R$ is Fr\'{e}chet differentiable at $u$ and its Fr\'{e}chet derivative is given by
\begin{equation*}
    (E\cdot F)^\prime(u) = E^\prime(u)F(u) + E(u)F^\prime(u).
\end{equation*}
\end{proposition}

\begin{example}\label{ex:frechet-lin}
    Let $\Bspace$ be a Banach space and $b \in \dual{\Bspace}$ fixed. Then the Fr\'{e}chet derivative of the linear functional $u \mapsto b[u]$ is $b[u]^\prime = b$ for all $u \in \Bspace$.
\end{example}

\begin{example}\label{ex:frechet-quad}
Let $\Hspace$ be a Hilbert space and $A \in \calL(\Hspace,\Hspace)$ a bounded linear operator. Then the Fr\'{e}chet derivative of the quadratic functional $u \mapsto \frac12\sp{Au,u}$ at any $u_0 \in \Hspace$ is given by
\begin{align*}
    \left. \frac12\sp{Au,u}^\prime(\cdot) \right|_{u=u_0} = \frac12\sp{(A+A^*)u_0,\cdot}.
\end{align*}
By the Riesz representation theorem this derivative can be identified with $\frac12(A + A^*)u_0$. If $A$ is self-adjoint then $\left. \frac12\sp{Au,u}^\prime \right|_{u=u_0} = Au_0$. 
The proof is left as an exercise.
\end{example}

We frequently face the situation of a functional $E:\Hspace\to\RI$ defined on a Hilbert space whose domain $\dom(E):=\{u\in\Hspace \st E(u)\neq\pm\infty\}$ is a Banach space continuously embedded and dense in $\Hspace$.
In this case, we can nevertheless define a Fr\'echet derivative as follows.
\begin{definition}[Hilbert gradient]\label{def:H-gradient}
Let $\Hspace$ be a Hilbert space and assume that $\Bspace:=\dom(E)\contemb\Hspace$ is a continuously embedded and dense Banach space.
If $E\vert_\Bspace$ is Fr\'echet differentiable at $u\in\Bspace$ and its Fr\'echet derivative lies in $\dual{\Hspace}$, we define its $\Hspace$-gradient $\nabla_\Hspace E(u)\in\Hspace$ via
\begin{align*}
    \langle \nabla_\Hspace E(u), v\rangle_\Hspace :=
    (E\vert_\Bspace)'(u)[v],\quad\forall v\in\Hspace.
\end{align*}
\end{definition}

\begin{remark}
If $\Bspace$ is a proper subspace of $\Hspace$, it is either dense or closed.

In the first case, it is easy to check that the continuous embedding $\Bspace\contemb\Hspace$ implies the continuous embedding $\Hspace^*\contemb\Bspace^*$.
By definition the Fr\'echet derivative of $E\vert_\Bspace$ in $u\in\Bspace$ is an element of $\Bspace^*$, the condition in \cref{def:H-gradient} requires it to lie in the \emph{smaller} space $\Hspace^*$.

If $\Bspace$ is a closed subspace of $\Hspace$, it is equipped with the same norm as $\Hspace$ and we can identify $\dual{\Bspace}$ with a subspace of $\dual{\Hspace}$ in the following way: Using Hahn-Banach we can extend $\phi\in\dual{\Bspace}$ to a bounded linear functional on $\Hspace$ with the same norm as $\phi$. 
This extension is not unique if $\Bspace$ is a proper subspace of $\Hspace$.
One could use any such extension of $(E\vert_\Bspace)'(u)$ to $\dual{\Hspace}$ to define a Hilbert gradient.
However, the situation that $\Bspace$ is a closed subspace of $\Hspace$ will be of no relevance in the following.
\end{remark}

\subsubsection{Weak and weak* convergence}

\begin{definition}[Weak convergence]
A sequence $\{u_n\}_{n\in\N} \subset \Bspace$ is said to converge weakly to $u \in \Bspace$,
\begin{equation*}
    u_n \weakto u,
\end{equation*}
if for all $\phi \in \dual{\Bspace}$ one has $\phi[u_n] \to \phi[u]$.
\end{definition}
\begin{definition}[Weak* convergence]
A sequence $\{\phi_n\}_{n\in\N} \subset \dual{\Bspace}$ is said to converge weakly* to $\phi \in \dual{\Bspace}$,
\begin{equation*}
    \phi_n \wsto \phi,
\end{equation*}
if for all $u \in {\Bspace}$ one has $\phi_n[u] \to \phi[u]$.
\end{definition}
Clearly, norm convergence implies weak and weak* convergence.

\begin{definition}
A subset $K \subset \Bspace$ is called weakly (resp. strongly) sequentially compact if any sequence $\{u_n\}_{n\in\N} \subset K$ contains a weakly (resp. strongly) convergent subsequence
\begin{equation*}
    u_{n_k} \weakto u \quad \text{(resp. $u_{n_k} \to u$)}, \quad u \in K.
\end{equation*}
A subset $K \subset \dual{\Bspace}$ is called weakly* (resp. strongly) sequentially compact if any sequence $\{\phi_n\}_{n\in\N} \subset K$ contains a weakly* (resp. strongly) convergent subsequence
\begin{equation*}
    \phi_{n_k} \wsto \phi \quad \text{(resp. $\phi_{n_k} \to \phi$)}, \quad \phi \in K.
\end{equation*}

\end{definition}
\begin{theorem}[Banach-Alaoglu]\label{thm:banach_alaoglu}
Let $\Bspace$ be a Banach space.
Every closed bounded subset of $\dual{\Bspace}$ is weakly* sequentially compact.
\end{theorem}
\begin{corollary}
Every closed bounded subset of a Hilbert space ${\Bspace}$ is weakly sequentially compact.
\end{corollary}

\subsubsection{Lower semicontinuity}

The next definition gives a weaker version of continuity, which is sufficient for studying minimisers.
\begin{definition}[Lower-semicontinuity]\label{def:lsc}
Let $\tau$ be a topology on $\Bspace$ (such as strong, weak, weak*). A functional $E \colon \Bspace \to \RI$ is called $\tau$-\emph{\lsc{}}  if for any $\tau$-convergent sequence $u_n \tauto u$  one has
\begin{equation*}
    E(u) \leq \liminf_{n \to \infty} E(u_n).
\end{equation*}
\end{definition}

\begin{example}\label{ex:l1-l2-lsc}
    The functional $\norm{\cdot}_1 : \ell^2 \rightarrow \RI$ with
    \begin{align*}
        \norm{u}_1 = \begin{cases}
                \sum_{j=1}^\infty \abs{u^j} & \text{if $u \in \ell^1$}, \\ +\infty & \text{otherwise},
               \end{cases}
    \end{align*}
    where $u^j$ is the $j$-th coordinate of $u$, is \lsc{} in $\ell^2$.
    \begin{proof}
        Let $\{u_n\}_{n \in \N} \subset \ell^2$ be a convergent sequence such that $u_n \to u \in \ell^2$. Then, clearly, $u^j_n \to u^j$ for all $j \in \N$.
        We now use Fatou's lemma to obtain the assertion.
    \begin{align*}
        \norm{u}_1 = \sum_{j=1}^\infty |u^j| = \sum_{j=1}^\infty \lim_{n\to\infty} |u^j_n| \leq \liminf_{n\to\infty} \sum_{j=1}^\infty |u^j_n| = \liminf_{n\to\infty} \norm{u_n}_1.
    \end{align*}
    Note that it possible that both the right-hand side and the left-hand side of this inequality are infinite.
    \end{proof}
\end{example}

\begin{proposition}\label{ex:norm-weakly-lsc}
    Let $\Bspace$ be a Banach space with dual $\dual{\Bspace}$. Then the functional $\phi \mapsto \norm{\phi}_*$ is weakly* \lsc{} on $\dual{\Bspace}$.
    In particular, for every Hilbert space $\Hspace$ the functional $u\mapsto\norm{u}$ is weakly \lsc{.}
    \begin{proof}
       By definition it holds
        \begin{equation*}
            \norm{\phi}_* = \sup_{\substack{u \in \Bspace \\ \norm{u} \leq 1}}\phi[u].
        \end{equation*}
        Let $\phi_n \wsto \phi \in \dual{\Bspace}$. We have
        \begin{eqnarray*}
            \liminf_{n \to \infty} \norm{\phi_n}_* &=& \liminf_{n \to \infty} \sup_{\substack{u \in \Bspace \\ \norm{u} \leq 1}} \phi_n[u] \\
            &\geq& \sup_{\substack{u \in \Bspace \\ \norm{u} \leq 1}} \liminf_{n \to \infty} \phi_n[u] \\
            &=& \sup_{\substack{u \in \Bspace \\ \norm{u} \leq 1}} \phi[u] \\
            &=& \norm{\phi}_*,
        \end{eqnarray*}
        where we used the ``$\liminf \sup \geq \sup \liminf$'' inequality.
    \end{proof}
\end{proposition}

\subsubsection{Optimality conditions}

\begin{definition}
An element $u_0 \in \Bspace$ is called a local minimiser of $E \colon \Bspace \to \RI$ if there exists $\eps >0$ such that
\begin{equation*}
    E(u_0) \leq E(u) \quad \text{for all $u$ s.t. $\norm{u-u_0} \leq \eps$}.
\end{equation*}
The local minimum is called strict if the inequality is strict for all $u \neq u_0$.
\end{definition}

The next result is known as Fermat's rule.

\begin{theorem}[Fermat's rule]\label{thm:FOC}
Let $\Hspace$ be a Hilbert space and $u_0 \in \Hspace$ a local minimiser of $E \colon \Hspace \to \RI$. 
For any $v\in\Hspace$ for which the function $\eps\mapsto E(u_0+\eps v)$ is differentiable the first variation $\delta E(u_0;v)$ exists and satisfies $\delta E(u_0;v)=0$.
\end{theorem}
\begin{proof}
Since the real-valued function $\eps\mapsto E(u_0+\eps v)$ has a local minimum at $\eps=0$ and is differentiable, it holds
\begin{align*}
    0 = \frac{\d}{\d\eps}\Big\vert_{\eps=0}E(u_0+\eps v) = \lim_{\eps\to 0}\frac{E(u_0+\eps v)-E(u_0)}{\eps} = \delta E(u_0;v).
\end{align*}
\end{proof}

\section{Rayleigh quotient and  Euler--Lagrange equation}

For the rest of this section $\Omega\subset\R^d$ is assumed to be a bounded domain with Lipschitz boundary.
In particular, this allows us to equip the space $H^1_0(\Omega)$ with the norm $\norm{u}_{H^1_0(\Omega)}:=\norm{\grad u}_{L^2(\Omega)}$.
At some places, we add the extra assumption that the boundary is even smooth.

We consider the following functional referred to as the \emph{Dirichlet energy}
\begin{equation*}
    u \mapsto \frac{1}{2} \int_\Omega\abs{\nabla u}^2\d x, \quad u \in H^1_0(\Omega).
\end{equation*}
This functional can be extended to $L^2(\Omega)$ by letting it be infinite on $L^2(\Omega) \setminus H^1_0(\Omega)$:
\begin{equation*}
    u \mapsto \calD_2(u) :=
    \begin{cases}
        \frac{1}{2} \int_\Omega\abs{\nabla u}^2\d x, \quad & u \in H^1_0(\Omega), \\
        +\infty, \quad & u \in L^2(\Omega) \setminus H^1_0(\Omega).
    \end{cases}
\end{equation*}
It can be easily seen that the functional $\calD_2$ is absolutely $2$-homogeneous, i.e. $\calD_2 (\alpha u) = \abs{\alpha}^2 \calD_2(u)$ for all $u \in L^2(\Omega)$ and $\alpha \in \R$.

\begin{proposition}\label{prop:lsc_dirichlet}
The functional $\calD_2 \colon L^2(\Omega) \to \RI$ is proper and \lsc{} on $L^2(\Omega)$.
\end{proposition}
\begin{proof}
Since $0\in H^1_0(\Omega)$ and $\calD_2(0) = 0 < \infty$, the functional is proper. 
Let $u_n \to u_\infty$ in $L^2(\Omega)$ and consider $a \defeq \liminf \norm{\grad u_n}_{L^2}$. Without loss of generality, we may assume $a < + \infty$, otherwise the statement is trivial. Hence, there is a subsequence $u_{n_k}$ such that $a = \lim_{k \to \infty} \norm{\grad u_{n_k}}_{L^2}$ and hence $\norm{u_{n_k}}_{H^1_0} = \norm{\grad u_{n_k}}_{L^2} \leq C$ for some $C>0$. Therefore, $u_{n_k}$ is bounded in $H^1_0(\Omega)$ and therefore contains weakly convergent subsequence (which we don't relabel to avoid a triple index). By the uniqueness of the limit, we have
\begin{equation*}
    u_{n_k} \weakto u_\infty \quad \text{weakly in $H^1_0(\Omega)$}.
\end{equation*}
Since Hilbert norms are weakly \lsc{,} see \cref{ex:norm-weakly-lsc}, we have
\begin{eqnarray*}
    \norm{\grad u_\infty}_{L^2} = \norm{u_\infty}_{H^1_0} &\leq&  \liminf_{k\to\infty} \norm{u_{n_k}}_{H^1_0} =  \liminf_{k\to\infty} \norm{\grad u_{n_k}}_{L^2} \\
    &=& \lim_{k \to \infty} \norm{\grad u_{n_k}}_{L^2} = a = \liminf_{n\to\infty} \norm{\grad u_n}_{L^2}.
\end{eqnarray*}

Since the function $t \mapsto t^2$ is continuous and nondecreasing on $\R_+$, this implies the claim.
\end{proof}

Consider the following \emph{Rayleigh quotient}
\begin{equation}\label{eq:Rayleigh-quotient-lin}
    u \mapsto \calR_2(u):= \frac{2\calD_2(u)}{\norm{u}_{L^2}^2}, \quad u \in L^2(\Omega).
\end{equation}
This functional is, clearly, $0$-homogeneous, i.e. its value does not change if the argument is multiplied by any constant $\alpha \neq 0$. Furthermore, for functions $u\in H^1_0(\Omega)$ the Dirichlet energy is finite and the Rayleigh quotient equals
\begin{align*}
    \calR_2(u) = \frac{\norm{\grad u}_{L^2}^2}{\norm{u}_{L^2}^2}.
\end{align*}
We will study stationary points of this Rayleigh quotient, i.e., points whose first variation is zero in sufficiently many directions. 
For this, we first compute its first variation.

\begin{proposition}\label{prop:first_var_rayleigh}
Let $u\in H^1_0(\Omega)\setminus\{0\}$ and $v\in C^\infty_c(\Omega)$.
The first variation of the functional~\eqref{eq:Rayleigh-quotient-lin} in $u$ in the direction of $v$ is given by
\begin{equation*}
    \delta\calR_2(u;v) = \frac{1}{\norm{u}_{L^2}^2} \left(\int_\Omega\grad u\cdot\grad v\d x - \frac{\norm{\grad u}_{L^2}^2}{\norm{u}_{L^2}^2} \int_\Omega u\, v\d x \right).
\end{equation*}
\end{proposition}
\begin{proof}
Using \cref{prop:chain_rule_first_var,prop:prod_rule_first_var} it sufficies to compute the first variations of the functionals $u\mapsto\calD_2(u)$ and $u\mapsto\norm{u}_{L^2}^2$:
\begin{align*}
    \delta(2\calD_2)(u;v) 
    &= 
    \frac{\d}{\d\eps}\Big\vert_{\eps=0} 
    2\calD_2(u+\eps v)
    =
    \frac{\d}{\d\eps}\Big\vert_{\eps=0} 
    {\int_\Omega\abs{\nabla u + \eps\nabla v}^2\d x}
    \\
    &=
    2 \int_\Omega\nabla u \cdot \nabla v \d x. \\
    \delta\norm{\cdot}_{L^2}^2(u;v) 
    &= 
    \frac{\d}{\d\eps}\Big\vert_{\eps=0} 
    \int_\Omega\abs{u + \eps v}^2\d x
    =
    2 \int_\Omega u \, v \d x.
\end{align*}
\end{proof}

A direct consequence of this result is
\begin{theorem}\label{thm:opt-cond-lin}
Let $\Omega\subset\R^d$ denote an open domain with smooth boundary.
An element $u \in H^1_0(\Omega)\setminus\{0\}$ satisfies $\delta\calR_2(u;v)=0$ for all $v\in C^\infty_c(\Omega)$ if and only if $u\in H^2(\Omega)\cap H^1_0(\Omega)\setminus\{0\}$ and satisfies the following eigenvalue problem
\begin{equation}\label{eq:euler-lagrange-lin}
    -\laplace u = \lambda u \quad \text{in }\Omega,
\end{equation}
where $\lambda \defeq \frac{\norm{\grad u}_{L^2}}{\norm{u}_{L^2}}$. 
\end{theorem}
\begin{proof}
By \cref{prop:first_var_rayleigh} we know that $\delta\calR_2(u;v)=0$ if and only if
\begin{align}\label{eq:weak_ev_problem}
    \int_\Omega \nabla u \cdot \nabla v \dx - \lambda \int_\Omega u\, v\d x = 0.
\end{align}
Under the hypothesis that $u\in H^2(\Omega)$ and hence $\laplace u\in L^2(\Omega)$, this is equivalent to
\begin{align*}
\int_\Omega \big(-\laplace u - \lambda u\big) v \dx = 0.
\end{align*}
By the density of $C^\infty_c(\Omega)$ in $L^2(\Omega)$ this is equivalent to $-\laplace u = \lambda u$.

It remains to shows that $u$ solving \labelcref{eq:weak_ev_problem} implies that $u\in H^2(\Omega)$. 
This follows from standard elliptic theory since $u$ is a weak solution to the Poisson equation with right hand side $\lambda u\in L^2(\Omega)$ and $\partial\Omega$ is smooth, see for instance \cite{savare1998regularity}.
\end{proof}
\begin{remark}\label{rem:regularity}
For smooth domains, using bootstrapping one gets that any weak solution $u\in H^1_0(\Omega)$ of $-\laplace u = \lambda u$ satisfies $u\in H^k\cap H^1_0(\Omega)$ for all $k\in\N$.  
Using Sobolev embeddings yields classical smoothness up to the boundary and in fact $u$ is even analytic.
\end{remark}

\cref{thm:opt-cond-lin} shows that \ef{s} of the negative Laplace operator are exactly the critical points of the Rayleigh quotient~\eqref{eq:Rayleigh-quotient-lin}. We will use this observation to generalise the concept of an \ef{} in \cref{ch:nonlin-ef}.

Now consider minimising the Rayleigh quotient~\eqref{eq:Rayleigh-quotient-lin}
\begin{equation}\label{eq:Rayleigh-opt-prob-lin}
    \inf_{u \in L^2(\Omega)} \frac{2\calD_2(u)}{\norm{u}_{L^2}^2}
    = 
    \inf_{u \in H^1_0(\Omega)} \frac{\norm{\grad u}_{L^2}^2}{\norm{u}_{L^2}^2}.
\end{equation}
This problem can be written equivalently as an optimisation problem on the unit sphere,
\begin{equation}\label{eq:Rayleigh-opt-prob-lin-constrained}
    \inf_{\substack{u\in H^1_0(\Omega) \\ \norm{u}_{L^2}^2=1}} \norm{\grad u}_{L^2},
\end{equation}
but not on the unit ball $\{u \colon \norm{u}\leq1\}$.

\begin{theorem}
    The optimisation problem~\eqref{eq:Rayleigh-opt-prob-lin} admits a minimiser.
\end{theorem}
\begin{proof}
This follows from the direct method, using \cref{cor:rellich,prop:lsc_dirichlet}. We will show it later in \cref{ch:nonlin-ef} in more generality.
\end{proof}

\begin{definition}[ground states]
    Any global minimiser of~\eqref{eq:Rayleigh-opt-prob-lin} is called a \emph{ground state} of the Dirichlet energy.
\end{definition}

Using standard elliptic regularity theory for PDEs one can prove positivity of ground states.
\begin{proposition}\label{prop:positive_gs}
Let $\Omega\subset\R^d$ denote an open domain with smooth boundary.
Any ground state of the Dirichlet energy is either zero, or positive everywhere, or negative everywhere on~$\Omega$.
\end{proposition}
\begin{proof}
If $u$ is a ground state, so is $\abs{u}$. This is because $\calD_2(\abs{u})\leq\calD_2(u)$ and $\norm{\abs{u}}_{L^2}=\norm{u}_{L^2}$.
By standard elliptic regularity theory any solution to $\lambda u = -\laplace u$ in $H^1_0(\Omega)$ has a representative which is continuous on $\overline{\Omega}$.
Furthermore, this representative satisfies the Harnack $\max_{\closure{B}(x,r)} u \leq C \min_{\closure{B}(x,r)} u$ (see \cite{trudinger1967harnack}) for any $x\in\Omega$ such that $B(x,2r)\subset\Omega$.
Applying this to the function $\abs{u}$---which is a solution thanks to \cref{thm:opt-cond-lin}---we see get that $\abs{u}>0$ has to hold in $\Omega$.  
By continuity $u>0$ (or $u<0$).
\end{proof}

\begin{theorem}\label{thm:uniqueness_gs}
Let $\Omega\subset\R^d$ denote an open domain with smooth boundary.
The minimiser of~\eqref{eq:Rayleigh-opt-prob-lin} is unique up to scalar multiplication.
\end{theorem}
\begin{proof}
The proof is a special case of \cite[Lemma 3.2]{kawohl2006positive}:

Let $u_1$ and $u_2$ be two different ground states and assume that both are positive and in $C^1(\Omega)$.
Then we can define $v:=(u_1^2+u_2^2)^\frac{1}{2}$ and, since $u_1,u_2$ are ground states, it holds
\begin{align}\label{eq:competitor_gs}
    \lambda := \frac{\int_\Omega\abs{\grad u_1}^2\d x}{\int_\Omega\abs{u_1}^2\d x} = \frac{\int_\Omega\abs{\grad u_2}^2\d x}{\int_\Omega\abs{u_2}^2\d x} 
    \leq 
    \frac{\int_\Omega\abs{\grad v}^2\d x}{\int_\Omega\abs{v}^2\d x}.
\end{align}
We can write
\begin{align*}
    \grad v = \frac{1}{v}(u_1\grad u_1 + u_2\grad u_2)
    =
    v\left(\frac{u_1^2\grad\log u_1 + u_2^2\grad\log u_2}{u_1^2+u_2^2}\right)
\end{align*}
and the brackets are a convex combination of $\grad\log u_1$ and $\grad\log u_2$.
Strict convexity of the function $t\mapsto t^2$ then implies
\begin{align*}
    \abs{\grad v}^2 
    &\leq 
    v^2
    \left(\frac{u_1^2\abs{\grad\log u_1}^2 + u_2^2\abs{\grad\log u_2}^2}{u_1^2+u_2^2}\right)
    \\
    &=
    v^2
    \left(\frac{\abs{\grad u_1}^2 + \abs{\grad u_2}^2}{u_1^2+u_2^2}\right)
    =
    \abs{\grad u_1}^2 + \abs{\grad u_2}^2,
\end{align*}
and the inequality is strict if $\grad\log u_1\neq\grad\log u_2$.
Integrating this and using \labelcref{eq:competitor_gs} along with $v^2 = u_1^2 + u_2^2$ then gives the contradiction:
\begin{align*}
    \lambda 
    &\leq 
    \frac{\int_\Omega\abs{\grad v}^2\d x}{\int_\Omega\abs{v}^2\d x}
    <
    \frac{\int_\Omega\abs{\grad u_1}^2\d x+\int_\Omega\abs{\grad u_2}^2\d x}{\int_\Omega\abs{u_1}^2\d x+\int_\Omega\abs{u_2}^2\d x}
    \\
    &=
    \frac{\int_\Omega\abs{u_1}^2\d x}{\int_\Omega\abs{u_1}^2\d x+\int_\Omega\abs{u_2}^2\d x} \lambda 
    +
    \frac{\int_\Omega\abs{u_2}^2\d x}{\int_\Omega\abs{u_1}^2\d x+\int_\Omega\abs{u_2}^2\d x} \lambda = \lambda.
\end{align*}
Hence, almost everywhere we must have $\grad\log u_1 = \grad\log u_2$ which implies $u_1=Cu_2$ for $C>0$.
\end{proof}
\begin{remark}
In fact assuming a smooth domain is not necessary for \cref{thm:uniqueness_gs} since one can exhaust any bounded domain with smooth sets and the minimal values of the Rayleigh quotients converge \cite{kawohl2006positive}.
\end{remark}

\section{Gradient flows and asymptotic profiles}
\label{sec:gradflow_profiles_linear}

In this section we will consider the gradient flow of the Dirichlet energy $\calD_2$.
Abstractly, the gradient flow of a function $E:\Hspace\to\RI$ on a Hilbert space is defined as solution to the differential equation
\begin{align}\label{eq:gateaux_gf}
    u^\prime(t) = -\grad E(u(t)),\quad u(0) = f,
\end{align}
where $\grad E$ denotes a suitable notion of gradient of the functional $E$, and $f\in\Hspace$ is an initial datum.
Depending on the regularity of $E$ different notions of gradients are appropriate, e.g., Fr\'echet-gradients, Gateaux-gradients, subgradient, or the Hilbert space gradient.

Since the Dirichlet energy is nowhere continuous on $L^2(\Omega)$, one either has to work with its $L^2(\Omega)$-gradient as defined in \cref{def:H-gradient} or with its subgradient.
We will get to general theory of subgradient flows of arbitrary convex (see \cref{def:convexity} in \cref{ch:nonlin-ef}) functionals $E:\Hspace\to\RI$ later in full generality but for now we stick to the $L^2(\Omega)$-gradient of the Dirichlet energy.
For this we first compute the Fr\'echet derivative of the restriction of $\calD_2$ onto $H^1_0(\Omega)$.
\begin{proposition}\label{prop:frechet_dirichlet}
Let $E:H^1_0(\Omega)\to\R$ be defined via $E(u):=\calD_2(u) = \frac{1}{2}\norm{\grad u}_{L^2}^2$.
Then $E$ is Fr\'echet differentiable in all $u\in H^1_0(\Omega)$ and its Fr\'echet derivative is given by
\begin{align*}
    E'(u)[v] := \int_\Omega \grad u \cdot \grad v \d x,\quad\forall v\in H^1_0(\Omega).
\end{align*}
\end{proposition}
\begin{proof}
We compute
\begin{align*}
    E(u + v) - E(u) 
    &=
    \frac{1}{2}\int_\Omega\abs{\grad u + \grad v}^2\d x - \frac{1}{2}\int_\Omega \abs{\grad u}^2\d x 
    \\
    &=
    \frac{1}{2}\int_\Omega \abs{\grad v}^2\d x
    + \int_\Omega \grad u \cdot \grad v \d x
\end{align*}
and hence
\begin{align*}
    \frac{\abs{E(u+v) - E(u) - \int_\Omega \grad u \cdot \grad v \d x}}{\norm{v}_{H^1_0}}
    =
    \frac{1}{2} \frac{\norm{v}_{H^1_0}^2}{\norm{v}_{H^1_0}} =  \frac{1}{2} \norm{v}_{H^1_0} \to 0,\quad\text{as }\norm{v}_{H^1_0} \to 0.
\end{align*}
It remains to show that the map $v\mapsto E'(u)[v]:=\int_\Omega\grad u\cdot\grad v\d x$ is linear and bounded. 
The linearity is trivial whereas the boundedness follows from Hölder's inequality:
\begin{align*}
    \abs{E'(u)[v]}
    \leq
    \int_\Omega \abs{\grad u}\abs{\grad v}\d x \leq \underbrace{\norm{\grad u}_{L^2}}_{=:C}\norm{\grad v}_{L^2} = C\norm{v}_{H^1_0}.
\end{align*}
This concludes the proof.
\end{proof}

Having the Fr\'echet derivative of the Dirichlet energy on its domain at hand, we can now compute the $L^2(\Omega)$-gradient from \cref{def:H-gradient}.

\begin{proposition}
The $L^2(\Omega)$-gradient of the Dirichlet energy $\calD_2:L^2(\Omega)\to\RI$ in $u\in H^2(\Omega)\cap H^1_0(\Omega)$ is given by
\begin{align*}
    \grad_{L^2}\calD_2(u) = -\laplace u.
\end{align*}
\end{proposition}
\begin{proof}
By \cref{def:H-gradient,prop:frechet_dirichlet} it holds for all $v \in H^1_0(\Omega)$ that
\begin{align*}
    \langle \grad_{L^2}\calD_2(u), v\rangle_{L^2} 
    = 
    (\mathcal{D}_2\vert_{H^1_0(\Omega)})'(u)[v]
    =
    \int_\Omega \nabla u \cdot \nabla v \d x
    =
    -\int_\Omega \laplace u \, v \d x.
\end{align*}
Since $H^1_0(\Omega)$ is dense and continuously embedded in $L^2(\Omega)$, we obtain the assertion.
\end{proof}

Hence, the $L^2(\Omega)$-gradient flow \eqref{eq:gateaux_gf} of the Dirichlet energy is the heat equation
\begin{align}\label{eq:heat_eq}
    u^\prime(t) = \laplace u(t)
\end{align}
with homogeneous boundary conditions, incorporated through $u(t)\in H^1_0(\Omega)$ for all $t>0$.

Using the spectral theorem for the Laplace operator \cref{thm:laplace-efs}, we can solve the heat equation analytically.
For this we express the initial datum $f\in L^2(\Omega)$ in terms of the orthonormal eigenbasis $\{u_i\}_{i\in\N}$ as
\begin{align*}
    f = \sum_{i=1}^\infty c_i u_i,
\end{align*}
where the coefficients are given by $c_i := \langle f,u_i\rangle_{L^2}$ for $i\in\N$.
We assume that the \ef{s} are ordered such that the \ev{s} satisfy $0<\lambda_1\leq\lambda_2\leq\dots$.
Using the ansatz
\begin{align*}
    u(t) = \sum_{i=1}^\infty c_i(t)u_i
\end{align*}
and plugging it into the heat equation, a formal computation yields
\begin{align*}
    u^\prime(t) - \laplace u(t) 
    =
    \sum_{i=1}^\infty
    c_i^\prime(t)u_i - \sum_{i=1}^\infty c_i(t)\laplace u_i 
    =
    \sum_{i=1}^\infty
    \big(c_i^\prime(t) + \lambda_i c_i(t)\big)u_i.
\end{align*}
Hence, the coefficient functions should satisfy 
\begin{align*}
    \begin{cases}
        c_i^\prime(t)=-\lambda_ic_i(t), \quad t>0,\\
        c_i(0) = c_i,
    \end{cases}
    \quad\forall i\in\N.
\end{align*}
The unique solution to these ODEs is $c_i(t)=c_i\exp(-\lambda_i t)$ and hence 
\begin{align}\label{eq:sol_heat_eq}
    u(t) = \sum_{i=1}^\infty c_i\exp(-\lambda_i t)u_i
\end{align}
solves the heat equation \eqref{eq:heat_eq}.
From this explicit representation one can derive the following properties of the solution.
\begin{theorem}[Asymptotic behavior of the heat equation]
    Let $f\in L^2(\Omega)$ and $i_0 = \min\{i\in\N\st \langle f,u_i\rangle_{L^2}\neq 0\}$.
    Then the solution of the heat equation \eqref{eq:heat_eq} with initial datum $f\in L^2(\Omega)$ satisfies:
    \begin{itemize}
        \item $u(t) \to 0$ strongly in $L^2(\Omega)$ as $t\to\infty$,
        \item $\langle f,u_{i_0}\rangle_{L^2}\exp(-\lambda_{i_0} t)\leq \norm{u(t)}_{L^2}\leq \norm{f}_{L^2}\exp(-\lambda_{i_0} t)$ for all $t\geq 0$,
        \item
        Furthermore, if $u_{i_0}$ is simple---meaning that the solutions of $\lambda_{i_0} u = -\laplace u$ form a one-dimensional subspace of $L^2(\Omega)$---then it holds
        $\frac{u(t)}{\norm{u(t)}_{L^2}}\to u_{i_0}$.
    \end{itemize}
\end{theorem}
\begin{remark}
According to \cref{thm:uniqueness_gs} the simplicity is in particular true for the first eigenfunction $u_1$ of the Laplacian, which we also showed to have a sign.
Hence, if $\sp{f,u_1}_{L^2}\neq 0$ (e.g., satisfied if $f$ has a sign) we get that $\frac{u(t)}{\norm{u(t)}_{L^2}}$ converges to the first eigenfunction $u_1$.
\end{remark}
\begin{remark}
If the simplicity is dropped, then still $\frac{u(t)}{\norm{u(t)}_{L^2}}$ converges to an eigenfunction, which then is a linear combination of different eigenfunctions from the same eigenspace.
The proof is left as an exercise.
\end{remark}
\begin{proof}
The second statement obviously implies the first one. 
Letting  $c_i := \langle f, u_i\rangle_{L^2}$ it holds
\begin{align*}
    \norm{f}_{L^2}^2 
    = 
    \langle f, f \rangle_{L^2} 
    = 
    \sum_{i,j=1}^\infty c_i c_j \underbrace{\langle u_i, u_j \rangle_{L^2}}_{=\delta_{ij}}
    =
    \sum_{i=1}^\infty c_i^2
    =
    \sum_{i=i_0}^\infty c_i^2.
\end{align*}
Similarly, we compute
\begin{align*}
    \norm{u(t)}_{L^2}^2 
    &= 
    \langle u(t), u(t)\rangle_{L^2} 
    =
    \sum_{i,j=1}^\infty c_i c_j\exp(-\lambda_i t)\exp(-\lambda_j t) \underbrace{\langle u_i, u_j\rangle_{L^2}}_{=\delta_{ij}}
    \\
    &=
    \sum_{i=1}^\infty c_i^2\exp(-2\lambda_i t)
    =
    \sum_{i=i_0}^\infty c_i^2\exp(-2\lambda_i t)
    \\
    &=
    \exp(-2\lambda_{i_0} t)\sum_{i=i_0}^\infty c_i^2\exp(-2(\lambda_i-\lambda_{i_0}) t)
    \\
    &=
    \exp(-2\lambda_{i_0} t) 
    \left(
    c_{i_0}^2 + 
    \sum_{i=i_0+1}^\infty c_i^2\exp(-2(\lambda_i-\lambda_{i_0}) t)
    \right).
\end{align*}
Using the two estimates $0 \leq \exp(-2(\lambda_i - \lambda_{i_0})t) \leq 1$ 
we get
\begin{align*}
    c_{i_0}^2 \exp(-2\lambda_{i_0}t) \leq \norm{u(t)}^2_{L^2} \leq \norm{f}_{L^2}^2\exp(-2\lambda_{i_0}t)
\end{align*}
which implies second statement after taking the square root.

Moreover, since $u_{i_0}$ is simple, we know that $\lambda_j > \lambda_{i_0}$ for all $j\geq i_0 + 1$ and one gets
\begin{align*}
    \frac{\norm{u(t)}_{L^2}^2}{\exp(-2\lambda_{i_0}t)} 
    &= 
    c_{i_0}^2 + \sum_{i=i_0+1}^\infty c_i^2\exp(-2(\lambda_i-\lambda_{i_0}) t)
    \\
    &\leq
    c_{i_0}^2 + \underbrace{\exp(-2(\lambda_{i_0+1}-\lambda_{i_0}) t)}_{\to 0 \text{ as }t\to\infty} \underbrace{\sum_{i=i_0+1}^\infty c_i^2}_{\leq\norm{f}_{L^2}^2}
    \to c_{i_0}^2 \quad\text{as }t\to\infty.
\end{align*}
This implies the third statement:
\begin{align*}
    \norm{\frac{u(t)}{\norm{u(t)}_{L^2}}-u_{i_0}}_{L^2}^2
    &=
    1 - 2 \frac{\langle u(t), u_{i_0}\rangle_{L^2}}{\norm{u(t)}_{L^2}}
    + 1
    =
    2 - 2
    \sum_{i=1}^\infty
    \frac{c_i\exp(-\lambda_i t)}{\norm{u(t)}_{L^2}}
    \langle u_i, u_{i_0} \rangle_{L^2}
    \\
    &=
    2 - 2 \frac{c_{i_0}\exp(-\lambda_{i_0}t)}{\norm{u(t)}_{L^2}}
    \to 0\quad\text{as }t\to\infty.
\end{align*}
\end{proof}

\chapter{Nonlinear Eigenvalue Problems in Hilbert Spaces}
\label{ch:nonlin-ef}

In this chapter we extend the results of \cref{ch:lin-ef} to the setting where the Dirichlet energy $\calD_2$ is replaced by a non-quadratic energy functional. 
Most importantly, we allow the Fr\'echet derivative (or more general the subdifferential) of the considered functionals to be a nonlinear, and even multi-valued, operator.
We start by recalling some necessary background material.

\section{Background on Convex Analysis}\label{sec:nonlin-background}

\begin{definition}[convexity]\label{def:convexity}
A subset $\calC \subset \Bspace$ is called \emph{convex} if for all $x,x' \in \calC$ and $\alpha \in (0,1)$ 
\begin{equation*}
    \alpha x + (1-\alpha)x' \in \calC.
\end{equation*}
A function $E \colon \Bspace \to \RI$ is called convex if for all $x,x' \in D$ and $\alpha \in (0,1)$ 
\begin{equation*}
    E(\alpha x + (1-\alpha)x') \leq \alpha E(x) + (1-\alpha) E(x').
\end{equation*}
It is called strictly convex if the inequality is strict for $x\neq x'$ and $\theta$-strongly convex if $E(\cdot) - \frac{\theta}{2} \norm{\cdot}^2$ is convex.
\end{definition}
\noindent The following implications hold
\begin{equation*}
    \text{strongly convex} \implies \text{strictly convex} \implies \text{convex}.
\end{equation*}

\begin{proposition}
A subset $\calC \subset \Bspace$ is convex if and only if its characteristic function is convex.
\end{proposition}

One calls $\FunA$ a concave function if $-\FunA$ is convex and analogously one can define strictly and $\theta$-strongly concave functions.

\subsubsection{Lower semicontinuity}

\begin{proposition}[{\cite[Cor. 2.2, p. 11]{ekeland-temam:1976}} or {\cite[p. 149]{bauschke-combettes:2011}}] \label{prop:LSC_convex}
A convex function $\FunA:\Bspace\to\RI$ is \lsc{} if and only if it is weakly \lsc{.}
\end{proposition}

\subsubsection{Convex conjugate}

In convex optimisation problems (i.e. those involving convex functions) the concept of \emph{convex conjugates} plays an important role.
\begin{definition}
Let $\FunFromTo{\FunA}{\Bspace}{\RI}$ be a functional. The functional $\FunFromTo{\conj{\FunA}}{\dual{\Bspace}}{\RI}$,
\begin{equation*}
\conj{\FunA}(\subgradA) = \sup_{\imageA \in \Bspace} [\scalprod{\subgradA}{\imageA} - \FunA(\imageA)],
\end{equation*}
\noindent is called the convex conjugate of $\FunA$.
\end{definition}

The convex conjugate is a lower semicontinuous and convex function.

\begin{theorem}[{\cite[Prop. 4.1]{ekeland-temam:1976}}] \label{thm:biconj}
For any functional $\FunFromTo{\FunA}{\Bspace}{\RI}$ the following inequality holds:
\begin{equation*}
\biconj{\FunA} \defeq \conj{(\conj{\FunA})}\vert_\Bspace \leq \FunA.
\end{equation*}
If $\FunA$ is proper, lower-semicontinuous (see \cref{def:lsc}) and convex, then
\begin{equation*}
\biconj{\FunA} = \FunA.
\end{equation*}
\end{theorem}

\subsubsection{Subgradients}
For convex functions one can generalise the concept of a derivative so that it also makes sense for non-differentiable functions.

\begin{definition}
A functional $\FunFromTo{\FunA}{\Bspace}{\RI}$ is called \emph{subdifferentiable} at $\imageA \in \Bspace$, if there exists an element $\subgradA \in \dual{\Bspace}$ such that
\begin{align*}
\FunA(\imageB) \geq \FunA(\imageA) + \innerp{\subgradA}{\imageB - \imageA}
\end{align*}
holds, for all $\imageB \in \Bspace$. Furthermore, we call $\subgradA$ a \emph{subgradient} at position $\imageA$. The collection of all subgradients at position $\imageA$, i.e.
\begin{align*}
\partial \FunA(\imageA) := \left\{ \subgradA \in \dual{\Bspace} \st \FunA(\imageB) \geq \FunA(\imageA) + \innerp{\subgradA}{\imageB - \imageA} \, \text{,} \, \forall \imageB \in \Bspace \right\} \, ,
\end{align*}
is called \emph{subdifferential} of $\FunA$ at $\imageA$.
We also define the domain of the subdifferential as
\begin{align*}
    \dom(\partial\FunA) = \{\imageA \in \dual{\Bspace} \st \partial\FunA(\imageA)\neq\emptyset\}.
\end{align*}
\end{definition}

It is clear that if a convex functional $\FunFromTo{\FunA}{\Bspace}{\RI}$ is proper, i.e. $\dom(\FunA) \neq \emptyset$, then for all $\imageA \not\in \dom(\FunA)$ the subdifferential is empty. 
In other words it holds $\dom(\partial\FunA)\subset\dom(\FunA)$.
A sufficient (but not necessary) condition for $E$ to have a subgradient at $u \in \dom(E)$ is given by

\begin{proposition}[{\cite[Prop. 5.2]{ekeland-temam:1976}}]
Let $\FunFromTo{\FunA}{\Bspace}{\RI}$ be a convex functional and $u \in \dom(E)$ such that $E$ is continuous at $u$. Then $\subdiff E(u) \neq \emptyset$.
\end{proposition}

\begin{theorem}[{\cite[Thm. 7.13]{aliprantis-border:2006}}]\label{thm:subdiff_closed}
Let $\FunFromTo{\FunA}{\Bspace}{\RI}$ be a proper convex function and $\imageA \in \dom(\FunA)$. Then $\partial \FunA (\imageA)$ is a weakly-* closed and convex subset of $\dual{\Bspace}$.
\end{theorem}

\begin{figure}
\centering
\includegraphics{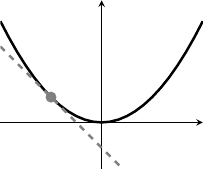}

\qquad
\includegraphics{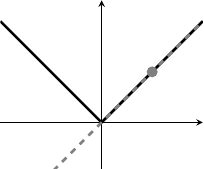}
\qquad
\includegraphics{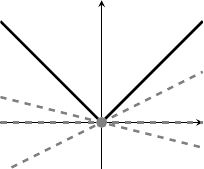}
\caption{Visualisation of the subdifferential. Linear approximations of the functional have to lie completely underneath the function. For points where the function is not differentiable there may be more than one such approximation.} \label{FIG:VARREG:SUBDIFFERENTIAL}
\end{figure}

For Fr\'{e}chet differentiable functions the subdifferential consists of just one element -- the Fr\'{e}chet derivative. For non-differentiable functionals the subdifferential is multivalued; we want to consider the subdifferential of the absolute value function as an illustrative example.
\begin{example}\label{exm:absval}
Let $\FunFromTo{\FunA}{\R}{\R}$ be the absolute value function $\FunA(\imageA) = |\imageA|$. Then, the subdifferential of $E$ at $\imageA$ is given by
\begin{align*}
\partial E(\imageA) = \begin{cases} \{1\} & \text{for $\imageA > 0$}\\ [-1, 1] & \text{for $\imageA = 0$}\\ \{-1\} & \text{for $\imageA < 0$} \end{cases} \, ,
\end{align*}
which is left as an exercise. A visual explanation is given in Figure \ref{FIG:VARREG:SUBDIFFERENTIAL}.

\end{example}

The subdifferential of a sum of two functions can be characterised as follows.
\begin{theorem}[{\cite[Prop. 5.6]{ekeland-temam:1976}}]\label{thm:subdiff_sum}
Let $\FunFromTo{\FunA}{\Bspace}{\RI}$ and $\FunFromTo{\FunB}{\Bspace}{\RI}$ be proper l.s.c. convex functions.
\begin{itemize}
    \item If $\dom(\FunA)\cap\dom(\FunA)\neq\emptyset$ it holds 
    \begin{align*}
    \partial(\FunA+\FunB) \supset \partial\FunA + \partial\FunB.
    \end{align*}
    \item If $\exists \imageA \in \dom(\FunA) \cap \dom(\FunB)$ such that $\FunA$ is continuous at $\imageA$, then
    \begin{align*}
        \partial(\FunA+\FunB) = \partial\FunA + \partial\FunB.
    \end{align*}
\end{itemize}
\end{theorem}

Using the subdifferential, one can characterise minimisers of convex functionals.
\begin{theorem}[Fermat's rule for subdifferentials]
An element $\imageA \in \Bspace$ is a minimiser of the functional $\FunFromTo{\FunA}{\Bspace}{\RI}$ if and only if $0 \in \partial \FunA(\imageA)$.
\begin{proof}
By definition, $0 \in \partial \FunA(\imageA)$ if and only if for all $\imageB \in \Bspace$ it holds
\begin{align*}
    \FunA(\imageB) \geq \FunA(\imageA) + \innerp{0}{\imageB-\imageA} = \FunA(\imageA) \, ,
\end{align*}
which is by definition the case if and only if $\imageA$ is a minimiser of $\FunA$.
\end{proof}
\end{theorem}

The next result connects subgradients and convex conjugates
\begin{theorem}[{\cite[Prop. 5.1]{ekeland-temam:1976}}] \label{thm:fenchel-young}
Let $E \colon \Bspace \to \RI$ be a convex function and $E^* \colon \Bspace^* \to \RI$ its convex conjugate. Then  $\subgradA \in \subdiff E(\imageA)$ if and only if
\begin{equation*}
 E(\imageA) + E^*(\subgradA) = \sp{\subgradA,\imageA}.
\end{equation*}
\end{theorem}
\begin{proof}
Left as an exercise.
\end{proof}

\subsubsection{Convex \texorpdfstring{$p$}{p}-homogeneous functionals}

\begin{definition}
A functional $\FunFromTo{\FunA}{\Bspace}{\RI}$ is called $p$-homogeneous, $p > 0$, if for all $u \in \Hspace$ and $t \neq 0$ one has
\begin{equation}\label{eq:p-hom}
    \FunA(tu) = \abs{t}^p \FunA(\imageA),
    \qquad
    \FunA(0) = 0.
\end{equation}
\end{definition}

Such functionals have some useful properties, we start with the simplest one:
\begin{proposition}\label{prop:positivity}
Let $\FunA:\Bspace\to\RI$ be absolutely $p$-homogeneous for $p>0$ and convex.
Then $\FunA \geq 0$.
\end{proposition}
\begin{proof}
Using the assumptions we can compute
\begin{align*}
    0 = \FunA(0) = \FunA(0\cdot u) = \FunA\left(\frac{1}{2}u + \frac{1}{2}(-u)\right) \leq \frac{1}{2}\FunA(u) + \frac{1}{2}\FunA(-u) = \FunA(u),\quad\forall u\in\Bspace.
\end{align*}
\end{proof}
Next, we show that non-zero convex functionals can only be absolutely $p$-homogeneous for $p\geq 1$.
\begin{proposition}
Let $\FunA:\Bspace\to\RI$ be convex and absolutely $p$-homogeneous for $p>0$.
Then either $\FunA\equiv 0$ on $\dom(\FunA)$ or $p\geq 1$.
\end{proposition}
\begin{proof}
Let us assume that $0<p<1$ and $\FunA\not\equiv 0$ on $\dom(\FunA)$.
Then there exists $\imageA\in\dom(\FunA)$ with $\FunA(\imageA)\neq 0$ and thanks to \cref{prop:positivity} it holds $\FunA(\imageA)\geq 0$.
If we assume that $\FunA(\imageA)>0$, then the strict concavity of the function $t\mapsto t^p$ and the absolute $p$-homogeneity and convexity of $\FunA$ imply that for $\lambda\in(0,1)$ and $s\neq t\in\R\setminus\{0\}$ it holds
\begin{align*}
    (\lambda s^p + (1-\lambda)t^p)\FunA(\imageA) 
    &<
    (\lambda s + (1-\lambda) t)^p \FunA(\imageA)
    \\
    &=
    \FunA((\lambda s + (1-\lambda) t)\imageA)
    \\
    &\leq 
    \FunA(\lambda s \imageA + (1-\lambda) t\imageA)
    \\
    &\leq
    \lambda\FunA(s\imageA) + (1-\lambda)\FunA(t\imageA)
    \\
    &=
    (\lambda s^p +(1-\lambda) t^p)\FunA(\imageA).
\end{align*}
This is a contradiction and hence $\FunA(\imageA)=0$ has to hold.
\end{proof}

Another important property of absolutely $p$-homogeneous functionals is the Euler rule.
\begin{theorem}[Euler's homogeneous function theorem]\label{thm:euler_homogeneous}
Let $\FunA:\Bspace\to\RI$ be absolutely $p$-homogeneous for $p>0$.
Then for any $\imageA\in\Bspace$ and $\subgradA\in\partial\FunA(\imageA)$ it holds
\begin{align}\label{eq:euler}
    p\FunA(\imageA) = \sp{\subgradA,\imageA}.
\end{align}
\end{theorem}
\begin{proof}
By definition of the subdifferential it holds
\begin{align*}
    \FunA(\imageA) + \sp{\subgradA,\imageB-\imageA} \leq \FunA(\imageB),\quad\forall \imageB\in\Bspace.
\end{align*}
Choosing $\imageB=cu$ for $c>0$ yields
\begin{align*}
    \FunA(\imageA) + (c-1) \sp{\subgradA,\imageA} \leq c^p \FunA(\imageA).
\end{align*}
which can be rewritten to
\begin{align*}
    (1-c^p)\FunA(\imageA) \leq (1-c)\sp{\subgradA,\imageA}.
\end{align*}
Hence, we get
\begin{align*}
    \frac{1-c^p}{1-c}\FunA(\imageA) &\leq \sp{\subgradA,\imageA},\quad \text{if } c<1,\\
    \frac{1-c^p}{1-c}\FunA(\imageA) &\geq \sp{\subgradA,\imageA},\quad \text{if } c>1.
\end{align*}
Using that $\lim_{c\to 1} \frac{1-c^p}{1-c}=p$ concludes the proof.
\end{proof}

Furthermore, the subdifferential is also homogeneous, however, of degree $p-1$:
\begin{proposition}\label{prop:subdifferential_homogeneous}
Let $\FunA:\Bspace\to\RI$ be absolutely $p$-homogeneous for $p>0$.
Then it holds
\begin{align*}
    \partial\FunA(c\imageA) = c \abs{c}^{p-2} \partial\FunA(\imageA),\quad\forall c\neq 0.
\end{align*}
\end{proposition}
\begin{proof}
By definition it holds
\begin{align*}
    \partial\FunA(c\imageA) 
    &= 
    \left\lbrace
    \subgradA\in\dual{\Bspace} \st \FunA(c\imageA) + \sp{\subgradA,\imageB-c\imageA} \leq
    \FunA(\imageB)\quad\forall\imageB\in\Bspace
    \right\rbrace
    \\
    &=
    \left\lbrace
    \subgradA\in\dual{\Bspace} \st \abs{c}^p\FunA(\imageA) + \sp{\subgradA,\imageB-c\imageA} \leq
    \FunA(\imageB)\quad\forall\imageB\in\Bspace
    \right\rbrace
    \\
    &=
    \left\lbrace
    \subgradA\in\dual{\Bspace} \st \FunA(\imageA) + \abs{c}^{-p}\sp{\subgradA,\imageB-c\imageA} \leq
    \abs{c}^{-p}\FunA(\imageB)\quad\forall\imageB\in\Bspace
    \right\rbrace
    \\
    &=
    \left\lbrace
    \subgradA\in\dual{\Bspace} \st \FunA(\imageA) + c\abs{c}^{-p}\sp{\subgradA,\imageB/c-\imageA} \leq
    \FunA(\imageB/c)\quad\forall\imageB\in\Bspace
    \right\rbrace
    \\
    &=
    \left\lbrace
    \subgradA\in\dual{\Bspace} \st \FunA(\imageA) + c^{-1}\abs{c}^{2-p}\sp{\subgradA,\imageB-\imageA} \leq
    \FunA(\imageB)\quad\forall\imageB\in\Bspace
    \right\rbrace
    \\
    &=
    \left\lbrace
    c\abs{c}^{p-2}\subgradB\in\dual{\Bspace} \st \FunA(\imageA) + \sp{\subgradB,\imageB-\imageA} \leq
    \FunA(\imageB)\quad\forall\imageB\in\Bspace
    \right\rbrace
    \\
    &= 
    c\abs{c}^{p-2} \partial\FunA(\imageA).
\end{align*}
\end{proof}
The following result shows that zero level set of a convex homogeneous functional is a linear spaces, referred to as \emph{nullspace}.

\begin{lemma}\label{lem:null-space-J}
Let $\FunA(\imageA) = \sum_{i=1}^n \FunA_i(\imageA)$, where each $\FunA_i(\imageA)$ is convex, absolutely $p_i$-homogeneous ($p_i > 0$), and \lsc{.}
Then the set
\begin{equation*}
\nullspace_{\FunA} \defeq \{ \imageA \in \Bspace \st \FunA(\imageA) = 0\} 
\end{equation*}
is a closed linear subspace of $\Bspace$.
\end{lemma}
\begin{proof}
Let $u,v \in \nullspace_{\FunA}$ be arbitrary. Then $\FunA_i(\imageA) = \FunA_i(\imageB) = 0$ for all $i=1,...,n$.
Hence, using also \cref{prop:positivity}, it holds for any $\alpha \in \R$ 
\begin{eqnarray*}
0 \leq \FunA_i(\alpha u + v) &=& 2^{p_i}  \FunA_i\left(\frac{\alpha u}{2} + \frac{v}{2}\right) \leq 2^{p_i}  \left(\frac12 \FunA_i\left(\frac{\alpha u}{2}\right) + \frac12\FunA_i\left(\frac{v}{2}\right)\right) \\
&=& \frac12 \FunA_i(\alpha u) + \frac12\FunA_i(\imageB) = \frac{\abs{\alpha}^{p_i}}{2} \FunA_i(\imageA) + \frac12\FunA_i(\imageB) = 0.
\end{eqnarray*}
Therefore, $\FunA_i(\alpha u + v) = 0$ for all $i$ and hence $\FunA(\alpha u + v) = 0$. Closedness of $\nullspace_{\FunA}$ follows from the lower-semicontinuity of $\FunA$.
\end{proof}

Next we prove several important properties of the nullspace and its interaction with the subdifferential, taken from \cite{bungert2019asymptotic}.

\begin{proposition}\label{prop:orth}
If $\FunA$ is absolutely $p$-homogeneous for $p>0$, it holds for any $\imageA\in\Bspace$ that $\partial \FunA(\imageA)\subset\nullspace_\FunA^\perp$.
\end{proposition} 
\begin{proof}
Let $u\in\Bspace$ and $\imageB_0\in\nullspace_\FunA$ be arbitrary.
Using the definition of the subdifferential together with \cref{thm:euler_homogeneous} yields
\begin{align*}
    0 
    &= 
    \FunA(\imageB_0) 
    \geq \FunA(\imageA) + \sp{\subgradA,\imageB_0 - u}
    \\
    &= 
    \FunA(\imageA) + \sp{\subgradA,\imageB_0} - p \FunA(\imageA) 
    = (1-p)\FunA(\imageA) + \sp{\subgradA,\imageB_0}
    ,\quad\forall\subgradA\in\partial \FunA(\imageA).
\end{align*}
Applying the same argument to $-\imageB_0\in\nullspace_\FunA$ yields
\begin{align*}
|\langle\subgradA,\imageB_0\rangle|\leq(p-1)\FunA(\imageA),\quad\forall u\in\Bspace,\;\forall\subgradA\in\partial \FunA(\imageA),\;\forall\imageB_0\in\nullspace_\FunA.    
\end{align*}
If $0<p<1$, this is a contradiction unless $\FunA(\imageA)=0$ in which case $\sp{\subgradA,\imageB}=0$ holds.
If $p=1$, we can already conclude that $\sp{\subgradA,\imageB}=0$.
In the case $p>1$ we argue as follows:
Thanks to \cref{prop:subdifferential_homogeneous} it holds $c^{p-1}\subgradA\in\partial \FunA(c \imageA)$ for all $c>0$. Replacing $\subgradA$ by $c^{p-1}\subgradA$ and $\imageA$ by $c\imageA$ in the inequality above yields
\begin{align*}
    |\langle c^{p-1}\subgradA,\imageB_0\rangle|\leq (p-1)\FunA(c\imageA)
    \iff
    |\langle\subgradA,\imageB_0\rangle|\leq c(p-1)\FunA(\imageA)
    ,\quad\forall\subgradA\in\partial \FunA(\imageA),
\end{align*}
by using the homogeneity of $\FunA$ and dividing by $c^{p-1}$. Letting $c\searrow 0$ concludes the proof.
\end{proof}

The following proposition states that the value of the functional $J$ and its subdifferential is invariant under addition of a nullspace element.

\begin{proposition}\label{prop:subdiff_nullspace}
Let $\FunA:\Bspace\to\RI$ be absolutely $p$-homogeneous for $p\geq 1$, convex, and lower semi-continuous.
For any $\imageB_0\in\nullspace_\FunA$ and $\imageA\in\Bspace$ it holds $\FunA(\imageA)=\FunA(\imageA+\imageB_0)$ and  for any $\subgradA\in\partial\FunA(\imageA)$ it holds $\partial \FunA(\imageA)=\partial \FunA(\imageA+\imageB_0)$.
\end{proposition}
\begin{proof}
We claim that the convex conjugate $\conj{\FunA}$ meets $\conj{\FunA}(\subgradA)=\infty$ whenever $\subgradA\notin\nullspace_\FunA^\perp$. 
To see this we note that for any $\imageB_0\in\nullspace_\FunA$ it holds
\begin{align*}
\conj{\FunA}(\subgradA)=\sup_{u\in\Bspace}\langle\subgradA,u\rangle-\FunA(\imageA)\geq \langle\subgradA,\imageB_0\rangle.
\end{align*}
If $\langle\subgradA,\imageB_0\rangle\neq 0$ we can use that $\nullspace_\FunA$ is a linear space and replace $\imageB_0$ by $c\imageB_0$ for $c\in\R$ to obtain that $\conj{\FunA}(\subgradA)=\infty$.
Using this together with the fact that $\FunA$ is lower-semicontinuous and convex and therefore equals its biconjugate, we obtain
\begin{align*}
\FunA(\imageA)=\biconj{\FunA}(\imageA)=\sup_{\subgradA\in\dual{\Bspace}}\langle\subgradA,u\rangle - \conj{\FunA}(\subgradA) = \sup_{\subgradA\in\nullspace_\FunA^\perp}\langle\subgradA,u\rangle - \conj{\FunA}(\subgradA),\quad\forall \imageA\in\Bspace.
\end{align*}
From here it is obvious that $\FunA(\imageA+\imageB_0)=\FunA(\imageA)$ for all $\imageB_0\in\nullspace_\FunA$.

For the second statement we observe that
\begin{align*}
    \partial\FunA(\imageA+\imageB_0) 
    &=
    \left\lbrace
    \subgradA \in \dual{\Bspace}
    \st
    \FunA(\imageA+\imageB_0) + 
    \sp{\subgradA,\imageB - (\imageA+\imageB_0)}
    \leq 
    \FunA(\imageB)\quad\forall\imageB\in\Bspace
    \right\rbrace
    \\
    &=
    \left\lbrace
    \subgradA \in \dual{\Bspace}
    \st
    \FunA(\imageA) + 
    \sp{\subgradA,\imageB - \imageA}
    \leq 
    \FunA(\imageB)\quad\forall\imageB\in\Bspace
    \right\rbrace
    =
    \partial\FunA(\imageA).
\end{align*}
\end{proof}

If $\Bspace=\Hspace$ is a Hilbert space, then $\nullspace_\FunA$ is a closed subspace of a Hilbert space and hence \emph{complemented} in $\spaceA$~\cite[Thm. 5.89]{aliprantis-border:2006}, i.e. there exists a closed subspace $\spaceA_0 \subset \spaceA$ such that $\spaceA_0 \cap \nullspace_\FunA = \{0\}$ and
\begin{equation}\label{eq:direct-sum}
\spaceA = \spaceA_0 \oplus \nullspace_\FunA.
\end{equation}
In fact, $\Hspace_0$ equals the orthogonal complement $\nullspace_\FunA^\perp$.
In a Banach space setting, one can consider the \emph{quotient space} $\Bspace_{\sim} := \Bspace / \nullspace_\FunA$ whose elements $[\imageA]$ are equivalence classes with respect to the equivalence relation $\imageA\sim\imageB$ if and only if $\imageA-\imageB\in\nullspace_\FunA$.
A norm on $\Bspace_\sim$ is given by $\norm{[\imageA]}_\sim\defeq\inf_{\imageB_0\in\nullspace_\FunA}\norm{\imageA+\imageB_0}$.
If $\Bspace=\Hspace$ is a Hilbert space, it holds that $\Hspace_0$ and $\Hspace_\sim$ are isomorphic.
\begin{remark}
One can use the previous insights to define an absolutely $p$-homogeneous, convex, and \lsc{} functional on the quotient space as
\begin{align*}
    \FunA_\sim([\imageA]) := \FunA(\imageA),\quad\forall[\imageA]\in\Bspace_\sim.
\end{align*}
This is well-defined by \cref{prop:subdiff_nullspace}.
Furthermore, the nullspace of $\FunA_\sim$ is given by $\nullspace_{\FunA_\sim}=\{[0]\}$.
\end{remark}

\begin{example}
Let $\Omega\subset\R^n$ be a domain with Lipschitz boundary and $\Hspace=L^2(\Omega)$. Let $p \geq \frac{2n}{n+2}$ (in this case the Rellich-Kondrachov embedding theorem ensures that $W^{1,p}(\Omega) \hookrightarrow L^2(\Omega)$; for $n=2$ we have $p \geq 1$). 
Then the functional
\begin{align*}
    \calJ_p : L^2(\Omega)\to\RI,\quad
    u \mapsto 
    \begin{cases}
    \frac{1}{p}\int_\Omega\abs{\grad u}^p\d x,\quad u \in W^{1,p}(\Omega)\cap L^2(\Omega), \\
    +\infty, \quad u \in L^2(\Omega)\setminus W^{1,p}(\Omega),
    \end{cases}
\end{align*}
is proper, absolutely $p$-homogeneous, convex, and lower semicontinuous. 
Its nullspace is given by $\nullspace_{\calJ_p} = \{u \in L^2(\Omega) \colon u\equiv const\}$.
\end{example}

If~\eqref{eq:p-hom} is satisfied with $p=1$, the functional is called absolutely one-homogeneous.  Some of their properties are listed below.

\begin{proposition}\label{prop:dJ(0)}
Let $\FunA:\Bspace\to\RI$ be a proper and absolutely one-homogeneous functional.
Then the convex conjugate $\conj{\FunA}:\dual{\Bspace}\to\RI$ is the characteristic function of the convex set $\subdiff\FunA(0)$, in other words
\begin{align*}
    \conj{\FunA}(\subgradA) = \chi_{\partial\FunA(0)}(\subgradA).
\end{align*}
\end{proposition}
\begin{proof}
Left as exercise.
\end{proof}

An obvious consequence of this result and \cref{thm:fenchel-young} is the following
\begin{proposition}\label{prop:char_of_dJ(u)}
Let $\FunA:\Bspace\to\RI$ be a proper and absolutely one-homogeneous functional.
Let $\imageA \in \Bspace$. Then it holds
\begin{align*}
    \subdiff\FunA(\imageA) = \left\lbrace\subgradA\in\dual{\Bspace}\st \sp{\subgradA,\imageA}=\FunA(\imageA),\,\sp{\subgradA,\imageB}\leq\FunA(\imageB)\;\forall\imageB\in\Bspace\right\rbrace.
\end{align*}
\end{proposition}

Furthermore, we obtain the following \emph{dual representation} of absolutely one-homogeneous functionals.
\begin{proposition}\label{prop:dual_one-hom}
Let $\FunA:\Bspace\to\RI$ be a proper, absolutely one-homogeneous, convex, and \lsc{} functional.
Then it holds
\begin{align*}
    \FunA(\imageA) = \sup_{\subgradA\in\partial\FunA(0)}\sp{\subgradA,\imageA}.
\end{align*}
\end{proposition}
\begin{proof}
Since $\FunA$ equals its biconjugate we obtain from \cref{prop:dJ(0)} that
\begin{align*}
    \FunA(\imageA) = \biconj{\FunA}(\imageA) = \sup_{\subgradA\in\dual{\Bspace}}\sp{\subgradA,\imageA}-\conj{\FunA}(\subgradA) =
    \sup_{\subgradA\in\dual{\Bspace}}\sp{\subgradA,\imageA}-\chi_{\partial\FunA(0)}(\subgradA)
    =
    \sup_{\subgradA\in\partial\FunA(0)}\sp{\subgradA,\imageA}.
\end{align*}
\end{proof}

\section{Rayleigh quotient and nonlinear optimality condition}

From now on we will consider a functional $\calJ:\Hspace\to\RI$ on a Hilbert space which is absolutely $p$-homogeneous with $p\geq 1$, convex, proper, and \lsc{.}
We will assume this throughout the whole chapter without repeating the assumption.
We shall study the following minimisation problem for the Rayleigh quotient (cf.~\eqref{eq:Rayleigh-opt-prob-lin})
\begin{equation}\label{eq:Rayleigh-opt-prob}
    \lambda_1 := \inf_{u \in \nullspace_\calJ^\perp\setminus\{0\}} \frac{p\calJ(u)}{\norm{u}^p},
\end{equation}
where $\nullspace_\calJ^\perp$ denotes the orthogonal complement of the null space $\nullspace_\calJ$.
The reason for performing the minimisation over $\nullspace_\calJ^\perp$ is to avoid that every element in $\nullspace_\calJ$ is a minimiser.
We will soon see that the value $\lambda_1$ can be interpreted as \emph{minimal eigenvalue} and any minimiser of \labelcref{eq:Rayleigh-opt-prob} as the corresponding eigenfunction.
Hence, as usual in physics, we call these minimal eigenfunction \emph{ground states}.

\begin{definition}[Ground states]\label{def:ground_states}
We refer to minimisers of \labelcref{eq:Rayleigh-opt-prob} as \emph{ground states}.
\end{definition}

Under reasonable conditions on $\calJ$ and its interplay with the Hilbert space topology, we can prove existence of ground states.

\begin{theorem}[Existence of ground states]\label{thm:existence_GS}
Assume that $\dom(\calJ)\cap\nullspace_\calJ^\perp\setminus\{0\}\neq\emptyset$.
If the sublevel sets of $u\mapsto\norm{u}+\calJ(u)$ are precompact, then the optimisation problem~\eqref{eq:Rayleigh-opt-prob} admits a minimiser.
Furthermore, if $u\in\Hspace$ is a minimiser, so is $C\cdot u$ for any $C\neq 0$.
\end{theorem}
\begin{proof}
The absolute $p$-homogeneity of $\calJ$ and $\norm{\cdot}^p$ implies
\begin{align*}
    \frac{p\calJ(C\cdot u)}{\norm{C\cdot u}^p} = \frac{p\abs{C}^p\calJ(u)}{\abs{C}^p\norm{u}^p}
    =
    \frac{p\calJ(u)}{\norm{u}^p},\quad C\neq 0.
\end{align*}
Hence, multiples of minimisers are minimisers.
In particular, we can study the equivalent minimisation problem 
\begin{align}\label{eq:equivalent_GS}
    \inf \left\lbrace\calJ(u) \st u\in\nullspace_\calJ^\perp\setminus\{0\},\; \norm{u} = 1\right\rbrace.
\end{align}
By assumption, there exists $u\in\nullspace_\calJ^\perp\setminus\{0\}$ with $\calJ(u)<\infty$ and hence also $\calJ(u/\norm{u}) = \calJ(u)/\norm{u}^p<\infty$.
Hence, we can choose a minimising sequence $(u_k)_{k\in\N}\subset\Hspace$ with $\norm{u_k}=1$ such that $\calJ(u_k)$ converges to the infimal value in \labelcref{eq:equivalent_GS} as $k\to\infty$.
Since sublevel sets of $\norm{\cdot}+\calJ$ are precompact, there exists a subsequence of $(u_k)_{k\in\N}$ (which we don't relabel) which converges to some $u\in\Hspace$ with $\norm{u}=1$.
Since $\nullspace_\calJ^\perp$ is closed it holds $u\in\nullspace_\calJ^\perp\setminus\{0\}$.
Since $\calJ$ is \lsc{}, we have
\begin{align*}
    \calJ(u) \leq \liminf_{n\to\infty}\calJ(u_k) = \inf \left\lbrace\calJ(u) \st u\in\nullspace_\calJ^\perp\setminus\{0\},\; \norm{u}= 1\right\rbrace.
\end{align*}
Hence, $u$ solves \labelcref{eq:equivalent_GS} and therefore also \labelcref{eq:Rayleigh-opt-prob}.
\end{proof}

\begin{example}
The precompactness of the sublevel sets of $\norm{\cdot}+\calJ$ is necessary. 
If, for instance, $\calJ(u) = \int_{0}^1\abs{u}\d x$ and $\Hspace=L^2((0,1))$, the compactness is violated. 
Furthermore, the sequence $u_k(x) = k\max(1-kx,0)$ satisfies $\calJ(u_k)=1$ and $\norm{u_k}_{L^2}=k$ and hence $\calJ(u_k)/\norm{u_k}_{L^2} \to 0$.
Therefore, $(u_k)_{k\in\N}$ is a minimising sequence for \labelcref{eq:Rayleigh-opt-prob} but does not converge.
\end{example}

\begin{example}
It can indeed happen that $\calJ$ is proper but $\lambda_1=\infty$. 
This is the case for $\calJ(u) = \chi_{\{0\}}(u)$.
\end{example}

Next we derive optimality conditions for~\eqref{eq:Rayleigh-opt-prob}, cf. \cref{thm:opt-cond-lin}.
\begin{theorem}\label{thm:opt-cond-nonlin}
It holds that $u\in\nullspace_\calJ^\perp\setminus\{0\}$ is a minimiser of~\eqref{eq:Rayleigh-opt-prob} if and only if it satisfies
\begin{equation}\label{eq:opt-cond-nonlin}
    \lambda_1 \norm{u}^{p-2} u \in \dJ(u),
\end{equation}
where $\lambda_1 \defeq \frac{p\calJ(u)}{\norm{u}^p}$.
\end{theorem}
\begin{proof}
Let us assume that $u$ minimises \labelcref{eq:Rayleigh-opt-prob}.
Then it holds for any $v\in\nullspace_\calJ^\perp\setminus\{0\}$ that
\begin{align*}
    \calJ(u) + \sp{\lambda_1\norm{u}^{p-2}u,v-u}
    &=
    \calJ(u) + \lambda_1 \left(\norm{u}^{p-2}\sp{u,v} 
    - \norm{u}^p
    \right)
    \\
    &\leq 
    \calJ(u)
    + \lambda_1
    \left(
    \norm{u}^{p-1}\norm{v}-\norm{u}^p
    \right)
    \\
    &\leq
    \calJ(u) + \lambda_1
    \norm{u}^{p-1}\norm{v}
    - p\calJ(u)
    \\
    &\leq 
    (1-p)\calJ(u) + 
    \lambda_1
    \left(
    \frac{p-1}{p}
    \norm{u}^p
    +
    \frac{1}{p}
    \norm{v}^p
    \right)
    \\
    &=
    (1-p)\calJ(u) + 
    (p-1)\calJ(u)
    +
    \frac{p\calJ(v)}{\norm{v}^p}
    \frac{1}{p}
    \norm{v}^p
    \\
    &=
    \calJ(v),
\end{align*}
where we used Young's inequality $ab \leq a^p/p + b^q/q$ for any $a,b \geq 0$, $p > 1$ and $q = \frac{p}{p-1}$ to get from the third to the forth line (for $p=1$, this step is trivial). 
By \cref{prop:subdiff_nullspace} we see that this is in fact true for any $v\in\Hspace\setminus\{0\}$.
For $v=0$ it trivially holds
\begin{align*}
    \calJ(u) + \sp{\lambda_1\norm{u}^{p-2}u,0-u}
    =
    \calJ(u) - \lambda_1\norm{u}^p
    =
    \calJ(u) - p\calJ(u)
    \leq 0.
\end{align*}
which lets us conclude that $\lambda_1\norm{u}^{p-2}u\in\partial\calJ(u)$.

For the converse implication we assume that $\lambda_1\norm{u}^{p-2}u\in\partial\calJ(u)$, where $\lambda_1$ is the minimal value of the Rayleigh quotient defined in \labelcref{eq:Rayleigh-opt-prob}.
Then we can use \cref{thm:euler_homogeneous} to obtain
\begin{align*}
    p\calJ(u) = \sp{\lambda_1\norm{u}^{p-2}u,u}=\lambda_1\norm{u}^p
\end{align*}
which is equivalent to
\begin{align*}
    \frac{p\calJ(u)}{\norm{u}^p} = \lambda_1
\end{align*}
and means that $u$ solves the minimisation problem in \labelcref{eq:Rayleigh-opt-prob}.
\end{proof}

We recall that \cref{thm:opt-cond-nonlin} implies that eigenfunctions with minimal eigenvalue correspond to minimisers of the Rayleigh quotient. 
Eigenfunctions with larger eigenvalues will, of course, not be minimisers anymore. 
Still, they can be interpreted as stationary points of the Rayleigh quotient in the sense of the  Clarke subdifferential~\cite{rockafellar:1980}.

This suggests the following definition. 
\begin{definition}\label{def:ef-nonlin}
Let $\calJ \colon \Hspace \to \R \cup \{+\infty\}$ a proper, convex, \lsc{,} and absolutely $p$-homogeneous functional.  An element $u \in \Hspace\setminus\{0\}$ is called an \ef{} w.r.t. $\calJ$ with \ev{} $\lambda \in \R$ if
\begin{equation}\label{eq:opt-con-Hilbert}
    \lambda \norm{u}^{p-2} u \in \dJ(u).
\end{equation}
\end{definition}

In the following proposition we collect important properties of nonlinear eigenfunctions and their eigenvalues which are trivial consequences of the results we already have.
\begin{proposition}
Let $u\in\Hspace\setminus\{0\}$ be and eigenfunction of $\calJ$ with eigenvalue $\lambda\in\R$.
Then it holds
\begin{enumerate}
    \item $u\in\nullspace_\calJ^\perp$;
    \item $\lambda\geq 0$;
    \item For all $c\in\R$ it holds that $cu$ is an eigenfunction of $\calJ$ with eigenvalue $\lambda$. 
\end{enumerate}
\end{proposition}

Eigenfunctions posses a remarkable property, namely that they are the minimal norm elements in their own subdifferential.
\begin{proposition}\label{prop:ef_minimal_norm}
Let $u\in\Hspace\setminus\{0\}$ satisfy $\lambda\norm{u}^{p-2}u\in\partial\calJ(u)$. Then it holds
\begin{align*}
    \lambda\norm{u}^{p-2}u = \argmin
    \left\lbrace
    \norm{\subgradA} \in \Hspace \st \subgradA \in \dJ(u)
    \right\rbrace.
\end{align*}
\end{proposition}
\begin{proof}
First we note that $\partial\calJ(u)$ is convex and closed, and therefore has a unique minimal norm element.

Let us use the abbreviation $\subgradA_0:=\lambda\norm{u}^{p-2}u\in\dJ(u)$ and choose an arbitrary $\subgradA\in\dJ(u)$.
Then it holds
\begin{align*}
    \norm{\subgradA_0}^2 
    &= 
    \sp{\subgradA_0,\lambda\norm{u}^{p-2}u} 
    = 
    \lambda\norm{u}^{p-2}\sp{\subgradA_0,u} 
    = 
    \lambda\norm{u}^{p-2} p\calJ(u) 
    \\
    &=
    \lambda\norm{u}^{p-2} \sp{\subgradA,u} 
    \leq 
    \lambda\norm{u}^{p-1} \norm{\subgradA}.
\end{align*}
It holds $\norm{\subgradA_0}=\lambda\norm{u}^{p-1}$ and we get $\norm{\subgradA_0}\leq\norm{\subgradA}$.
Since $\subgradA\in\partial\calJ(u)$ was arbitrary, the proof is complete. 
\end{proof}

\section[Eigenfunctions of  one-homogeneous functionals]{Nonlinear eigenfunctions of absolutely one-homoge\-neous functionals}

In this section we continue with studying properties of nonlinear eigenfunctions associated with absolutely one-homogeneous functions. 
For such functionals the eigenvalue problem has particularly nice geometric interpretations.

\begin{definition}
Let $\calJ \colon \Hspace \to \R \cup \{+\infty\}$ a proper, convex, \lsc{,} and absolutely one-homogeneous functional. The following set
\begin{equation*}
    \dualball \defeq \dJ(0)
\end{equation*}
is called the \emph{dual unit ball} of $\calJ$. 
\end{definition}

The terminology is justified by \cref{prop:char_of_dJ(u)} which states that $\dualball$ has the following characterization:
\begin{align*}
    \dualball 
    &=
    \left\lbrace
    \subgradA \in {\Hspace} \st \sp{\subgradA,u}\leq\calJ(u)\;\forall u\in\Hspace
    \right\rbrace
    \\
    &=
    \left\lbrace
    \subgradA \in {\Hspace} \st \sup_{\substack{u\in\Hspace\\\calJ(u)\leq 1}}\sp{\subgradA,u}\leq 1
    \right\rbrace.
\end{align*}
It is easily checked that, just as $\calJ$, also the mapping
\begin{align}\label{eq:dual_seminorm}
\subgradA\mapsto\calJ_*(\subgradA):=\sup_{\substack{u\in\Hspace\\\calJ(u)\leq 1}}\sp{\subgradA,u}    
\end{align}
is absolutely one-homogeneous and convex. 
Hence, with a grain of salt, owing to the fact that $\calJ$ and $\calJ_*$ can take infinite values, one can interpret $\calJ$ as semi-norm and $\calJ_*$ as the associated dual semi-norm.

The reader is encouraged to check that $\dualball \subset \Hspace$ is a closed convex set.
Furthermore, the characterization of the subdifferential $\dJ(u)$ from \cref{prop:char_of_dJ(u)} can be reformulated as
\begin{align}\label{eq:one-hom_subdiff_dualball}
    \dJ(u) =  
    \left\lbrace\subgradA\in\dualball\st\sp{\subgradA,u}=\calJ(u)\right\rbrace,
\end{align}
meaning that subgradients of an absolutely one-homogeneous functional are distinguished points in the dual unit ball.

\begin{example}[One-homogeneous functionals and their dual unit balls]\label{ex:dual_unit_ball}
Some examples for absolutely one-homogeneous functionals and their duality are given below.
\begin{enumerate}
    \item If $\calJ(u):=\norm{u}_{p}:=\left(\sum_{i=1}^d\abs{u_i}^p\right)^{1/p}$ for $p\in[1,\infty]$ is the $p$-norm of a vector $u\in\R^n$, then $\calJ(\subgradA)_\ast=\norm{\subgradA}_{p'}$, where $p':=\frac{p}{p-1}$ for $p>1$, $1':=\infty$, and $\infty':=1$.
    \item Let $K\subset\Hspace$ be a set which satisfies $-K=K$ and define $\calJ(u):=\sup_{\subgradA\in K}\sp{\subgradA,u}$. Then $\calJ$ is absolutely one-homogeneous and convex, and the dual unit ball of $\calJ$ is given by the closed convex hull of $K$, i.e., $\dualball=\overline{\operatorname{conv}}(K)$.
    \item Let $\calJ$ be given by the total variation of a function on a domain $\Omega\subset\R^2$, i.e.,
    \begin{align*}
        \calJ : L^2(\Omega) \to \RI,\quad
        \calJ(u) := \sup\left\lbrace
        \int_\Omega u \div\phi \d x \st \phi\in C^\infty_c(\Omega;\R^2),\;\sup_\Omega\abs{\phi}\leq 1
        \right\rbrace.
    \end{align*}
    Then the previous example tells us that $\calJ$ is absolutely one-homogeneous and convex. 
    Furthermore, the set $K:=\{\div\phi\st\phi\in C^\infty_c(\Omega;\R^2),\;\sup_\Omega\abs{\phi}\leq 1 \}$ is convex, hence the dual unit ball of $\calJ$ coincides with its closure $\overline{K}$ in $L^2(\Omega)$.
    This closure is characterized in \cite{bredies2016pointwise} as
    \begin{align*}
        \dualball = \overline{K} = \left\{\div\phi\st\phi\in W^{1,2}_0(\operatorname{div};\R^2),\;\sup_\Omega\abs{\phi}\leq 1 \right\},
    \end{align*}
    where the space $W^{1,2}_0(\operatorname{div};\R^2)$ is defined as the closure of the space $C^\infty_c(\Omega;\R^2)$ with respect to the norm
    \begin{align*}
        \norm{\phi}_{W^2(\div)}^2 := \norm{\abs{\phi}}_{L^2(\Omega)}^2 + \norm{\div\phi}_{L^2(\Omega)}^2.
    \end{align*}
\end{enumerate}
\end{example}

The following result gives a useful geometric characterization of \ef{s}.
\begin{theorem}\label{thm:ef-equiv-cond}
Let $\calJ \colon \Hspace \to \R \cup \{+\infty\}$ a proper, convex, \lsc{,} and absolutely one-homogeneous functional.
If element $u \in \Hspace$ is an \ef{} of $\calJ$ with \ev{} $\lambda \in \R$ then the following inequality holds for $\subgradA \defeq \lambda u/\norm{u}\in\dualball$:
\begin{align}\label{ineq:ef-equiv-cond}
    \sp{\subgradA,\subgradA-\eta} \geq 0 \quad \forall \, \subgradB \in \dualball.
\end{align}
Conversely, if \labelcref{ineq:ef-equiv-cond} holds for some $\subgradA\in\dualball$, then it holds $\subgradA\in\dJ(\subgradA)$.
\end{theorem}
\begin{proof}
By \cref{prop:subdifferential_homogeneous} $\lambda u/\norm{u}\in\dJ(u)$ implies $\subgradA\in\dJ(\subgradA)$ which implies
\begin{align*}
    \sp{\subgradA,\subgradA}
    =
    \calJ(\subgradA)
    =
    \sup_{\subgradB\in K_\calJ}\sp{\subgradB,\subgradA}
    \geq 
    \sp{\subgradB,\subgradA},\quad\forall\subgradB\in K_\calJ.
\end{align*}
Conversely, if $\subgradA\in\partial\calJ(\subgradA)$ and \labelcref{ineq:ef-equiv-cond} holds then
\begin{align*}
    J(\subgradA)
    \geq 
    \sp{\subgradA,\subgradA}
    \geq 
    \sp{\subgradA,\subgradB},\quad\forall\subgradB\in K_\calJ.
\end{align*}
Taking the supremum over $\subgradB\in K_\calJ$ shows
\begin{align*}
    \sp{\subgradA,\subgradA} = J(\subgradA)
\end{align*}
and hence $\subgradA\in\partial\calJ(\subgradA)$.
\end{proof}

Now we investigate for which elements $\subgradB \in \dualball$ the characterizing inequality~\eqref{ineq:ef-equiv-cond} is actually an equality.

\begin{proposition}\label{prop:minsub}
Let $u\in\dom(\partial J)$ and $\subgradA:=\argmin\left\lbrace\norm{\subgradB}\st \subgradB\in\partial J(u)\right\rbrace$ be an eigenvector with eigenvalue 1. Then it holds
\begin{align}\label{eq:minsub_general}
\langle \subgradA,\subgradA-\subgradB\rangle =0,\quad\forall \subgradB\in\partial J(u),
\end{align}
which can be reformulated as $\partial J(u)\subset\partial J(\subgradA).$
\end{proposition}
\begin{proof}
The non-negativity of the left-hand side follows directly from the assumption that $\subgradA$ is an eigenvector and thus fulfills~\eqref{ineq:ef-equiv-cond}. For the other inequality, we let $\subgradB\in\partial J(u)$ arbitrary and consider $\subgradC:=\lambda \subgradB+(1-\lambda)\subgradA$ for $\lambda\in(0,1)$, which is in $\partial J(u)$ as well due to convexity. Using \eqref{ineq:ef-equiv-cond} and the minimality of $\norm{\subgradA}$, we get
\begin{align*}
\lambda\sp{\subgradA,\subgradB} + (1-\lambda)\norm{\subgradA}^2 = \sp{\subgradA,\subgradC} \leq \norm{\subgradA}^2 \leq \norm{\subgradC}^2 = \lambda^2 \norm{\subgradB}^2 + (1-\lambda)^2 \norm{\subgradA}^2 + 2\lambda(1-\lambda) \sp{\subgradA,\subgradB}.
\end{align*}
Dividing this by $\lambda(1-\lambda)$ one finds
\begin{align*}
\frac{1}{1-\lambda}\sp{\subgradA,\subgradB} + \frac{1}{\lambda}\norm{\subgradA}^2\leq\frac{\lambda}{1-\lambda}\norm{\subgradB}^2 + \frac{1-\lambda}{\lambda}\norm{\subgradA}^2 + 2\sp{\subgradA,\subgradB},
\end{align*}
which can be simplified to 
\begin{align*}
\frac{1}{1-\lambda}\sp{\subgradA,\subgradB} \leq \frac{\lambda}{1-\lambda}\norm{\subgradB^2}-\norm{\subgradA}^2 + 2\sp{\subgradA,\subgradB}.
\end{align*}
Letting $\lambda$ tend to zero and reordering shows $\sp{ \subgradA,\subgradA-\subgradB} \leq 0$, hence equality holds.
\end{proof}

Next we prove that subdifferentials are contained in the relative boundary of the dual unit ball $\dualball$.

\begin{proposition}\label{prop:subdiff_bdry}
Let $\calJ \colon \Hspace \to \R \cup \{+\infty\}$ a proper, convex, \lsc{,} and absolutely one-homogeneous functional and let $\subgradA\in\dJ(u)$ for $u\in\nullspace_\calJ^\perp\setminus\{0\}$.
Then $\subgradA\in\partial_\mathrm{rel}\dualball$ where the \emph{relative boundary} of a convex set $C\subset\Hspace$ is defined as $\partial_\mathrm{rel}C:=\closure{C}\setminus\{u\in C\st \exists c>1 \text{ with }cu\in C\}$.
\end{proposition}
\begin{proof}
Since $\dualball$ is closed and convex, by \labelcref{eq:one-hom_subdiff_dualball} it suffices to show that for $c>1$ we have $c\subgradA\notin\dualball$.
Assuming the contrary, we have
\begin{align*}
    c\sp{\subgradA,u}=\sp{c\subgradA,u}\leq\calJ(u) = \sp{\subgradA,u}
\end{align*}
which is a contradiction since $\sp{\subgradA,u}=\calJ(u)\neq 0$.
\end{proof}

\begin{remark}[Geometric interpretation]
\cref{thm:ef-equiv-cond} in conjunction with \cref{prop:subdiff_bdry} can be interpreted as relating the geometry of the ambient Hilbert space and the dual unit ball $\dualball$ associated to the functional $\calJ$.
They state that eigenvectors (normalized to lie in boundary of the dual unit ball) are precisely those vectors which admit an orthogonal hyperplane tangent to the dual unit ball.
This is illustrated in \cref{fig:eigenvectors}, showing four different dual unit balls and all eigenvectors (up to normalization).
\begin{figure}[h!]
\centering
\begin{tikzpicture}[scale=.9]
		\node at (0,0) [circle,fill,inner sep=1pt]{};
		\draw (0,0) circle [x radius=2, y radius=1];
		\draw[->,thick](0,0)--(0,1);
        \draw[->,thick](0,0)--(0,-1);
        \draw[->,thick](0,0)--(2,0);
        \draw[->,thick](0,0)--(-2,0);
\end{tikzpicture}
\hfill
\begin{tikzpicture}[scale=.9]
		\node at (0,0) [circle,fill,inner sep=1pt]{};
        \draw (-1.,-1)--(2,-1)--(1,1)--(-2,1)--(-1.,-1);
        \draw[->,thick](0,0)--(0,1);
        \draw[->,thick](0,0)--(0,-1);
        \draw[->,thick](0,0)--(1,1);
        \draw[->,thick](0,0)--(-1,-1);        
        \draw[->,thick](0,0)--(2,-1);
        \draw[->,thick](0,0)--(-2,1);
        \draw[->,thick](0,0)--(1.2,0.6);
        \draw[->,thick](0,0)--(-1.2,-0.6);       
\end{tikzpicture}
\hfill
\begin{tikzpicture}[scale=.9]
		\node at (0,0) [circle,fill,inner sep=1pt]{};
		\draw (1,0.5)--(0.5,1)--(-0.5,1)--(-1,0.5)--(-1,-0.5)
		--(-0.5,-1)--(0.5,-1)--(1,-0.5)--(1,0.5);
		\draw[->,thick](0,0)--(0,1);
		\draw[->,thick](0,0)--(0.5,1);
		\draw[->,thick](0,0)--(1,0);
		\draw[->,thick](0,0)--(1,0.5);		
		\draw[->,thick](0,0)--(0.75,0.75);
		\draw[->,thick](0,0)--(1,-0.5);
		\draw[->,thick](0,0)--(0.5,-1);			
		\draw[->,thick](0,0)--(0.75,-0.75);
		
		\draw[->,thick](0,0)--(0,-1);
		\draw[->,thick](0,0)--(-0.5,-1);
		\draw[->,thick](0,0)--(-1,0);
		\draw[->,thick](0,0)--(-1,-0.5);		
		\draw[->,thick](0,0)--(-0.75,-0.75);
		\draw[->,thick](0,0)--(-1,0.5);
		\draw[->,thick](0,0)--(-0.5,1);			
		\draw[->,thick](0,0)--(-0.75,0.75);
\end{tikzpicture}
\hfill
\begin{tikzpicture}[scale=.9]
		\node at (0,0) [circle,fill,inner sep=1pt]{};
        \draw (-2,-1)--(-0.5,-1)--(2,1)--(0.5,1)--(-2,-1); 
		\draw[->,thick](0,0)--(2,1);	
   		\draw[->,thick](0,0)--(-0.2927,0.3659);      
   		\draw[->,thick](0,0)--(-2,-1);	
   		\draw[->,thick](0,0)--(0.2927,-0.3659);        
\end{tikzpicture}
\caption{Four different dual unit balls $\dualball$ with all eigenvectors.\label{fig:eigenvectors}}
\end{figure}
\end{remark}

\cref{thm:ef-equiv-cond} implies that, if they exist, maximal elements in the dual unit ball are eigenfunctions.
\begin{corollary}
Let $\subgradA\in\argmax\{\norm{\subgradB}\st\subgradB\in\dualball\}$.
Then it holds $\subgradA\in\dJ(\subgradA)$.
\end{corollary}
\begin{proof}
By the maximality of $\subgradA$ it follows
\begin{align*}
    \sp{\subgradA,\subgradB} \leq \norm{\subgradA}\norm{\subgradB} \leq \norm{\subgradA}^2\quad\forall\subgradB\in\dualball
\end{align*}
and hence \cref{thm:ef-equiv-cond} implies the assertion.
\end{proof}

Conversely, we can also show that ground states, which we defined in \cref{def:ground_states} and showed to be eigenfunctions in \cref{thm:opt-cond-nonlin}, are essentially minimal elements on the relative boundary of $\dualball$.
\begin{proposition}
Let $u\in\Hspace$ be a ground state with minimal eigenvalue $\lambda>0$.
Then $\subgradA:=\lambda u\in\dualball$ satisfies $\norm{\subgradA}\leq\norm{\subgradB}$ for all $\subgradB\in\dualball$ for which there exists $v\in\nullspace_\calJ^\perp\setminus\{0\}$ with $\subgradB\in\dJ(v)$. 
\end{proposition}
\begin{proof}
Since $\subgradA\in\dJ(u)$ (and by \cref{prop:subdifferential_homogeneous} also $\subgradA\in\dJ(\subgradA)$) and $\subgradB\in\dJ(v)$, we obtain
\begin{align*}
    \norm{\subgradA} = \frac{\norm{\subgradA}^2}{\norm{\subgradA}} = \frac{\calJ(\subgradA)}{\norm{\subgradA}} = \frac{\calJ(u)}{\norm{u}} \leq \frac{\calJ(v)}{\norm{v}} = \frac{\sp{\subgradB,v}}{\norm{v}}\leq\norm{\subgradB}.
\end{align*}
\end{proof}

We can use \cref{thm:ef-equiv-cond} to characterize eigenfunctions with respect to many different absolutely one-homogeneous functionals. 
Since we know from \cref{prop:ef_minimal_norm} that eigenfunctions necessarily are subgradients of minimal norm in their subdifferential, it suffices to check condition \labelcref{ineq:ef-equiv-cond} for such subgradients.

The first example deals with dual unit balls given by a sufficiently well-behaved polyhedron, see \cref{fig:minsub}.

\begin{proposition}[Non-distorted polyhedron] \label{polyhedraleigenvectors}
Let $K\subset\R^n$ be a closed convex polyhedron satisfying $-K=K$ and let $\calJ(u):=\sup_{\subgradA\in K}\sp{\subgradA,u}$.
Let $u\in\R^n$ and assume that $\subgradA\in\dJ(u)$ satisfies
\begin{align}\label{eq:minsub}
\langle{\subgradA},{\subgradA}-\subgradB\rangle=0,\quad\forall \subgradB\in\dJ(u).
\end{align}
Then ${\subgradA}\in\dJ(\subgradA)$.
\end{proposition}
\begin{proof}
By \labelcref{eq:one-hom_subdiff_dualball} and \cref{prop:subdiff_bdry} we infer that $\dJ(u)$ -- being the intersection of $K$ and the hypersurface $\{\subgradB\in\R^n\st\langle \subgradB,u\rangle=J(u)\}$ -- must coincide with a facet $F$ of the polyhedron.  
Due to \eqref{eq:minsub}, the set $\mathcal{S}:=\{\subgradB\in\R^n\st\langle \subgradB,{\subgradA}\rangle=\norm{{\subgradA}}^2\}$ defines a hypersurface through ${\subgradA}$ such that $F\subset\mathcal{S}$ and $\mathcal{S}$ is orthogonal to ${\subgradA}$. Since $K$ is convex, all other points in $K$ lie on one side of $\mathcal{S}$ which implies that $\mathcal{S}$ is supporting $K$ and hence $\langle{\subgradA},{\subgradA}-\subgradB\rangle\geq0$ for all $\subgradB\in K$. With \cref{thm:ef-equiv-cond} we conclude the assertion.
\end{proof}
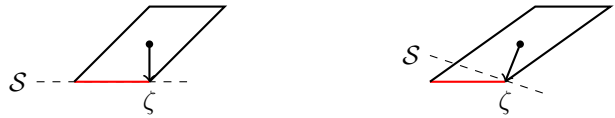
\begin{figure}[h!]
\centering
\begin{tikzpicture}
		\node at (0,0) [circle,fill,inner sep=1pt]{};	
		\draw (-1.5,-0.5) node[anchor = east] {$\mathcal{S}$};
		\draw[dashed] (-1.5,-0.5)--(0.5,-0.5);	
		\draw[red, thick] (-1,-0.5)--(-0,-0.5);
        \draw[thick] (-0,-0.5)--(1,0.5)--(0,0.5)--(-1,-0.5);
        \draw[->,thick] (0,0)--(-0,-0.5);
        \draw (0,-0.5) node[anchor=north] {${\subgradA}$};
\end{tikzpicture}
\hspace*{2cm}
\begin{tikzpicture}
		\node at (0,0) [circle,fill,inner sep=1pt]{};
		\draw (-1.2,-0.15) node[anchor = east] {$\mathcal{S}$};
		\draw[dashed] (-1.2,-0.15)--(0.3,-0.65);		
		\draw[red, thick] (-1.2,-0.5)--(-0.2,-0.5);
        \draw[thick] (-0.2,-0.5)--(1.2,0.5)--(0.2,0.5)--(-1.2,-0.5);
        \draw[->,thick] (0,0)--(-0.2,-0.5);
        \draw (-0.2,-0.5) node[anchor=north] {${\subgradA}$};
\end{tikzpicture}
\caption{\textbf{Left:} \eqref{eq:minsub} is met since $\mathcal{S}$ is supporting, \textbf{Right:} here \eqref{eq:minsub} is violated, i.e., $\mathcal{S}$ is not supporting and ${\subgradA}$ is no eigenvector\label{fig:minsub}}
\end{figure}

The next example concerns the $L^1$-norm, regarded as absolutely one-homogeneous functional on $L^2(\Omega)$.

\begin{proposition}[$L^1$-norm]
Let $\Hspace:=L^2(\Omega)$ and
\begin{align*}
    \calJ(u):=
    \begin{cases}
    \int_\Omega\abs{u}\d x,\quad&\text{if }u\in L^1(\Omega)\cap L^2(\Omega),\\
    +\infty,\quad&\text{else}.
    \end{cases}
\end{align*}
Then it holds $\dualball=\{\subgradB\in L^2(\Omega)\st \norm{\subgradB}_{L^\infty}\leq 1\}$ and for every $u\in\dom(\dJ)$ it holds
$$\dJ(u)=\left\{\subgradA\in \dualball\st\text{$\subgradA=\pm 1$ a.e. where $u\gtrless 0$ and $\subgradA\in[-1,1]$ a.e. where $u=0$}\right\}.$$
Then the subgradient of minimal norm is characterized by
\begin{align*}
\subgradA
&:=
\argmin
\left\lbrace
\norm{\subgradA}_{L^2} \st \subgradA \in \dJ(u)
\right\rbrace
\\
&=
\left\{\subgradA\in \dualball\st\text{$\subgradA=\pm 1$ a.e. where $u\gtrless 0$ and $\subgradA=0$ a.e. where $u=0$}\right\}.
\end{align*}
and satisfies $\sp{\subgradA,\subgradA-\subgradB}_{L^2}\geq0$ for all $\subgradB\in \dualball.$
\end{proposition}
\begin{proof}
We encourage the reader to check the claims about the form of $\dualball$, $\dJ(u)$, and the subgradient of minimal norm.
Since $\subgradA$, being the subgradient of minimal norm and hence taking only the values $0$, $-1$, or $+1$, it satisfies $|\subgradA|=|\subgradA|^2$.
Hence, we calculate with Cauchy-Schwarz
$$\sp{\subgradA,\subgradA-\subgradB}_{L^2}=\int_\Omega|\subgradA|^2\d x-\int_\Omega \subgradA\subgradB\d x\geq\int_\Omega|\subgradA|^2\d x-\int_\Omega |\subgradA|\d x\geq0.$$
\end{proof}

Also for the $L^\infty$-norm we can prove that subgradients of minimal norm are eigenfunctions.

\begin{proposition}[$L^\infty$-norm]
Let $\Hspace:=L^2(\Omega)$ and
\begin{align*}
    \calJ(u):=
    \begin{cases}
    \esssup_\Omega\abs{u}\d x,\quad&\text{if }u\in L^\infty(\Omega)\cap L^2(\Omega),\\
    +\infty,\quad&\text{else}.
    \end{cases}
\end{align*}
Then it holds $\dualball=\{\subgradB\in L^2(\Omega)\st \norm{\subgradB}_{L^1}\leq 1\}$ and for $u\in\dom(\dJ)\setminus\{0\}$ it holds
$$\dJ(u)=\left\lbrace \subgradA\in \dualball\st \norm{\subgradA}_{L^1}=1,\, \operatorname{sign}(\subgradA)=\operatorname{sign}(u)\text{ a.e. in }{\Omega}_{\max} \text{ and } \subgradA=0\text{ a.e. in }\Omega\setminus\Omega_{\max}\right\rbrace,$$
where $\Omega_{\max}:=\{x\in\Omega\st |u(x)|=\calJ(u)\}$ is defined up to a Lebesgue nullset.

Then the subgradient of minimal norm of $u\neq 0$ is characterized by
\begin{align*}
\subgradA
&:=
\argmin
\left\lbrace
\norm{\subgradA}_{L^2} \st \subgradA \in \dJ(u)
\right\rbrace
\\
&=
\left\{\subgradA\in \dualball\st\text{$\subgradA=\operatorname{sign}(u)/\abs{\Omega_{\max}}$ a.e. in $\Omega_{\max}$ and $\subgradA=0$ a.e. in $\Omega\setminus\Omega_{\max}$}\right\}
\end{align*}
and satisfies $\sp{\subgradA,\subgradA-\subgradB}_{L^2}\geq0$ for all $\subgradB\in \dualball.$
\end{proposition}
\begin{proof}
Again, we encourage the reader to check the claims about the form of $\dualball$ and $\dJ(u)$.
W.l.o.g. we can assume that $u\neq 0$ which implies
$$1=\norm{\subgradB}_{L^1}^2\leq|\Omega_{\max}|\norm{\subgradB}_{L^2}^2,\quad\forall \subgradB\in\dJ(u).$$
Furthermore, the tentative subgradient of minimal norm satisfies $\norm{\subgradA}_{L^2}^2=1/|\Omega_{\max}|$ which implies $\norm{\subgradA}_{L^2}\leq\norm{\eta}_{L^2}$. 
By uniqueness we infer that $\subgradA$ is the subgradient of minimal norm.

To show that it is an eigenfunction, we compute
\begin{align*}
\sp{\subgradA,\subgradA-\subgradB}_{L^2}=\norm{\subgradA}_{L^2(\Omega)}^2-\frac{1}{|\Omega_{\max}|}\int_{\Omega_{\max}}\operatorname{sign}(u)\subgradB\dx\geq \frac{1}{|\Omega_{\max}|}-\frac{1}{|\Omega_{\max}|}\int_{\Omega_{\max}}|\subgradB|\dx\geq 0.
\end{align*}
\end{proof}

In our last example we consider the functional
\begin{align}\label{eq:divergence}
\calJ(u):=\int_\Omega|\div u|:=\sup\left\lbrace-\int_\Omega u\cdot\nabla v\dx\st v\in C^\infty_c(\R^n;[-1,1])\right\rbrace,\;u\in L^2(\Omega,\R^n).
\end{align}
It holds that $u\in\dom(J)$, i.e., $\calJ(u)<\infty$, if and only if the distribution $\div u$ can be represented as a finite Radon measure (cf.~\cite{briani2011gradient}). 

The dual unit ball can be shown to be given by
\begin{align}\label{eq:set_K_div}
    \dualball:=\left\lbrace-\nabla v\st v\in H^1_0(\Omega;[-1,1])\right\rbrace.
\end{align}
Hence, a subgradient $\subgradA$ of $\calJ$ in $u$ fulfills $\subgradA\in\dualball$ and $\int_\Omega \subgradA\cdot u\d x=\int_\Omega|\div u|\d x$.
To understand the meaning of $v$ in \eqref{eq:set_K_div}, we perform a integration by parts for smooth $u$ to obtain
$$\int_\Omega \subgradA\cdot u\d x=-\int_\Omega\nabla v\cdot u\d x=\int_\Omega v\,\div u\d x.$$
Therefore, $v$ should be chosen as 
\begin{align}\label{eq:v_smooth}
v(x)\in 
\begin{cases}
\{1\},\quad&\div u(x)>0,\\
\{-1\},\quad&\div u(x)<0,\\
[-1,1],\quad&\div u(x)=0,
\end{cases}\quad x\in\Omega.
\end{align}
For more irregular $u$, one considers the polar decomposition of the measure $\div u$, given by $\div u=\theta_u\,|\div u|,$ where $\theta_u(x):=\frac{\div u}{|\div u|}(x)$ denotes the Radon-Nikodym derivative, which exists $|\div u|$-almost everywhere. This allows us to define the sets
\begin{align}
\E_u^\pm:=\{x\in\Omega\st\theta_u(x)\text{ exists and }\theta_u(x)=\pm 1\}.
\end{align}
Note that trivially it holds $\E_u^+\cup\E_u^-\subset\supp(\div u)$  and one can easily prove that $\supp(\div u)=\overline{\E_u^+\cup\E_u^-}$. In \cite{briani2011gradient} the authors showed that, in full analogy to \eqref{eq:v_smooth}, the subdifferential of $\calJ$ in $u\in\dom(\dJ)$ can be characterized as
\begin{align}\label{eq:subdiff}
\dJ(u)=\left\{-\nabla v\st v\in H^1_0(\Omega;[-1,1]),\,v=\pm1\;|\div u|\text{-a.e. on }\E_u^\pm\right\}.
\end{align}
If $\subgradA=-\nabla v$ is a subgradient of $\calJ$ in $u$ we refer to $v$ as a \emph{calibration} of $u$.

\begin{remark}
Since $v\in H^1_0(\Omega,[-1,1])$ can be defined pointwise everywhere except from a set with zero $H^1$-capacity and the measure $\div u$ vanishes on such sets (cf.\cite[Ch.~6]{adams1999function} and \cite{chen2009gauss}, respectively), the pointwise definition of $v$ in~\eqref{eq:subdiff} makes sense.
\end{remark}
Using \eqref{eq:subdiff} the subgradient of minimal norm is defined by
\begin{align}\label{eq:minimal_subgradient}
\subgradA:=-\nabla\argmin\left\{\int_\Omega|\nabla v|^2\d x\st v\in H^1_0(\Omega;[-1,1]),\,v=\pm1\;|\div u|\text{-a.e. on }\E_u^\pm\right\},
\end{align}
which implies that $\subgradA=-\nabla v$ where $v$ minimizes the Dirichlet energy and, in particular, is harmonic on the open set $\Omega_H:=\Omega\setminus(\overline{\E_u^+\cup\E_u^-})$. 

\begin{proposition}[Divergence semi-norm]\label{prop:eigenvectors_div}
Let $\calJ:L^2(\Omega;\R^n)\to\RI$ be defined by \labelcref{eq:divergence}.
Let $u\in\dom(\dJ)$ and $\subgradA\in\dJ(u)$ be the subgradient of minimal norm, given by \labelcref{eq:minimal_subgradient}.
Then it holds $\sp{\subgradA,\subgradA-\subgradB}_{L^2}\geq 0$ for all $\subgradB\in\dualball$.
\end{proposition}

To prove the proposition, we start with an approximation lemma. 

\begin{lemma}\label{lem:approx}
Let $0<\eps<1$ and define
\begin{align}
\psi_\eps(t)=
\begin{cases}
0,&|t|>1,\\
1,&|t|<1-\eps,\\
\frac{1}{\eps}(1-|t|),&1-\eps\leq|t|\leq1,
\end{cases}
\qquad\phi_\eps(t)=\int_0^t\psi_\eps(\tau)\d\tau.
\end{align}
If $v\in H^1_0(\Omega,[-1,1])$ then $v_\eps:=\phi_\eps\circ v\in H^1_0(\Omega,[-1,1])$ and $v_\eps\to v$ strongly in $H^1$ as $\eps\searrow0$. 
\end{lemma}
\begin{proof}
The membership to $H^1_0(\Omega,[-1,1])$ follows directly from the chain rule for compositions of Lipschitz and Sobolev functions. For strong convergence it suffices to show $\nabla v_\eps\to\nabla v$ in $L^2(\Omega)$. Using $\nabla v_\eps=\psi_\eps(v)\nabla v$ this follows from
\begin{align*}
\norm{\nabla v_\eps-\nabla v}_{L^2}^2=&\norm{(\psi_\eps(v)-1)\nabla v}_{L^2}^2=\int_{\{1-\eps\leq|v|\leq 1\}}\underbrace{\left[\frac{1-\eps-|v(x)|}{\eps}\right]^2}_{\leq1}|\nabla v|^2\dx\\
\leq&\int_{\{1-\eps\leq|v|\leq1\}}|\nabla v|^2\dx\to\int_{\{|v|=1\}}|\nabla v|^2\dx=0,\quad\eps\searrow0.
\end{align*}
\end{proof}

\begin{proof}[Proof of Proposition \ref{prop:eigenvectors_div}]
We write $\subgradA=-\nabla v$ with $v$ as in \eqref{eq:minimal_subgradient} and define $v_\eps$ as in \cref{lem:approx}. 
Any $\subgradB\in \dualball$ can be written as $\subgradB=-\nabla w$ with $w\in H^1_0(\Omega,[-1,1])$ and thanks to \cref{lem:approx} we have
\begin{align}\label{eq:proof1}
\sp{\subgradA,\subgradA-\subgradB}_{L^2}=\int_\Omega\nabla v\cdot\nabla(v-w)\d x=\lim_{\eps\searrow0}\int_\Omega\nabla v_\eps\cdot\nabla(v-w)\d x.
\end{align}
As a first step, we replace $v-w$ by a smooth function $z\in C^\infty_c(\Omega)$ with compact support, to obtain with the product rule:
\begin{align*}
\int_\Omega\nabla v_\eps\cdot\nabla z\dx=&\int_\Omega\psi_\eps(v)\nabla v\cdot\nabla z\dx\\
=&\int_\Omega\nabla v\cdot(\nabla[\psi_\eps(v)z]-\psi_\eps'(v)z\nabla v)\dx\\
=&-\int_\Omega\psi_\eps'(v)z|\nabla v|^2\dx
\end{align*}
For the last equality we used that $\psi_\eps(v)z$ is a test function for the minimization of the Dirichlet energy in \eqref{eq:minimal_subgradient} and that the first variation of the Dirichlet integral in direction of the test function vanishes consequently.

Strongly approximating $v-w$ by smooth and compactly supported functions in $H^1_0(\Omega)$ and using above-noted calculation, results in
\begin{align}\label{eq:proof2}
\int_\Omega\nabla v_\eps\cdot\nabla(v-w)\dx=-\int_\Omega\psi_\eps'(v)(v-w)|\nabla v|^2\dx.
\end{align}
Note that
$$\psi_\eps'(t)=\frac{1}{\eps}
\begin{cases}
1,&t\in[-1,-1+\eps]\\
-1,&t\in[1-\eps,1],\\
0,&$else$,
\end{cases}
$$
which we can use to calculate
\begin{align}\label{eq:proof3}
-&\int_\Omega\psi_\eps'(v)(v-w)|\nabla v|^2\dx\notag\\
=&\int_{\{-1\leq v\leq-1+\eps\}}\frac{1}{\eps}(w-v)|\nabla v|^2\dx+\int_{\{1-\eps\leq v\leq1\}}\frac{1}{\eps}(v-w)|\nabla v|^2\dx\notag\\
\geq&-\int_{\{1-\eps\leq|v|\leq 1\}}|\nabla v|^2\dx\notag\\
\to&-\int_{\{|v|=1\}}|\nabla v|^2\dx=0,\quad\eps\searrow0.
\end{align}
Putting \labelcref{eq:proof1,eq:proof2,eq:proof3} together, we have shown that 
$$\sp{\subgradA,\subgradA-\subgradB}_{L^2}=\lim_{\eps\searrow 0}\int_\Omega\nabla v_\eps\cdot\nabla(v-w)\dx =\lim_{\eps\searrow 0}-\int_\Omega\psi_\eps'(v)(v-w)|\nabla v|^2\d x\geq 0,$$
as desired.
\end{proof}

\section{Gradient flows of absolutely \texorpdfstring{$p$}{p}-homogeneous functionals}

We consider the following (sub-)gradient flow.
\begin{subequations}\label{eq:grad-flow}
\begin{align}[left=\empheqlbrace]
    & u'(t) + \dJ(u(t)) \ni 0, \quad t>0, \\
    & u(0) = f,
\end{align}
\end{subequations}
where $f \in \Hspace$ is the initial datum.

The next theorem due to Brezis ensures the existence of a solution to this gradient flow. We define the following single-valued operator 
\begin{align}\label{eq:subgrad_minimal_norm}
u \mapsto \partial^0 \calJ(u) \defeq \argmin\left\lbrace\norm{\subgradA}\st\subgradA\in\dJ(u)\right\rbrace,\quad u \in \dom(\dJ),
\end{align}
which associates each $u$ with the subgradient with minimal norm in $\dJ(u)$. This operator is single-valued since $\dJ(u)$ is a closed convex set.

\begin{theorem}[Brezis \cite{brezis1973ope}]\label{thm:brezis}
Let $\calJ:\H\to\R\cup\{\infty\}$ be proper, convex, and lower semicontinuous and let $f\in\overline{\dom(\calJ)}$.  Then there exists exactly one continuous map $u:[0,\infty)\to\H$---called $t\mapsto u(t)$ the gradient flow of $\calJ$ with initial datum $f$---which is Lipschitz continuous on $[\delta,\infty),\;\delta>0$ and right-differentiable on $(0,\infty)$ such that 
\begin{itemize}
\item $u(0)=f$,
\item $\subgradA(t):=-\partial_t^+ u(t)=\partial^0 \calJ(u(t))$ for all $t>0$,
\item \labelcref{eq:grad-flow} holds for almost every $t>0$,
\item $t\mapsto\subgradA(t)$ is right continuous for all $t>0$ and $t\mapsto\norm{\subgradA(t)}$ is non-increasing,
\item $t\mapsto \calJ(u(t))$ is convex, non-increasing and Lipschitz continuous on $[\delta,\infty),\;\delta>0$ with
\begin{align}\label{eq:derivative_J}
\frac{\d^+}{\d t}\calJ(u(t))=-\norm{\subgradA(t)}^2,\quad t>0.
\end{align}
\end{itemize}
Here $\d^+/\d t$ and $\partial_t^+$ denote right-derivatives and will be replaced by standard derivative symbols in the sequel.
\end{theorem}

We will need the following result due to Brezis that provides an integrability condition of the time derivative at zero.

\begin{proposition}[Brezis, Proposition 3 in  \cite{brezis1971proprietes}]\label{prop:finite_energy_integrable}
Under the conditions of \cref{thm:brezis} it holds that
\begin{align*}
f\in\dom(\calJ) 
\iff
t\mapsto\subgradA(t) \in L^2((0,\delta),\Hspace)\;\quad\forall\delta>0.
\end{align*}
\end{proposition}

\begin{corollary}\label{cor:integration-of-grad-flow}
If the initial datum satisfies $\calJ(f)<\infty$, the following representation (obtained by formally integrating the gradient flow~\eqref{eq:grad-flow}) is rigorous:
\begin{align}\label{eq:int-grad-flow}
    f- u(t) = -\int_{0}^t \partial_s u(s)\d s = \int_0^t \subgradA(s)\d s, \quad t>0.
\end{align}
\end{corollary}

Another useful statement is the following invariance property of the gradient flow.
Later we will use that with $C$ being the nullspace $\nullspace_\calJ$ or the non-negative orthant.

\begin{proposition}[{\cite[Thm. 1.1]{barthelemy1996invariance}, see also \cite[Prop. 4.5]{brezis1973ope}}]
Let $t\mapsto u(t)$ denote the gradient flow of $\calJ$ with initial datum $f\in\closure{\dom(\calJ)}$ and let $C\subset\Hspace$ be a convex set with associated orthogonal projection operator $\Pi_C$.
If $\calJ(\Pi_C u) \leq \calJ(u)$ and $f\in C$ then $u(t)\in C$ for all $t>0$.
\end{proposition}

\begin{proposition}[Conservation of mass]\label{prop:conserv_mass}
Let $u(t)$ solve the gradient flow \eqref{eq:grad-flow} with data $f\in\dom(\calJ)$ and let $P_{\nullspace_\calJ}:\Hspace\to\nullspace_\calJ$ denote the orthogonal projection onto $\nullspace_\calJ$. Then it holds $P_{\nullspace_\calJ}{u(t)}=P_{\nullspace_\calJ}{f}$.
\end{proposition}
\begin{proof}
Using~\eqref{eq:int-grad-flow} and \cref{prop:orth}, we deduce $P_{\nullspace_\calJ}{\subgradA(s)}=0$ for all $s>0$. Hence, also $P_{\nullspace_\calJ}(u(t)-f)=0$ which implies the statement due to linearity of the projection.
\end{proof}

\begin{proposition}[{\cite[Thm. 5]{bruck:1975}}]\label{prop:properties_GF}
Let $u(t)$ denote the solution of the gradient flow \eqref{eq:grad-flow} corresponding to the absolutely $p$-homogeneous, convex, and \lsc{} functional $\calJ:\Hspace\to\RI$ and let $f\in\overline{\dom(\calJ)}$. 
Then it holds
\begin{alignat}{2}
\label{conv:u-to-0}
u(t)\to P_{\nullspace_\calJ} f ,\quad&\calJ(u(t))\to 0,\quad&&t\to\infty,\\
\label{eq:diss_u}
\frac{\d }{\d t}\frac{1}{2}\norm{u(t)-P_{\nullspace_\calJ} f}^2&=-p\calJ(u(t)),\quad&&t>0.
\end{alignat}
\end{proposition}
\begin{proof}
Fix $t_0>0$ and let $g \colon [0,t_0] \to \R$ be defined as
\begin{equation*}
    g(t) \defeq \norm{u(t)}^2 - \norm{u(t_0)}^2 - \frac12 \norm{u(t)-u(t_0)}^2.
\end{equation*}
By \cref{thm:brezis}, $u(t)$ is Lipschitz continuous on $[\delta,\infty)$ for any $\delta > 0$. Thus, $g(t)$ is also Lipschitz on $[\delta,\infty)$ and therefore differentiable a.e. (Rademacher's theorem). Hence we get $g'(t) =  \sp{u'(t), u(t) +u(t_0)}$ a.e.

Since $\calJ(u(t))$ is non-decreasing (\cref{thm:brezis}) and $-u'(t) \in \dJ(u(t))$ a.e., we have
\begin{align*}
    \calJ(u(t)) &\geq \calJ(u(t_0)) = \calJ(-u(t_0))  \\
        &\geq \calJ(u(t)) + \sp{-u'(t),-u(t_0)-u(t)} \\
        &= \calJ(u(t)) + g'(t)    
\end{align*}
a.e. in $(0,t_0]$. Therefore, $g'(t)\leq0$ a.e. on in $(0,t_0]$ and $g$ is non-increasing. Hence, $0 = g(t_0) \leq g(t)$ for all $0 < t \leq t_0$ and
\begin{equation*}
    \norm{u(t)-u(t_0)}^2 \leq 2 (\norm{u(t)}^2 - \norm{u(t_0)}^2), \quad 0 < t \leq t_0.
\end{equation*}
Therefore, $\norm{u(t)}$ is non-increasing on $(0,\infty)$ and therefore converges. This implies that the  right-hand side converges to $0$ as $t,t_0 \to \infty$ and so does $\norm{u(t)-u(t_0)}$.  Therefore, $\{u(t)\}_{t>0}$ is Cauchy and converges strongly to some $u_\infty \in \Hspace$. 

For any $v_0\in\nullspace_\calJ$ we have
\begin{align*}
    \calJ(u(t)) + \sp{-u'(t),v_0-u(t)} \leq \calJ(v_0) \leq \calJ(u(t))
\end{align*}
and hence
\begin{align*}
    h(t) := \sp{-u'(t),u(t)-v_0} \geq 0.
\end{align*}
It also holds
\begin{align*}
    \frac{\d}{\d t}\frac{1}{2}\norm{u(t)-v_0}^2 = -h(t) \leq 0
\end{align*}
and hence the limit $\lim_{t\to\infty}\frac12\norm{u(t)-v_0}^2$ exists.
Therefore $h\in L^1((0,\infty))$ and there exists a sequence $(t_k)_{k\in\N}$ which satisfies $t_k\to\infty$ and $h(t_k) \to 0$.

Using the lower semicontinuity of $\calJ$ and the uniqueness of the limit of $u(t)$ we obtain from
\begin{align*}
    \calJ(u(t_k)) +  \sp{-u'(t_k),v_0-u(t_k)} \leq \calJ(v_0) 
\end{align*}
that
\begin{align*}
    \calJ(u_\infty) \leq \calJ(v_0)
\end{align*}
and hence $u_\infty\in\nullspace_\calJ$.
By \cref{prop:conserv_mass} it follows $u_\infty = P_{\nullspace_\calJ} f$.

For \eqref{eq:diss_u} one uses the chain rule together with \eqref{eq:euler} to obtain
$$\frac{\d}{\d t}\frac{1}{2}\norm{u(t)-P_\nullspace f}^2=\left\langle\partial_tu(t),u(t)-P_\nullspace f\right\rangle=-\langle\subgradA(t),u(t)-P_\nullspace f\rangle=-p \calJ(u(t)).$$
Hence, we get
\begin{align*}
    p\calJ(u(t)) = -\frac{\d}{\d t}\frac12\norm{u(t)-P_{\nullspace_\calJ}f}^2 = h(t)
\end{align*}
and a subsequence converges to zero. Since $t\mapsto\calJ(u(t))$ is non-increasing in fact the whole sequence converges to zero.
\end{proof}

We now show that the gradient flow~\eqref{eq:grad-flow} allows one to decompose any datum $f \in \Hspace$ using subgadients of the functional $\calJ$.

This immediately gives us the following result.
\begin{theorem}
Let $u(t)$ denote the solution of the gradient flow \eqref{eq:grad-flow} corresponding to the absolutely $p$-homogeneous, convex, and \lsc{} functional $\calJ:\Hspace\to\RI$ and let $f\in{\dom(\calJ)}$.  
Then the following decomposition holds
\begin{equation}\label{eq:ef-decomp}
    f =  P_{\nullspace_\calJ} f + \int_0^\infty \zeta(t) \d t.
\end{equation}
\end{theorem}

This decomposition looks very similar to the linear decomposition in a basis of \ef{s}, however, to justify the analogy, we need to make sure that $\zeta(t)$ are indeed \ef{s} according to \cref{def:ef-nonlin}. We have a complete answer in the $1$-homogeneous case.
\begin{theorem}\label{thm:flow-decomp-efs}
Let $\calJ$ be proper, convex, \lsc{} and absolutely $1$-homogeneous. Then the decomposition~\eqref{eq:ef-decomp} is a decomposition into \ef{s} if and only if for any $t>0$ it holds that
\begin{equation}\label{eq:ef-cond-all-t}
    \sp{\subgradA(t),\subgradA(t)-\subgradB} \geq 0 \quad \forall \, \subgradB \in \dualball,
\end{equation}
where $\dualball$ is the dual unit ball of $\calJ$.

Sufficient for this is that every subgradient of minimal norm $\subgradA$ satisfies $\sp{\subgradA,\subgradA-\subgradB}\geq 0$ for all $\subgradB\in\dualball$.
\end{theorem}
\begin{proof}
Cf. \cref{thm:ef-equiv-cond}.
\end{proof}

We will need the following  result.
\begin{proposition}\label{prop:subgradients-accumulate}
Let $\calJ$ be proper, convex, \lsc{} and absolutely $1$-homogeneous and let~\eqref{eq:ef-cond-all-t} be satisfied. Then for all $0<s \leq t$ it holds that $\subgradA(s) \in \dJ(u(t)))$.
\end{proposition}
\begin{proof}
Since $\subgradA(t)$ is an \ef{} by \cref{thm:flow-decomp-efs}, we have that
\begin{equation*}
    \sp{\subgradA(t),\subgradA(t)-\subgradA(s)} \geq 0
\end{equation*}
for all $0 < s \leq t$. Therefore, using \cref{thm:brezis} we get
\begin{equation*}
    0 \geq \sp{-\subgradA(t),\subgradA(t)-\subgradA(s)} = \frac{\d}{\d t} (\calJ(u(t))) - \sp{u'(t),\subgradA(s)}.
\end{equation*}
Integrating this from $s$ to $t$ yields
\begin{equation*}
    \calJ(u(t)) - \calJ(u(s)) - \sp{u(t),\subgradA(s)} + \calJ(u(s)) \leq 0,
\end{equation*}
which is equivalent to $\calJ(u(t)) \leq \sp{u(t),\subgradA(s)}$. Since $\subgradA(s)$ is a subgradient, we have that also $\calJ(u(t)) \geq \sp{u(t),\subgradA(s)}$, equality holds, and $\subgradA(s) \in \dJ(u(t))$.
\end{proof}

If the condition of \cref{thm:flow-decomp-efs} is satisfied then the \ef{s} $\zeta(t)$ have the following orthogonality property.
\begin{proposition}
Let $\calJ$ be proper, convex, \lsc{} and absolutely $1$-homogeneous and let~\eqref{eq:ef-cond-all-t} be satisfied. Then for any $0<r\leq s\leq t$ we have that 
\begin{align}\label{eq:orthogonality}
\sp{\subgradA(t),\subgradA(s)-\subgradA(r)}=0.
\end{align}
\end{proposition}
\begin{proof}
By \cref{prop:subgradients-accumulate}, it holds that $\subgradA(r),\subgradA(s),\subgradA(t)\in\dJ(u(t))$. Furthermore, $\norm{\subgradA(t)}$ is minimal in $\dJ(u(t))$. Hence, the assertion follows directly by employing \cref{prop:minsub} with $u:=u(t)$, $\subgradA:=\subgradA(t)$, and $\subgradB \in \{\subgradA(s),\subgradA(r)\}$ to obtain $\sp{\subgradA(t),\subgradA(s)}  = \norm{\subgradA(t)}^2 = \sp{\subgradA(t),\subgradA(r)}$.
\end{proof}

In case of a spectral decomposition into eigenvectors we can also prove a strong short time regularity of the gradient flow.
\begin{proposition}
Let $f\in\dom(\calJ)$ and \eqref{eq:ef-cond-all-t} be satisfied. Then it holds $\calJ(f-u(t))\to 0$ as $t\searrow 0$.
\end{proposition} 
\begin{proof}
Using the definition of the gradient flow, we write $f-u(t)=\int_0^t \subgradA(s)\d s$. Since all the subgradients $\subgradA(s)$ for $s>0$ are eigenvectors, this implies together with Jensen's inequality that
$$0\leq \calJ(f-u(t))\leq\int_0^t\calJ(\subgradA(s))\d s=\int_0^t\norm{\subgradA(s)}^2\d s.$$
According to \cref{prop:finite_energy_integrable} the map $s\mapsto\norm{\subgradA(s)}$ is in $L^2(0,\delta)$ for all $\delta>0$ since $f\in\dom(\calJ)$. Hence, the dominated convergence theorem yields $\lim_{t\searrow 0}\calJ(f-u(t))=0$, as desired.
\end{proof}

\section{Asymptotic profiles}\label{sec:asymptotic_profiles}

In this section we continue the study of the gradient flow
\begin{align}
    \begin{cases}
    u'(t) = -\subgradA(t),\quad \subgradA(t)\in\dJ(u(t)),\quad t>0,\\
    u(0) = f,
    \end{cases}
\end{align}
of an absolutely $p$-homogeneous, proper, convex, and \lsc{} functional $\calJ:\Hspace\to\RI$.

\subsection{First-order: extinction times and exact convergence rates}
\label{ss:extinction_times}

From \cref{prop:properties_GF} we know that the gradient flow converges to $u_\infty := P_\nullspace f$, i.e., the orthogonal projection of the initial data on the nullspace of $\calJ$.

We would like to understand how quick the convergence takes place.
To this end, we define the \emph{extinction time} of the gradient flow as
\begin{align}\label{eq:extinction_time}
    \tex(f) \defeq \inf\{t>0\st u(t)=u_\infty\}.
\end{align}
As it turns out, the value of the extinction time---if it is finite---depends on the minimal eigenvalue $\lambda_1$ of $\dJ$, see \labelcref{eq:Rayleigh-opt-prob}. 
Also the speed of convergence of $u(t)$ to its limit $u_\infty := P_\nullspace f$ directly depends on this eigenvalue.
An important condition will be that this eigenvalue is positive. 
Taking \labelcref{eq:Rayleigh-opt-prob} into account this is equivalent to the following:
\begin{align}\label{eq:poincare}
    \exists C>0 \st \norm{u-P_\nullspace u}^p \leq C \calJ(u),\quad\forall u\in\Hspace.
\end{align}
Obviously, the minimal constant such that \labelcref{eq:poincare} is satisfied, is given by $C := \frac{p}{\lambda_1}$.
In the sequel, we will refer to \labelcref{eq:poincare} as \emph{Poincar\'e-type inequality}.

\begin{example}[Poincar\'e inequality]
A classical example of \labelcref{eq:poincare} is the Poincar\'e(--Wirtinger) inequality
\begin{align}
    \norm{u - u_\Omega}_{L^2(\Omega)}^2 \leq C \norm{\nabla u}_{L^2(\Omega)}^2,\quad\forall u \in W^{1,2}(\Omega),
\end{align}
where $u_\Omega := \frac{1}{\abs{\Omega}}\int_\Omega u \d x$ and $\Omega\subset\R^n$ is a Lipschitz domain.
\end{example}

\begin{example}[Sobolev inequality]
More generally, Sobolev inequalities of the form
\begin{align}
    \norm{u - u_\Omega}_{L^2(\Omega)}^2 \leq 
    \begin{cases}
    C \tv(u) ,\quad&\forall u \in BV(\Omega), \; n \in \{1,2\}\\
    C \norm{\nabla u}_{L^p(\Omega)}^p,\quad&\forall u \in W^{1,p}(\Omega), \; p \geq \frac{2n}{n+2}.
    \end{cases}
\end{align}
\end{example}

An important quantity in the following will be the following Rayleigh quotient of the gradient flow:
\begin{align}\label{eq:Lambda}
    \Lambda(t) := \frac{p\calJ(u(t))}{\norm{u(t)-u_\infty}^p}.
\end{align}
Note that since $u_\infty = P_\nullspace f = P_\nullspace u(t)$ it holds $u(t)-u_\infty \in \nullspace_\calJ^\perp$.
Furthermore, by \cref{thm:brezis} we know that $t\mapsto\Lambda(t)$ is right differentiable for all $0<t<\tex$. 
We denote the right derivative by $\Lambda'(t)$.

\begin{proposition}\label{prop:Lambda_decreasing}
It holds that $\Lambda'(t) \leq 0$ for all $0<t<\tex$.
\end{proposition}
\begin{proof}
Using \cref{thm:brezis} and the quotient rule, we compute
\begin{align*}
\Lambda'(t)&=-p^2J(u(t))\norm{u(t)-u_\infty}^{-p-1}\frac{\d}{\d t}\norm{u(t)-u_\infty}-p\norm{u(t)-u_\infty}^{-p}\norm{\subgradA(t)}^2\\
&=-p^2J(u(t))\norm{u(t)-u_\infty}^{-p-2}\frac{\d}{\d t}\frac{1}{2}\norm{u(t)-u_\infty}^2-p\norm{u(t)-u_\infty}^{-p}\norm{\subgradA(t)}^2\\
&=p\norm{u(t)-u_\infty}^{-p}\left[\frac{\langle\subgradA(t),u(t)-u_\infty\rangle^2}{\norm{u(t)-u_\infty}^2}-\norm{\subgradA(t)}^2\right]\leq 0,
\end{align*}
\end{proof}

Next we study, how quick $u(t)$ approaches $u_\infty$.
For this we first compute the time derivative of suitable powers of $\norm{u(t)-u_\infty}$.
\begin{proposition}\label{prop:dissipation}
It holds that
\begin{align}\label{eq:diss_u_Lambda}
\begin{cases}
\frac{1}{2-p}\frac{\d}{\d t}\norm{u(t)-u_\infty}^{2-p}&=-\Lambda(t),\phantom{\norm{u(t)-u_\infty}^2}\quad\; p\neq 2,\\[10pt]
\hspace*{10pt}\frac{1}{2}\frac{\d}{\d t}\norm{u(t)-u_\infty}^{2}&=-\Lambda(t)\norm{u(t)-u_\infty}^2,\quad p=2.
\end{cases}
\end{align}
\end{proposition}
\begin{proof}
Using \eqref{eq:diss_u} and the definition of $\Lambda$ in \eqref{eq:Lambda} we can compute
$$\frac{\d }{\d t}\frac{1}{2}\norm{u(t)-u_\infty}^2=-pJ(u(t))=-\Lambda(t)\norm{u(t)-u_\infty}^p.$$
In the case $p=2$ this already proves \eqref{eq:diss_u_Lambda}. For $p>2$ we can use this equality and calculate
\begin{align*}
    \frac{1}{2-p}\frac{\d}{\d t}\norm{u(t)-u_\infty}^{2-p}
    &=
    \norm{u(t)-u_\infty}^{1-p}\frac{\d}{\d t}\norm{u(t)-u_\infty}
    \\
    &=
    \norm{u(t)-u_\infty}^{-p}\frac{\d}{\d t}\frac{1}{2}\norm{u(t)-u_\infty}^2
    =
    -\Lambda(t)
\end{align*}
which concludes the proof also in this case.
\end{proof}

Using \cref{prop:dissipation} we can prove upper bounds for the convergence rate to $u(t)$ to~$u_\infty$.
\begin{proposition}[Upper bounds]\label{prop:upper_bounds}
It holds that
\begin{subequations}\label{ineq:decay_u_gen}
\begin{alignat}{2}
\norm{u(t)-u_\infty}^{2-p}&\leq\norm{f-u_\infty}^{2-p}-(2-p)\lambda_1 t,\qquad &&p<2,\label{ineq:decay_u_p<2}\\
\norm{u(t)-u_\infty}^{2\phantom{-p}}&\leq\norm{f-u_\infty}^2\exp(-2\lambda_1t),\qquad &&p=2,\label{ineq:decay_u_p=2}\\
\norm{u(t)-u_\infty}^{p-2}&\leq\frac{1}{\norm{f-u_\infty}^{2-p}+(p-2)\lambda_1 t},\qquad &&p>2.\label{ineq:decay_u_p>2}
\end{alignat}
\end{subequations}
\end{proposition}
\begin{proof}
For the proof one uses $\Lambda(t)\geq\lambda_1$ in \eqref{eq:diss_u_Lambda} and integrates the resulting inequality from $0$ to $t$ (using Gronwall's Lemma for $p=2$).

For illustration purposes we prove the case $p<2$.
From \cref{prop:dissipation} we get
\begin{align*}
    \frac{1}{2-p} \norm{u(t)-u_\infty}^{2-p} - \frac{1}{2-p}\norm{f-u_\infty}^{2-p} \leq - \lambda_1.
\end{align*}
Since $p<2$ we can multiply with $2-p$ and reorder to get \labelcref{ineq:decay_u_p<2}.
\end{proof}

A first immediate consequence of these upper bounds is that the extinction time $\tex$ in \labelcref{eq:extinction_time} is finite if $p<2$ and the first eigenvalue $\lambda_1$ is positive (which is equivalent to the validity of the Poincar\'e-type inequality \labelcref{eq:poincare}.

\begin{corollary}
If $p<2$ and $\lambda_1>0$, then it holds that
\begin{align}
    \tex \leq \frac{\norm{f-u_\infty}^{2-p}}{(2-p)\lambda_1} < \infty.
\end{align}
\end{corollary}

Interestingly, we can also prove a lower bound for the extinction time, however, only in the absolutely one-homogeneous case.

\begin{proposition}
If $\calJ:\Hspace\to\RI$ is absolutely one-homogeneous, it holds
\begin{align}
    \tex \geq \calJ_*(f-u_\infty),
\end{align}
where $\calJ_*$ denotes the dual semi-norm defined in \labelcref{eq:dual_seminorm}. 
\end{proposition}
\begin{proof}
Using $u(\delta)-u_\infty = \int_\delta^\infty \subgradA(t)\d t$ for $\delta>0$ it holds
\begin{align*}
    \calJ_*(f-u_\infty) 
    &= 
    \sup_{u\in\Hspace}\lim_{\delta\downarrow 0}\frac{\sp{u(\delta)-u_\infty,u}}{\calJ(u)}
    =
    \sup_{u\in\Hspace}\lim_{\delta\downarrow 0}\frac{1}{\calJ(u)} \int_\delta^{\tex} \sp{\subgradA(t),u}\d t 
    \\
    &\leq 
    \lim_{\delta\downarrow 0}(\tex-\delta) = \tex.
\end{align*}
\end{proof}

One might wonder when the upper and lower bounds of the extinction time coincide. 
As it turns out this is the case if the initial datum is a ground state.

\begin{proposition}
Let $\calJ:\Hspace\to\RI$ be absolutely one-homogeneous and $\lambda_1>0$.
It holds that $\calJ_*(f-u_\infty)=\frac{\norm{f-u_\infty}}{\lambda_1}$ if $\lambda_1 \frac{f-u_\infty}{\norm{f-u_\infty}} \in \dJ(f)$.
\end{proposition}
\begin{proof}
W.l.o.g. we can assume that $u_\infty = P_\nullspace f = 0$.
It holds
\begin{align}\label{ineq:dual_seminorm_poincare}
    \calJ_*(f) = \sup_{u\in\H}\frac{\sp{f,u}}{\calJ(u)} \leq {\norm{f}}\sup_{u\in\H}\frac{\norm{u}}{\calJ(u)} \leq \frac{\norm{f}}{\lambda_1}.
\end{align}
If $f$ is a ground state, we also have
\begin{align*}
    \calJ_*(f) \geq \frac{\sp{f,f}}{\calJ(f)} = \frac{\norm{f}^2}{\lambda_1\norm{f}}=\frac{\norm{f}}{\lambda_1}.
\end{align*}
This shows the equality.
\end{proof}

For homogeneities $p\geq 2$ the extinction time is always zero or infinite.
This is a simple consequence of the following proposition which, analogously to \cref{prop:upper_bounds}, constructs lower bounds for the convergence rate of $u(t)$ to $u_\infty$.
\begin{proposition}[Lower bounds]\label{prop:lower_bounds}
For all $\delta>0$ and $t\geq\delta$ it holds that
\begin{alignat}{2}
\label{ineq:lower_bound_p<2}
\norm{u(t)-u_\infty}^{2-p}&\geq\norm{u(\delta)-u_\infty}^2-(2-p)\Lambda(\delta)(t-\delta),\qquad&&p<2,\\
\label{ineq:lower_bound_p=2}
\norm{u(t)-u_\infty}^{2\phantom{p-}}&\geq\norm{u(\delta)-u_\infty}^2\exp\left(-2\Lambda(\delta)(t-\delta)\right),\qquad&&p=2,\\
\label{ineq:lower_bound_p>2}
\norm{u(t)-u_\infty}^{p-2}&\geq\frac{1}{\norm{u(\delta)-u_\infty}^{2-p}+(p-2)\Lambda(\delta)(t-\delta)},\qquad&&p>2,
\end{alignat}
where $\delta=0$ is admissible if $\calJ(f)<\infty$.
\end{proposition}
\begin{proof}
The proof works just as the proof of \cref{prop:upper_bounds}, utilizing that $\Lambda(t) \leq \Lambda(\delta)$ for all $t\geq \delta$ according to \cref{prop:Lambda_decreasing}.
\end{proof}

In the regime of finite extinction time we can, \emph{a posteriori}, improve the upper and lower bounds from \cref{prop:upper_bounds,prop:lower_bounds}.
\begin{proposition}[Improved bounds]\label{prop:improved_bounds}
Assume that $p<2$ and $\tex<\infty$.
Then it holds that
\begin{align}\label{eq:improved_bounds}
    (2-p)\lambda_1(\tex-t)
    \leq 
    \norm{u(t)-u_\infty}^{2-p}
    \leq 
    (2-p)\Lambda(t)(\tex-t).
\end{align}
\end{proposition}
\begin{proof}
Integrating \labelcref{eq:diss_u_Lambda} from $t$ to $\tex$ yields
\begin{align*}
    \frac{1}{2-p}\norm{u(t)-u_\infty}^{2-p}=\int_t^{\tex} \Lambda(s)\d s. 
\end{align*}
The two estimates $\lambda_1\leq\Lambda(s)\leq\Lambda(t)$ for all $t\leq s\leq\tex$ conclude the proof.
\end{proof}

\subsection{Second-order: asymptotic profiles}

Having understood the ``first-order'' asymptotic behavior of gradient flows, i.e., the convergence to $u_\infty$ in finite or infinite time, we are now interested in the much more interesting ``second-order'' asymptotic profiles (cf.~\cref{sec:gradflow_profiles_linear}).
Our goal will be to prove that under certain conditions the gradient flow of $\calJ$ satisfies
\begin{align*}
    \frac{u(t) - u_\infty}{\norm{u(t) - u_\infty}} \to w\quad\text{as }t\to\tex\qquad\text{where }\lambda w\in\dJ(w).
\end{align*}
Our main reason to think that this is true is---apart from the linear case studied in \cref{sec:gradflow_profiles_linear}---is that eigenfunctions are invariant under the gradient flow.
This is stated in the following theorem.
\begin{theorem}[Invariance of eigenfunctions]\label{thm:invariance}
Assume that $\calJ:\Hspace\to\RI$ is absolutely $p$-homogeneous and let $f\in\Hspace$ satisfy $\lambda\norm{f}^{p-2}f\in\dJ(f)$ for some $\lambda\geq 0$.
Let $a:[0,\infty)\to\R$ be the unique solution to the initial value problem
\begin{align}\label{eq:ODE_a}
    \begin{cases}
    a(0) = 1,\\
    a'(t) = -\lambda\norm{f}^{p-2} a(t)^{p-1}.
    \end{cases}
\end{align}
Then $u(t) := \max(a(t),0) f$ solves the gradient flow of $\calJ$.
\end{theorem}
\begin{proof}
\sloppy Trivially it holds $u(0)=\max(a(0),0)f=f$.
Furthermore, the function $t\mapsto \max(a(t),0)$ is non-increasing, continuous, and continuously differentiable for all $t>0$ such that $a(t)>0$.
Hence, it holds
\begin{align*}
    u'(t) = a'(t)f = -\lambda\norm{f}^{p-2} a(t)^{p-1}f \in -a(t)^{p-1}\dJ(f) = -\dJ(a(t)f) = -\dJ(u(t)).
\end{align*}
For all $t>0$ such that $a(t)>0$ we have
\begin{align*}
    u'(t) = 0 \in -\dJ(0).
\end{align*}
Since eigenfunctions are subgradients of minimal norm (see \cref{prop:ef_minimal_norm}) the proof is complete.
\end{proof}
The temporal profiles $a(t)$ which solve \labelcref{eq:ODE_a} are given by
\begin{align}
    a(t) = 
    \begin{cases}
    (1-(2-p)\lambda t)^\frac{1}{2-p},\quad &p\neq 2,\\
    \exp(-\lambda t),\quad &p=2,
    \end{cases}
\end{align}
and nicely illustrate that all previous results on extinction times and convergence rates are sharp for eigenfunctions.
Obviously, the explicit solution provided in in \cref{thm:invariance} has the property that $u_\infty=0$ and hence
\begin{align*}
    \frac{u(t)-u_\infty}{\norm{u(t)-u_\infty}} = \frac{\max(a(t),0)f}{\norm{\max(a(t),0)f}} = \frac{f}{\norm{f}}\quad\forall 0 \leq t < \tex.
\end{align*}
The function $w:=\frac{f}{\norm{f}}$ is an asymptotic profile and satisfies $\lambda w \in \dJ(w)$.
For proving existence of such asymptotic profiles in the case of general initial data, we cannot rely on linear spectral techniques as used in \cref{sec:gradflow_profiles_linear}.
Instead we have to use a variational approach to tackle the problem. 

The core ingredient will be to understand and quantify the evolution of the Rayleigh quotient $\Lambda(t)$, defined in \labelcref{eq:Lambda}.
As shown in \cref{prop:Lambda_decreasing}, it decreases over time. 
However, this is not enough and we need a sharper estimate of its decrease in order to prove existence of asymptotic profiles.

The following is the key result of this section:
\begin{proposition}\label{prop:decrease_Lambda_sharp}
Defining, $w(t) := \frac{u(t)-u_\infty}{\norm{u(t)-u_\infty}}$ it holds that
\begin{align}\label{ineq:diff_ineq_Lambda}
\frac{1}{p}\Lambda'(t)+\frac{1}{\norm{u(t)-u_\infty}^{2-p}}\norm{\frac{\subgradA(t)}{\norm{u(t)-u_\infty}^{p-1}}-\Lambda(t)w(t)}^2\leq 0    
\end{align}
\end{proposition}
\begin{proof}
Using the quotient rule, we compute
\begin{align*}
    &\phantom{=}\frac{1}{p}\Lambda'(t) 
    \\
    &= 
    \frac{\d}{\d t}\frac{\calJ(u(t))}{\norm{u(t)-u_\infty}^p}
    \\
    &=
    \frac{-\norm{u(t)-u_\infty}^p\norm{\subgradA(t)}^2 - p\calJ(u(t))\norm{u(t)-u_\infty}^{p-2}\sp{u'(t),u(t)-u_\infty}}{\norm{u(t)-u_\infty}^{2p}}
    \\
    &=
    -\frac{\norm{\subgradA(t)}^2}{\norm{u(t)-u_\infty}^{p}} + \frac{p\calJ(u(t))\sp{\subgradA(t),u(t)-u_\infty}}{\norm{u(t)-u_\infty}^{p+2}}
    \\
    &=
    -\frac{\norm{\subgradA(t)}^2}{\norm{u(t)-u_\infty}^{p}} + \Lambda(t)\frac{\sp{\subgradA(t),w(t)}}{\norm{u(t)-u_\infty}}  
    \\
    &=
    -
    \frac{1}{\norm{u(t)-u_\infty}^{2-p}}
    \left(
    \norm{\frac{\subgradA(t)}{\norm{u(t)-u_\infty}^{p-1}}}^2
    -
    2\left\langle\frac{\subgradA(t)}{\norm{u(t)-u_\infty}^{p-1}},\Lambda(t)w(t)\right\rangle
    +
    \norm{\Lambda(t)w(t)}^2
    \right)
    \\
    &=
    -\frac{1}{\norm{u(t)-u_\infty}^{2-p}}
    \norm{
    \frac{\subgradA(t)}{\norm{u(t)-u_\infty}^{p-1}}
    - \Lambda(t) w(t)
    }^2,
\end{align*}
where we utilized
\begin{align*}
    \norm{\Lambda(t)w(t)}^2
    &=
    \Lambda(t)^2
    =
    \Lambda(t) \frac{p\calJ(u(t)}{\norm{u(t)-u_\infty}^{p}}
    =
    \Lambda(t)
    \frac{\sp{\subgradA(t),u(t)-u_\infty}}{\norm{u(t)-u_\infty}^{p}}
    \\
    &=
    \Lambda(t)
    \frac{\sp{\subgradA(t),w(t)}}{\norm{u(t)-u_\infty}^{p-1}}
    =
    \left\langle
    \frac{\subgradA(t)}{\norm{u(t)-u_\infty}^{p-1}},\Lambda(t)w(t)
    \right\rangle
\end{align*}
and hence
\begin{align*}
    \frac{1}{\norm{u(t)-u_\infty}^{2-p}}\left\langle
    \frac{\subgradA(t)}{\norm{u(t)-u_\infty}^{p-1}},\Lambda(t)w(t)
    \right\rangle
    =
    \Lambda(t) \frac{\sp{\subgradA(t),w(t)}}{\norm{u(t)-u_\infty}}.
\end{align*}
\end{proof}

From this proposition we can derive two important insights: firstly, $t\mapsto\Lambda(t)$ is non-increasing (which we already know from \cref{prop:Lambda_decreasing}), and secondly that the term in the squared norm goes to zero along a subsequence as $t\to\tex$. 
This is last bit summarized in the next proposition, for which we need the following lemma.
\begin{lemma}\label{lem:liminf}
Let $f,g:(a,b)\to[0,\infty)$ be measurable functions with $\int_s^b f(t)\d t=\infty$ for all $a\leq s<b$ and $\int_a^b f(t)g(t)\d t<\infty$.
Then it holds $\liminf_{t\to b}g(t)=0$.
\end{lemma}
\begin{proof}
Assume the contrary, i.e., $\liminf_{t\to b}g(t)=\alpha>0$.
Then, for all $\eps>0$ there exists $s\in(a,b)$ such that for all $t\in(s,b)$ it holds $g(t)\geq \alpha - \eps$.
For $0<\eps<\alpha$ it hence holds
\begin{align*}
    \infty > \int_a^b f(t)g(t)\d t \geq \int_s^b f(t)g(t)\d t \geq (\alpha-\eps)\int_s^b f(t)\d t = \infty
\end{align*}
which is a contradiction.
\end{proof}

Now we are ready to prove the proposition.
\begin{proposition}\label{prop:cvgc_subgrad}
For all $0<s<\tex$ it holds that 
\begin{align}\label{eq:integral_subgrads}
    \int_s^{\tex} \frac{1}{\norm{u(t)-u_\infty}^{2-p}}\norm{\frac{\subgradA(t)}{\norm{u(t)-u_\infty}^{p-1}}-\Lambda(t)w(t)}^2 \d t < \infty.
\end{align}
Furthermore it holds that
\begin{align}\label{eq:cvgc_subgrad}
    \liminf_{t\to\tex}\norm{\frac{\subgradA(t)}{\norm{u(t)-u_\infty}^{p-1}}-\Lambda(t)w(t)} = 0.
\end{align}

\end{proposition}
\begin{proof}
The first statement simply follows from integrating \labelcref{ineq:diff_ineq_Lambda} from $s$ to $t$, using that $\calJ\geq 0$, and letting $t\to\tex$.

To show \labelcref{eq:cvgc_subgrad} we want to use \cref{lem:liminf} and shall make a case distinction based on $p$.
For $p<2$ we get from \cref{prop:improved_bounds} that
\begin{align*}
    \frac{1}{\norm{u(t)-u_\infty}^{2-p}} \geq \frac{1}{(2-p)\Lambda(\delta)(\tex-t)}
\end{align*}
for any $0<\delta\leq t$ and $\int_s^{\tex} \frac{1}{\tex-t}\d t=\infty$.
For $p=2$ we have
\begin{align*}
    \liminf_{t\to\tex}\frac{1}{\norm{u(t)-u_\infty}^{2-p}} = 1
\end{align*}
and since $\tex=\infty$ it holds $\int_s^\infty 1\d t=\infty$.
For $p>2$ we use \labelcref{ineq:lower_bound_p>2} in \cref{prop:lower_bounds} to obtain
\begin{align*}
    \frac{1}{\norm{u(t)-u_\infty}^{2-p}} \geq \frac{1}{c+ Ct}
\end{align*}
for suitable constants $c,C>0$ which satisfies $\int_s^\infty \frac{1}{c+Ct}\d t=\infty$. 

Hence, in all cases we can use \cref{lem:liminf} to conclude that \labelcref{eq:cvgc_subgrad} holds true.
\end{proof}

Now we are in position to prove that asymptotic profiles, i.e., limits of $w(t):=\frac{u(t)-u_\infty}{\norm{u(t)-u_\infty}}$, are eigenfunctions.
As for the proof of existence of ground states in \cref{thm:existence_GS}, we will have to pose a compactness condition on $\calJ$.

\begin{theorem}\label{thm:asymptotic_profiles}
Assume that $\calJ:\Hspace\to\RI$ is absolutely $p$-homogeneous, proper, convex, \lsc{,} and the sublevel sets of $u\mapsto\norm{u}+\calJ(u)$ are precompact.
Let $u:[0,\infty)\to\Hspace$ be the gradient flow of $\calJ$ with initial datum $f\in\closure{\dom(\calJ)}$, limit $u_\infty:=\lim_{t\to\infty}u(t)$, and extinction time $\tex\in[0,\infty]$.
Furthermore, denote $w(t):=\frac{u(t)-u_\infty}{\norm{u(t)-u_\infty}}$.

Then there exists a non-decreasing sequence $(t_n)_{n\in\N}\subset[0,\infty)$ such that $\lim_{n\to\infty} t_n = \tex$ and $\lim_{n\to\infty} w(t_n) = w_\infty$.
Furthermore, it holds $\Lambda_\infty w_\infty \in \dJ(w_\infty)$ where $\Lambda_\infty := \lim_{t\to\tex}\Lambda(t) =  \inf_{0<t<\tex}\Lambda(t)$.
\end{theorem}
\begin{proof}
By means of \cref{prop:cvgc_subgrad}, there exists a non-decreasing sequence $(t_n)_{n\in\N}\subset[0,\infty)$ such that $\lim_{n\to\infty} t_n = \tex$ and
\begin{align}\label{eq:cvgc_subgrad_sequence}
    \lim_{n\to\infty} \frac{\subgradA(t_n)}{\norm{u(t_n)-u_\infty}^{p-1}} - \Lambda(t_n) w(t_n) = 0.
\end{align}
By definition of the gradient flow it holds that $\subgradA(t) \in \dJ(u(t))$ for all $t>0$ which is equivalent to
\begin{align*}
    \calJ(u(t)) + \sp{\subgradA(t),v-u(t)} \leq \calJ(v)\quad\forall v\in\Hspace,\;t>0.
\end{align*}
Using this for $v := \norm{u(t_n)-u_\infty}u$, where $u\in\Hspace$ is arbitrary, we obtain
\begin{align*}
    &\phantom{=}
    \calJ(w(t_N)) + 
    \left\langle
    \frac{\subgradA(t)}{\norm{u(t_n)-u_\infty}^{p-1}},u - w(t_n)
    \right\rangle
    \\
    &=
    \frac{1}{\norm{u(t_n)-u_\infty}^p}
    \big(
    \calJ(u(t_n))
    +
    \left\langle
    \subgradA(t_n),\norm{u(t_n)-u_\infty}u - u(t_n)
    \right\rangle
    \big)
    \\
    &\leq 
    \frac{1}{\norm{u(t_n)-u_\infty}^p}
    \calJ(\norm{u(t_n)-u_\infty}u)
    =
    \calJ(u),
\end{align*}
where we exploited the homogeneity of $\calJ$.
Now we note that
\begin{align*}
    \norm{w(t_n)} + \calJ(w(t_n)) = 1 + \frac{\calJ(u(t_n))}{\norm{u(t_n)-u_\infty}^p} = 1+ \frac{1}{p}\Lambda(t_n) \leq 1+ \frac{1}{p} \Lambda(\delta)
\end{align*}
for all $\delta>0$ and all $n$ sufficiently large such that $t_n\geq \delta$.
Hence, the precompactness of the sublevel sets of $u\mapsto\norm{u}+\calJ(u)$ implies that $(w(t_n))_{n\in\N}$ possesses a subsequence (which we do not relabel) which converges to some $w_\infty\in\Hspace$ and additionally $\norm{w_\infty}=1$. 
Furthermore, the fact that $t\mapsto\Lambda(t)$ is non-increasing and bounded from below by zero implies that it has a limit and it holds
\begin{align*}
    \Lambda_\infty := \lim_{t\to\tex}\Lambda(t) = \inf_{0<t<\tex}\Lambda(t).
\end{align*}
The lower semicontinuity of $\calJ$ and the discussed convergences now imply that
\begin{align*}
    \calJ(w_\infty) + \sp{\Lambda_\infty w_\infty, u-w_\infty}
    &\leq 
    \liminf_{n\to\infty} 
    \calJ(w(t_n)) + \sp{\Lambda(t_n) w(t_n), u-w(t_n)}
    \\
    &=
    \liminf_{n\to\infty} 
    \calJ(w(t_n)) + \left\langle{\frac{\subgradA(t_n)}{\norm{u(t_n)-u_\infty}^{p-1}}, u-w(t_n)}\right\rangle
    \\
    &\leq
    \calJ(u).
\end{align*}
Since $u\in\Hspace$ was arbitrary, this means that $\Lambda_\infty w_\infty\in\dJ(w_\infty)$, as claimed.
\end{proof}

Using the following classical result from the theory of abstract evolution equations, we can prove uniqueness of the asymptotic profile for non-negative gradient flows.
This is particularly relevant in situations where the ground state is known to be the unique non-negative eigenfunction (cf.~\cref{prop:positive_gs,thm:uniqueness_gs}).
In such cases our next results show that the asymptotic profile is a ground state.

\begin{lemma}[Crandall--Benilan \cite{crandall1980regularizing}]\label{lem:crandall}
Assume that $\calJ:\Hspace\to\RI$ is absolutely $p$-homogeneous, proper, convex, and \lsc{,} and assume that $(\Hspace,\geq)$ is a partially ordered space.
Let $S(t)[f]$ denote the the gradient flow of $\calJ$ with initial datum $f$ at time $t\geq 0$ and assume that 
\begin{align}\label{eq:monotonicity}
    S(t)[f] \geq S(t)[g]\quad\text{if }f\geq g.
\end{align}
If $p>2$ and $f\geq 0$, then it holds that
\begin{align}
    \partial_t S(t)[f] \geq -\frac{1}{p-2} \frac{S(t)[f]}{t},\quad\forall t>0.
\end{align}
\end{lemma}
\begin{proof}
It is an easy exercise to check that the solution map satisfies the homogeneity property $S(t)[cf] = c S(c^{p-2}t)[f]$ for all $c>0$.
For $h>0$ and $\lambda := 1+h/t$ we hence get
\begin{align*}
    S(t+h)[f] - S(t)[f]
    &=
    S(\lambda t)[f] - S(t)[f]
    \\
    &=
    \lambda^\frac{1}{2-p}S(t)[\lambda^\frac{1}{p-2}f] - S(t)[f]
    \\
    &=
    \lambda^\frac{1}{2-p}\left(
    S(t)[\lambda^\frac{1}{p-2}f] - S(t)[f]
    \right)
    +
    \left(
    \lambda^\frac{1}{2-p}-1
    \right)
    S(t)[f].
\end{align*}
Since $p>2$ we get that $\lambda^\frac{1}{p-2}\geq 1$ and hence $S(t)[\lambda^\frac{1}{p-2}f] - S(t)[f]\geq 0$ by the assumption on $S(t)$ and $f\geq 0$.
Dividing by $h>0$ and plugging in the definitions of $\lambda$ and $S(t)[f]$ yields
\begin{align*}
    \frac{S(t+h)[f] - S(t)[f]}{h} \geq \frac{1}{h}\left(
    \left(1+\frac{h}{t}\right)^\frac{1}{2-p}-1
    \right)
    S(t)[f].
\end{align*}
Letting $h\downarrow0$ and using L'Hôpital's rule yields the assertion.
\end{proof}

\begin{theorem}[Uniqueness of the profile]\label{thm:uniqueness}
Under the assumptions of \cref{thm:asymptotic_profiles}, assume that $(\Hspace,\geq)$ is a partially ordered space, $f\geq 0$, and the gradient flow satisfies \labelcref{eq:monotonicity}.
Furthermore, assume that $\lambda_1>0$ and that $u_\infty = 0$.
Then there exists at most one asymptotic profile.
\end{theorem}
\begin{proof}
The proof techniques for $p\leq 2$ and $p>2$ are somewhat different. We begin with the latter case. 
Using \cref{lem:crandall}, for positive solutions one has
$$\partial_tu(t)\geq-\frac{1}{p-2}\frac{u(t)}{t}.$$
This estimate immediately implies that for $\tilde{w}(t):=t^\frac{1}{p-2}u(t)$ it holds 
$$\partial_t\tilde{w}(t)=\frac{1}{p-2}t^{\frac{1}{p-2}-1}u(t)+t^\frac{1}{p-2}\partial_t u(t)\geq 0.$$
Hence, $\tilde{w}(t)$ increases monotonously and, thus, has at most one accumulation point $\tilde{w}_*$. Thanks to \cref{prop:upper_bounds,prop:lower_bounds} and $\lambda_1>0$ we get that
$$\norm{w(t)-c\tilde{w}(t)}\to 0,\quad t\to\infty$$ 
for a suitable $c>0$ and we obtain that $w(t)$ has at most one asymptotic profile $w_*=c\tilde{w}_*$.

In the case $p\leq 2$ we utilize a simple transformation in order to exploit the uniqueness of the asymptotic profile in the case $p>2$. Without loss of generality, we assume that $f\in\dom(\partial J)$ which implies that $t\mapsto J(u(t))$ is globally Lipschitz continuous on the interval $(0,\infty)$. We let $q>2p$ and define $r:=q/p>2$. Then the functional $\calJ(u)^r/r$ is absolutely $q$-homogeneous, convex, and lower semi-continuous. Furthermore, by application of the chain rule, one can easily see that its gradient flow is given by $\tilde{u}(\tau)=u(\varphi(\tau))$, where $u$ solves the gradient flow with respect to $\calJ$ and the time reparametrization $\varphi$ is defined by the initial value problem
$$
\begin{cases}
\varphi'(\tau)&=\calJ(u(\varphi(\tau)))^{r-1},\quad\tau>0,\\
\varphi(0)&=0.
\end{cases}
$$
Existence and uniqueness of $\varphi$ on the positive reals follow from standard arguments utilizing that $t\mapsto \calJ(u(t))^{r-1}$ is globally Lipschitz continuous. To see this it is sufficient to note that this map is Lipschitz continuous on $[0,
\tex]$ and zero on $(\tex,\infty)$. Furthermore, $\tau\mapsto\varphi(\tau)$ is invertible on the interval $(0,\tex)$ since it is strictly increasing, and it holds $\varphi^{-1}(t)\to\infty$ as $t\to\tex$. The latter statement follows from the fact that there is no $\tau>0$ such that $\varphi(\tau)=\tex$ since otherwise $\tilde{u}$ would have finite extinction time which is impossible according to \cref{prop:lower_bounds}. 

From the first part of the theorem we know that the rescalings of $\tilde{u}$, given by 
$$\tilde{w}(\tau)=\frac{\tilde{u}(\tau)}{\norm{\tilde u(t)}}$$
have at most one limit and without loss of generality we assume that $\tilde{w}(\tau)$ converges to some $\tilde{w}_*$ as $\tau\to\infty$. 
Consequently, we have for $t\in(0,\tex)$ that
\begin{align}\label{eq:transformation_profiles}
w(t)=\frac{u(t)}{\norm{u(t)}}=\frac{\tilde{u}(\varphi^{-1}(t))}{\norm{\tilde u(\varphi^{-1}(t))}} = \tilde w(\varphi^{-1}(t)),
\end{align}
which converges to $\tilde{w}_*$ since $\varphi^{-1}(t)\to\infty$ as $t\to\tex$. 
\end{proof}

\begin{remark}\label{rem:KB-spaces}
Partially ordered spaces (more precisely, Banach lattices~\cite{Meyer-Nieberg}) with the property that every norm bounded monotone sequence (in the sense of the partial order) is called a Kantorovich--Banach space~\cite[Def. 2.4.1]{Meyer-Nieberg}. 
Thus, since $\tilde w(t)$ is non-decreasing (in the sense of the partial order on $\Hspace$) and bounded in norm, it converges if $\Hspace$ is a Kantorovich--Banach space. \cite[Thm. 2.5.6]{Meyer-Nieberg} gives a characterisation of Kantorovich--Banach spaces. In particular, Banach lattices that are weakly sequentially complete are Kantorovich--Banach spaces. Hilbert spaces are weakly sequentially complete.
\end{remark}

\section{Examples}

\subsection{Distance function}

Let $\Omega\subset\R^n$ be an open and bounded domain with a piecewise smooth boundary $\partial\Omega$. 
We define the function space
\begin{align}
    W^{1,\infty}_0(\Omega):=\{u\in W^{1,\infty}(\Omega)\st u=0\text{ on }\partial\Omega\}
\end{align}
which consists of all Lipschitz continuous functions, vanishing on $\partial\Omega$ (and, at the same time, is a Sobolev space). 
In this section we consider the following functional 
\begin{align}\label{eq:fctl_v2}
\calJ(u)=
\begin{cases}
\norm{\nabla u}_\infty,\quad &u\in W^{1,\infty}_0(\Omega),\\
+\infty,\quad &u\in L^2(\Omega)\setminus W^{1,\infty}_0(\Omega),
\end{cases}
\end{align}
which coincides with the Lipschitz constant if $u\in W^{1,\infty}_0(\Omega)$. We are interested in its   ground states, 
i.e., minimizers of the following  Rayleigh quotient
\begin{align}\label{eq:ground_state}
u^\ast\in\argmin_{u\in W^{1,\infty}_0(\Omega)}\frac{\calJ(u)}{\norm{u}_2}.
\end{align}
We will prove that---up to multiplicative constants---they coincide with the distance function of the boundary $\partial\Omega$ of the domain which is defined as
\begin{align}\label{eq:distance_fct}
d_\Omega (x):=\dist(x,\partial\Omega):=\inf_{y\in\partial\Omega}|x-y|.
\end{align}
The following lemma is essential for our arguments.
\begin{lemma}\label{lem:lipschitz}
For every function $u\in W^{1,\infty}_0(\Omega)$ and every $x,y\in\Omega$ it holds
\begin{align*}
    \abs{u(x)-u(y)} \leq \norm{\grad u}_\infty \abs{x-y}.
\end{align*}
\end{lemma}
\begin{proof}
Since $u=0$ on $\partial\Omega$, the function $\tilde u:\R^n\to\R$, defined through $\tilde u(x)=u(x)$ for $x\in\Omega$ and $\tilde u(x)=0$ for $x\in\R^n\setminus\Omega$, satisfies $\norm{\grad \tilde u}_\infty = \norm{\grad u}_\infty$. 
Furthermore, it holds for all $x,y\in\Omega$
\begin{align*}
    u(x)-u(y)
    &=
    \tilde u(x)-\tilde u(y)
    =
    \int_0^1 \frac{\d}{\d t}\tilde u(tx+(1-t)y)\d t
    \\
    &=
    \int_0^1 \nabla \tilde u(tx+(1-t)y)\cdot(x-y) \d t
    \leq 
    \norm{\grad\tilde u}_\infty \abs{x-y}
    =
    \norm{\grad u}_\infty \abs{x-y}.
\end{align*}
\end{proof}
Using this lemma, we prove the following properties of the distance function.
\begin{proposition}\label{prop:dist-fun-maximal}
For every function $u\in W^{1,\infty}_0(\Omega)$ and every $x\in\Omega$ it holds
\begin{align}\label{ineq:distance_estimate}
    \abs{u(x)} \leq \norm{\grad u}_\infty d_\Omega(x).
\end{align}
Consequently, the distance function satisfies
\begin{align}\label{eq:grad_distance}
    \norm{\grad d_\Omega}_\infty = 1.
\end{align}
\end{proposition}
\begin{proof}
Letting $y\in\partial\Omega$ denote a point such that $\abs{x-y}=d_\Omega(x)$ \cref{lem:lipschitz} implies
\begin{align*}
    \abs{u(x)} = \abs{u(x) - 0} = \abs{u(x) - u(y)} \leq \norm{\grad u}_\infty \abs{x-y} = 
    \norm{\grad u}_\infty d_\Omega(x). 
\end{align*}
Applying this to $u=d_\Omega \geq 0$, we get that $\norm{\grad d_\Omega}_\infty \geq 1$. The opposite inequality is a trivial computation of the Lipschitz constant.
\end{proof}

We actually have a local result, too.

\begin{proposition}\label{prop:dist-fun-grad}
The distance function~\eqref{eq:distance_fct} satisfies
\begin{equation*}
    \abs{\grad d_\Omega(x)} = 1 \quad \text{a.e. in $\Omega$}.
\end{equation*}
\end{proposition}
\begin{proof}
By \cref{prop:dist-fun-maximal}, we have that $\abs{\grad d_\Omega(x)} \leq 1$ a.e.
Let us show that the opposite inequality holds, too. Let $x \in \Omega$ be arbitrary and 
\begin{equation*}
    y \in \argmin_{y' \in \partial\Omega} \abs{x-y'}
\end{equation*}
be the closest point on the boundary to $x$ (it does not have to unique, but it exists). Let
\begin{equation*}
    y(t) \defeq x + t(y-x), \quad t \in [0,1],
\end{equation*}
be a parametrisation of the shortest path between $x$ and $y$.  Then we have for any $t \in [0,1]$
\begin{equation*}
    d_{\Omega}(y(t)) = \int_t^1 \abs{\frac{d y}{\d t}(t')} \d t'.
\end{equation*}
Differentiating this with respect to $t$, we get
\begin{equation*}
    \sp{\grad d_{\Omega}(y(t)),\frac{d y}{\d t}(t)} = - \abs{\frac{d y}{\d t}(t)}.
\end{equation*}
With Cauchy-Schwarz we finally get
\begin{equation*}
    \abs{\grad d_{\Omega}(y(t))}\abs{\frac{d y}{\d t}(t)} \geq \abs{\sp{\grad d_{\Omega}(y(t)),\frac{d y}{\d t}(t)}} = \abs{\frac{d y}{\d t}(t)}
\end{equation*}
and $\abs{\grad d_{\Omega}y(t))} \geq 1$ for all $t \in [0,1]$. Taking $t=0$, we get $\abs{\grad d_{\Omega}(x)} \geq 1$. Since $x$ was arbitrary, this implies the assertion.
\end{proof}

\cref{prop:dist-fun-maximal} immediately implies that the distance function minimizes the Rayleigh quotient.
\begin{theorem}[Ground states are distance functions]\label{thm:ground_states_distance}
All solutions $u^\ast$ to  \eqref{eq:ground_state} are multiples of the distance function to $\partial\Omega$, given by \eqref{eq:distance_fct}.
\end{theorem}
\begin{proof}
Let $u\in W^{1,\infty}_0(\Omega)$ be arbitrary.
Squaring and integrating \labelcref{ineq:distance_estimate} yields
\begin{align*}
    \int_\Omega\abs{u(x)}^2\d x \leq \norm{\grad u}_\infty^2 \int_\Omega d_\Omega(x)^2 \d x.
\end{align*}
Taking square roots and using \labelcref{eq:grad_distance} implies
\begin{align*}
    \frac{\norm{\grad d_\Omega}_\infty}{\norm{d_\Omega}_2}=
    \frac{1}{\norm{d_\Omega}_2} \leq \frac{\norm{\grad u}_\infty}{\norm{u}_2}
\end{align*}
which shows that $d_\Omega$ (and all its non-zero scalar multiples) minimizes \labelcref{eq:ground_state}.

Regarding uniqueness, if $u$ is no multiple of the distance function, then \labelcref{ineq:distance_estimate} is a strict inequality on a set of non-zero Lebesgue measure.
As before, by integrating we would get
\begin{align*}
    \frac{\norm{\grad d_\Omega}_\infty}{\norm{d_\Omega}_2}< \frac{\norm{\grad u}_\infty}{\norm{u}_2}
\end{align*}
which contradicts the minimality of $d_\Omega$.
\end{proof}

\begin{remark}
Clearly, this result applies to Rayleigh quotients of the form $\frac{\norm{\grad u}_\infty}{\norm{u}_p}$ for any $L^p$-norm with $1 \leq p < \infty$.
For the $L^\infty$-norm multiples of the distance function are still minimisers, however, not the only ones.
\end{remark}

Next we state that every non-negative eigenfunction coincides with a ground state, i.e., is a multiple of the distance function to $\partial\Omega$.
\begin{proposition}[Uniqueness of non-negative eigenfunction]\label{prop:positivity_efs}
Any non-negative function $u\neq 0$ which meets $\lambda u\in\partial\calJ(u)$ is a ground state. 
\end{proposition}
\begin{proof}
Let us assume that we have a non-negative eigenfunction $u\neq 0$ on $\Omega$ which is no ground state. We can normalize in such a way that $\calJ(u)=1$.  We recall that $u\leq d_{\Omega}$ holds pointwise almost everywhere in $\Omega$. 
Define the set 
$$\Omega_\varepsilon:=\left\lbrace x\in\Omega\st d_{\Omega}(x)>u(x)+\varepsilon,\; u(x)>\varepsilon\right\rbrace.$$
Since $u$ is an eigenfunction, it holds $\langle\lambda u,v\rangle \leq\calJ(v)$ for all $v\in L^2(\Omega)$, where $\lambda=1/\norm{u}_2^2$. Testing this with $v=d_{\Omega}$, using the definition of $\Omega_\varepsilon$ and the fact that $d_{\Omega}\geq u$, we obtain
\begin{align*}
    \norm{u}_2^2
    &\geq\int_\Omega u(x)d_{\Omega}(x)\dx
    \geq\int_{\Omega_\varepsilon}u(x)(u(x)+\varepsilon)\dx + \int_{\Omega\setminus \Omega_{\varepsilon}}u(x)d_{\Omega}(x)\dx \\
    &\geq\int_\Omega u(x)^2\dx +\varepsilon\int_{\Omega_\varepsilon}u(x)\dx\\
    &\geq\norm{u}_2^2+\varepsilon^2|\Omega_\varepsilon|,
\end{align*}
which tells us that $|\Omega_\varepsilon|=0$. Letting $\varepsilon$ tend to zero we infer as before that almost everywhere in $\Omega$ it holds $u=d_{\Omega}$ or $u=0$. Since however both $u$ and $d_{\Omega}$ are continuous functions and by assumption $u\neq 0$, we find that $u=d_{\Omega}$ holds almost everywhere in $\Omega$. 
\end{proof}

The following result is a consequence of \cref{thm:asymptotic_profiles,thm:uniqueness,thm:ground_states_distance,prop:positivity_efs}.
\begin{theorem}[Asymptotic profiles]\label{thm:profiles_are_distance_fcts}
Let $u(t)$ be the solution of the gradient flow \eqref{eq:grad-flow} with respect to $\calJ$ and datum $f \in \closure{\dom(\calJ)}$ such that $f\geq 0$ a.e. Denote the finite extinction time of the flow by $T_{ex}$ and assume that $u_\infty=0$. Then $u(t)/\norm{u(t)}_2$ converges strongly in $L^2(\Omega)$ to a multiple of the distance function as $t\nearrow T_{ex}$.
\end{theorem}
\begin{proof}
Since $\dom(\calJ)=W^{1,\infty}_0(\Omega)$ is compactly embedded in $L^2(\Omega)$, we can apply \cref{thm:asymptotic_profiles} to conclude that a subsequence of $u(t)/\norm{u(t)}_2$ converges to an eigenfunction. One can show that the flow~\eqref{eq:grad-flow} preserves non-negativity and satisfies \labelcref{eq:monotonicity}, hence this \ef{} is non-negative. Furthermore, by \cref{thm:uniqueness}  this eigenfunction is unique and the whole sequence converges to it if conditions described in \cref{rem:KB-spaces} are met.  From Proposition~\ref{prop:positivity_efs} and Theorem~\ref{thm:ground_states_distance} we conclude that this eigenfunction has to be a multiple of the distance function.
\end{proof}

\clearpage

\subsection{Asymptotic behaviour of nonlinear PDEs}

In this section we shall apply our findings to some nonlinear parabolic partial differential equations which admit a gradient flows structure.
In particular, we will study equations of fast diffusion / porous medium type as well as the total variation flow.

\paragraph{Porous medium and fast diffusion}

Let us first study porous medium and fast diffusion equations which take the form
\begin{align}\label{eq:FDE_PME}
    \partial_t u = \Delta u^{[m]},
\end{align}
where $m>0$ is a constant and $u^{[m]} := \sign(u)\abs{u}^{m}$. 
Obviously, for non-negative solutions it holds $u^{[m]}=u^m$.
For $0<m<1$ one refers to \labelcref{eq:FDE_PME} as fast diffusion equation, for $m=1$ it coincides with the linear heat equation \labelcref{eq:heat_eq}, and for $m>1$ it is called porous medium or slow diffusion equation. 
These terminologies stem from the fact that for $m\leq 1$ solutions have infinite speed of propagation whereas for $m>1$ it is finite. 
In the latter regime the PDE is used, for instance, for modelling the spread of contaminants in soil; hence the name. 

As it turns out, when posed on a bounded domain $\Omega\subset\R^n$ the PDE \labelcref{eq:FDE_PME} has a Hilbertian gradient flow structure.
For simplicity we limit the discussion to the case of homogeneous Dirichlet boundary conditions.
In this case, one introduces the Hilbert space $H^{-1}(\Omega)$ as the dual space of the Hilbert space $H^1_0(\Omega)$.
We can define the convex, and absolutely $m+1$-homogeneous functional
\begin{align}\label{eq:FDE_PME_func}
    \calJ(u) := 
    \begin{cases}
    \frac{1}{m+1}\int_\Omega\abs{u}^{m+1}\d x, \quad&\text{if }u\in L^{m+1}(\Omega)\cap H^{-1}(\Omega),\\
    \infty,\quad&\text{else}.
    \end{cases}
\end{align}
We will show that \labelcref{eq:FDE_PME} is the gradient flow of \labelcref{eq:FDE_PME_func} in the Hilbert space $H^{-1}(\Omega)$.
For this we first have to understand this space better.
\begin{lemma}\label{lem:isomorphism_H1}
Let $\Omega\subset\R^n$ such that the Poincar\'e inequality
\begin{align}\label{ineq:poincare}
    \norm{u}_{L^2(\Omega)} \leq C\norm{\grad u}_{L^2(\Omega)},\quad\forall u\in H^1_0(\Omega),
\end{align}
holds true.
Then $-\laplace : H^1_0(\Omega)\to H^{-1}(\Omega)$ is an isometric isomorphism of $H^{-1}(\Omega)$ and the set distributions $\{-\laplace u \st u\in H^1_0(\Omega)\}$.
\end{lemma}
\begin{proof}
Thanks to the Poincar\'e inequality the space $H^1_0(\Omega)$ can be equipped with the inner product $\sp{u,v}_{H^1_0(\Omega)}:=\int_\Omega\grad u\cdot\grad v\d x$.

Obviously, for $u\in H^1_0(\Omega)$ the map
\begin{align*}
    (-\laplace u)(v) := \int_\Omega\grad u\cdot\grad v\d x,\quad v\in H^1_0(\Omega)
\end{align*}
is linear and bounded which shows $-\laplace:H^1_0(\Omega)\to H^{-1}(\Omega)$.
Furthermore, this map is injective since $\nabla u=0$ for $u\in H^1_0(\Omega)$ implies that $u=0$.

By the Riesz representation theorem, for any $w\in H^{-1}(\Omega)$ there exists $v\in H^1_0(\Omega)$ such that
\begin{align*}
    w(u) = \sp{v,u}_{H^1_0} = \int_\Omega\grad v\cdot\grad u\d x,\quad\forall u\in H^1_0(\Omega).
\end{align*}
Choosing $u\in C^\infty_c(\Omega)$ shows that $w = -\laplace v$ in the sense of distributions.
Hence, $-\laplace$ is also surjective and therefore an isomorphism.
The fact that $-\laplace$ is isometric follows from the following computation for $w:=-\laplace v\in H^{-1}(\Omega)$:
\begin{align*}
    \norm{w}_{H^{-1}(\Omega)} 
    &= \sup_{u\in H^1_0(\Omega)}\frac{w(u)}{\norm{u}_{H^1_0(\Omega)}}
    =
    \sup_{u\in H^1_0(\Omega)}\frac{\int_\Omega\grad v\cdot\grad u\d x}{\norm{u}_{H^1_0(\Omega)}}
    \\
    &=
    \sup_{u\in H^1_0(\Omega)}\frac{\sp{v,u}_{H^1_0(\Omega)}}{\norm{u}_{H^1_0(\Omega)}} = \norm{v}_{H^1_0(\Omega)}.
\end{align*}
\end{proof}
Thanks to \cref{lem:isomorphism_H1} we can identify $H^{-1}(\Omega)$ with $\{-\laplace v\st v\in H^1_0(\Omega)\}$ and equip it with the inner product
\begin{align}\label{eq:inner_prod_H-1}
    \sp{u,v}_{H^{-1}(\Omega)} 
    \defeq 
    \sp{(-\laplace)^{-1}(u),(-\laplace)^{-1}(v)}_{H^1_0(\Omega)}
    =
    \sp{(-\laplace)^{-1}u,v}_{L^2(\Omega)}
\end{align}
For the rest of this section we assume that $L^{m+1}(\Omega)\subset\subset H^{-1}(\Omega)$ which is true whenever $m \geq \frac{n-2}{n+2}$, see \cite{littig2015porous}.
In this case the functional in \labelcref{eq:FDE_PME_func} is proper.
We can use the above identification to compute its subdifferential.

\begin{proposition}
The subdifferential of $\calJ$ in \labelcref{eq:FDE_PME_func} with respect to $H^{-1}(\Omega)$ in $u\in H^{-1}(\Omega)$ is given by
\begin{align}\label{eq:subdiff_FDE_PME}
    \partial\calJ(u) =
    \begin{cases}
    \left\lbrace
    -\laplace u^{[m]} 
    \right\rbrace,\quad&\text{if }u\in L^{m+1}(\Omega)\cap H^{-1}(\Omega),\;-\laplace u^{[m]}\in H^{-1}(\Omega)\\
    \emptyset,\quad&\text{else}.
    \end{cases}
\end{align}
\end{proposition}
\begin{proof}
By definition of the inner product \labelcref{eq:inner_prod_H-1} and using Höder's and Young's inequalities it holds
\begin{align*}
    &\phantom{=}
    \calJ(u) + \sp{-\laplace u^{[m]},v-u}_{H^{-1}(\Omega)}
    \\
    &=
    \calJ(u) + \sp{u^{[m]},v-u}_{L^2(\Omega)}
    \\
    &=
    \frac{1}{m+1}\int_\Omega\abs{u}^{m+1}\d x 
    +
    \int_\Omega u^{[m]}v\d x - \int_\Omega\abs{u}^{m+1}\d x
    \\
    &\leq 
    -\frac{m}{m+1}\int_\Omega\abs{u}^{m+1}\d x
    +
    \left(\int_\Omega\abs{u^{[m]}}^\frac{m+1}{m}\d x\right)^\frac{m}{m+1}\left(\int_\Omega\abs{v}^{m+1}\d x\right)^\frac{1}{m+1}
    \\
    &\leq
    \frac{1}{m+1}\int_\Omega\abs{v}^{m+1}\d x = \calJ(v),\qquad\forall v\in L^{m+1}(\Omega)\cap H^{-1}(\Omega).
\end{align*}
This shows that $-\laplace u^{[m]}\in\partial\calJ(u)$.
It is an exercise to show that $\partial\calJ(u)=\{-\laplace u^{[m]}\}$.
\end{proof}
This allows us to interpret the PDE \labelcref{eq:FDE_PME} as gradient flow of the functional in \labelcref{eq:FDE_PME_func}.
\begin{remark}[Boundary conditions]
It might be not immediately clear why the interpretation of \labelcref{eq:FDE_PME} as gradient flow of $\calJ$ in $H^{-1}(\Omega)$ imposes homogeneous Dirichlet boundary conditions.
However, an inspection of the subdifferential \labelcref{eq:subdiff_FDE_PME} shows that $\partial\calJ(u)$ being non-empty implies that $-\laplace u^{[m]}\in H^{-1}(\Omega)$.
Then \cref{lem:isomorphism_H1} shows that $u^{[m]}\in H^1_0(\Omega)$ and hence $u$ has homogeneous Dirichlet boundary conditions.
\end{remark}
Now we turn to the nonlinear eigenvalue problem which is associated with the PDE \labelcref{eq:FDE_PME} and takes the form
\begin{align}\label{eq:EVP_FDE_PME}
    \lambda u \norm{u}^{m-1}_{H^{-1}(\Omega)} = -\laplace u^{[m]}.
\end{align}
Note that for $m=1$ we recover the Laplacian eigenvalue problem \labelcref{eq:helmholtz}.
The associated Rayleigh quotient is 
\begin{align}\label{eq:rayleigh_FDE_PME}
    \calR(u) = \frac{(m+1)\calJ(u)}{\norm{u}_{H^{-1}(\Omega)}^{m+1}} = \frac{\int_\Omega\abs{u}^{m+1}\d x}{\left(\int_\Omega (-\laplace)^{-1}(u) u\d x\right)^\frac{m+1}{2}}.
\end{align}
The existence of a positive minimal eigenvalue and a corresponding ground state is assured if $L^{m+1}(\Omega)$ compactly embeds in $H^{-1}(\Omega)$.
\begin{theorem}
Let $\Omega\subset\R^n$ be a domain such that the Poincar\'e inequality \labelcref{ineq:poincare} is satisfied and let $m>\frac{n-2}{n+2}$.
The Rayleigh quotient \labelcref{eq:rayleigh_FDE_PME} admits a minimiser and every minimiser solves \labelcref{eq:EVP_FDE_PME} with positive minimal eigenvalue.
\end{theorem}
\begin{proof}
The proof follows from our abstract results (cf. \cref{thm:existence_GS}) and we just sketch the missing details, following \cite{littig2015porous}.
Establishing lower semicontinuity of $\calJ$ w.r.t. the topology of $H^{-1}(\Omega)$ is left as an exercise.
The Sobolev embedding $H^1_0(\Omega)$ into $L^p(\Omega)$ for suitable values of $p$ implies an embedding of $L^q(\Omega)$ into $H^{-1}(\Omega)$ where $q = \frac{p}{p-1}$.
Choosing the appropriate values of $p$ shows that for $m\geq \frac{n-2}{n+2}$ we have
\begin{align}\label{ineq:poincare_FDE_PME}
    \norm{u}_{H^{-1}(\Omega)}\leq C\norm{u}_{L^{m+1}(\Omega)},\quad\forall u\in H^{-1}(\Omega)
\end{align}
and for $m>\frac{n-2}{n+2}$ we even have a compact embedding of $L^{m+1}(\Omega)$ into $H^{-1}(\Omega)$.
Hence, the Rayleigh quotient \labelcref{eq:rayleigh_FDE_PME} is bounded from below by a positive constant which shows positivity of the first eigenvalue.
Now the sublevel sets of $\norm{u}+\calJ(u)$ are relatively compact and  \cref{thm:existence_GS} shows the existence of ground states.
\end{proof}

An easy consequence of our results on extinction times from \cref{ss:extinction_times} is that for $m\geq 1$ the porous medium / heat equation has infinite extinction time whereas for $\frac{n-2}{n+2}\leq m < 1$ it has finite extinction time, thanks \labelcref{ineq:poincare_FDE_PME}.

Similarly, we also get existence of asymptotic profiles by applying \cref{thm:asymptotic_profiles}.
Since the PDE \labelcref{eq:FDE_PME} is also positivity preserving, see \cite{littig2015porous} for the proof, we can collect all statements that we get in the following theorem.
\begin{theorem}
Let $\Omega\subset\R^n$ be a domain such that the Poincar\'e inequality \labelcref{ineq:poincare} is satisfied and let $m>\frac{n-2}{n+2}$.
Let $u:[0,\infty)\to H^{-1}(\Omega)$ solving
\begin{align*}
    \begin{cases}
    \partial_t u = \laplace u^{[m]},\quad&\text{in }\Omega,\\
    u = 0,\quad&\text{on }\partial\Omega,
    \end{cases}
\end{align*}
denote the gradient flow of \labelcref{eq:FDE_PME_func} with non-negative initial datum $f\in H^{-1}(\Omega)$.
If $m<1$ then the extinction time $\tex$ is finite and if $m\geq 1$ it is infinite. 
Furthermore, $w(t):=\frac{u(t)}{\norm{u(t)}_{H^{-1}(\Omega)}}$ converges to $w\in H^{-1}(\Omega)$ as $t\to\tex$, where $w\geq 0$ solves $\lambda w=-\laplace w^{[m]}$.
\end{theorem}
\begin{remark}
Indeed one can even prove that the unique non-negative solutions to the eigenvalue problem \labelcref{eq:EVP_FDE_PME} are ground states, i.e., minimisers of the Rayleigh quotient \labelcref{eq:rayleigh_FDE_PME}.
Hence, the previous theorem even implies convergence of $w(t)$ to a ground state.
\end{remark}

\paragraph{The total variation flow}

In this section we study the total variation flow which formally is given by the PDE
\begin{align}\label{eq:TVF}
    \partial_t u = \div\left(\frac{\grad u}{\abs{\grad u}}\right).
\end{align}
Apart from having applications in physics (e.g., crystal growth) it has been studied extensively in the field of image processing.
There it is typically used for denoising images while preserving sharp edges.
The rationale behind this is that in regions where the gradient $\grad u$ is large the diffusivity in \labelcref{eq:TVF}, given by $\frac{1}{\abs{\grad u}}$, is very small.
Conversely, where the gradient is small or close to zero the diffusivity becomes extremely large. 
This leads to solutions which are piecewise constant which is a good model for cartoon images.

For interpreting the total variation flow as gradient flow one has to make sense of the term $\frac{\grad u}{\abs{\grad u}}$ where $\grad u=0$.
Doing a formal computation, \labelcref{eq:TVF} looks like the gradient flow of the functional
\begin{align*}
    u \mapsto \int_\Omega\abs{\grad u}\d x.
\end{align*}
It might be tempting to define this functional on $W^{1,1}(\Omega)\cap L^2(\Omega)$. 
However---due to the fact that the unit ball in $W^{1,1}(\Omega)$ does not have any compactness properties---one replaces it with the space $BV(\Omega)$, being defined as follows:
\begin{align}
    BV(\Omega) = \{u\in L^1(\Omega) \st \tv(u)<\infty\}.
\end{align}
Here, the so-called total variation is defined by
\begin{align}\label{eq:TV}
    \tv(u) = \sup\left\lbrace \int_\Omega u\div \phi \d x \st \phi\in C^\infty_c(\Omega;\R^n),\;\sup_{\Omega}\abs{\phi}\leq 1 \right\rbrace.
\end{align}
Integrating by parts it is obvious that $\tv(u)=\int_\Omega\abs{\grad u}\d x$ if $u\in W^{1,1}(\Omega)$ but the space $BV(\Omega)$ is much larger.
Furthermore, $\tv$ is convex, absolutely $1$-homogeneous, and---being a supremum over linear functionals---lower semicontinuous.
The total variation flow \labelcref{eq:TVF} is now defined as gradient flow of the total variation w.r.t. $L^2(\Omega)$.

As seen in \cref{ex:dual_unit_ball} its subdifferential is given by
\begin{align}
    \partial\tv(u) = 
    \left\lbrace
    -\div \phi \st \phi\in W^{1,2}_0(\operatorname{div};\R^2),\;\sup_\Omega\abs{\phi}\leq 1,\;
    \int_\Omega-\div\phi u\d x = \tv(u)
    \right\rbrace.
\end{align}
For smooth $u$ with $\grad u \neq 0$ on $\Omega$ it holds that $-\div\phi\in\partial\tv(u)$, where $\phi=\frac{\grad u}{\abs{\grad u}}$ and we define the $1$-Laplacian operator
\begin{align}
    -\div\left(\frac{\grad u}{\abs{\grad u}}\right) := \partial^0\tv(u)
\end{align}
as the subgradient of minimal norm of $u$.

For sufficiently nice and bounded $\Omega\subset\R^n$ the space $BV(\Omega)$ is continuously embedded in $L^{\frac{n}{n-1}}(\Omega)$ ($L^\infty(\Omega)$ if $n=1$).
Furthermore, for all $p<\frac{n}{n-1}$ one has a compact embedding into $L^p(\Omega)$.
In particular, for $n\in\{1,2\}$ we have the Poincar\'e inequality
\begin{align}\label{ineq:poincare_TV}
    \norm{u - u_\Omega}_{L^2(\Omega)} \leq C \tv(u),\quad\forall u\in L^2(\Omega),
\end{align}
where $u_\Omega:=\frac{1}{\abs{\Omega}}\int_\Omega u\d x$.
In this case the first eigenvalue
\begin{align}\label{eq:rayleigh_TV}
    \lambda_1 := \inf_{u\in L^2(\Omega)}\frac{\tv(u)}{\norm{u-u_\Omega}_{L^2(\Omega)}}
\end{align}
is positive thanks to~\eqref{ineq:poincare_TV}. 
The associated eigenvalue problem is
\begin{align}\label{eq:EVP_TV}
    \lambda u = -\div\left(\frac{\grad u}{\abs{\grad u}}\right).
\end{align}

The existence of a ground state is a delicate matter unless $n=1$. Since the embedding $\bv(\Omega) \hookrightarrow L^2(\Omega)$ is only continuous but not compact for $\Omega \subset \R^2$, the sublevel sets of $u \mapsto \norm{u}_{L^2} + \tv(u)$ are not precompact and \cref{thm:existence_GS} does not apply.

The situation is better for the $\tv$ flow. It is known that if the initial datum $f \in L^\infty(\Omega)$ then the flow is uniformly bounded in $L^\infty(\Omega)$, cf. \cite{andreu2002some}
\begin{equation*}
    \norm{u(t)}_{L^\infty} \leq C, \quad t \geq 0.
\end{equation*}
Using the compactness of the level sets of $u \mapsto \norm{u}_{L^\infty} + \tv(u)$ in $L^2(\Omega)$,  one can modify the proof of \cref{thm:asymptotic_profiles} so that its statement still holds.

Therefore, thanks to \labelcref{ineq:poincare_TV}, the total variation flow in one or two space dimension has a finite extinction time. 
Furthermore, it also admits asymptotic profiles.

It makes sense to associate the total variation flow with homogeneous Neumann boundary conditions since it is mass preserving. 
Indeed, since the projection of $u$ onto the space of constant functions, which is the nullspace of $\tv$, is given by the mean $u_\Omega$, we know from \cref{prop:conserv_mass} that the TV flow is mass preserving.

Note that since we are studying a Neumann type boundary condition, we cannot expect uniqueness of the profile.
From our results we get the following theorem, first proved in \cite{andreu2002some}.
\begin{theorem}
Let $\Omega\subset\R^n$ with $n\in\{1,2\}$ be a domain such that the Poincar\'e inequality \labelcref{ineq:poincare_TV} is satisfied.
Let $u:[0,\infty)\to L^2(\Omega)$ solving
\begin{align*}
    \begin{cases}
    \partial_t u = -\div\left(\frac{\grad u}{\abs{\grad u}}\right),\quad&\text{in }\Omega,\\
    \frac{\partial u}{\partial\nu} = 0,\quad&\text{on }\partial\Omega,
    \end{cases}
\end{align*}
denote the gradient flow of \labelcref{eq:TV} with initial datum $f\in L^\infty(\Omega)$.
The extinction time $\tex$ is finite.
Furthermore, there exists a non-decreasing sequence $(t_n)_{n\in\N}\subset\R_+$ with $t_n \to \tex$ such that $w(t_n):=\frac{u(t_n)-u_\Omega}{\norm{u(t_n)-u_\Omega}_{L^{2}(\Omega)}}$ converges to $w\in L^{2}(\Omega)$ as $t\to\tex$, which solves $\lambda w=-\div\left(\frac{\grad u}{\abs{\grad u}}\right)$.
\end{theorem}

One can also study the TV flow with homogeneous Dirichlet boundary conditions, by using the absolutely $1$-homogeneous functional
\begin{align*}
    \calJ(u) := 
    \begin{cases}
    \tv(u) + \int_{\partial\Omega}\abs{u}\d\mathcal{H}^{n-1},\quad&\text{if }u\in BV(\Omega),\\
    \infty,\quad&\text{else}.
    \end{cases}
\end{align*}
This functional is finite for $u\in BV(\Omega)$ since BV functions have integrable boundary traces (see \cite[Theorem 3.87, 3.88]{ambrosio2000functions}).
Quite rightly, one might ask why one does not just impose $u=0$ on $\partial\Omega$.
The reason is that the trace operator $u\mapsto u\vert_{\partial\Omega}$ is \emph{not continuous} w.r.t. the weak-star topology on $BV(\Omega)$.
Therefore, one has to relax the boundary conditions by penalizing the norm of $u$ on $\partial\Omega$.
For more information about the Dirichlet total variation flow we refer the interested reader to \cite{kawohl2007dirichlet}.


\chapter{Nonlinear Power Methods}

In this section we shall discuss how to solve numerically gradient flows and compute nonlinear eigenfunctions. 
For this purpose we will introduce a time implicit Euler discretization.
It turns out that under this discretization the normalizations of the gradient flow, which we have showed to converge to eigenfunctions in \cref{sec:asymptotic_profiles}, solve a nonlinear power method.


Consider an implicit Euler discretization of the gradient flow~\eqref{eq:grad-flow}
\begin{subequations}\label{eq:grad-flow-euler}
\begin{align}[left=\empheqlbrace]
    & \frac{u_{k+1} - u_k}{\tau_k} + \subgradA_{k+1} = 0, \quad k \geq 0, \\
    & \subgradA_{k+1} \in \dJ(u_{k+1}), \quad k \geq 0, \\
    & u^0 = f,
\end{align}
\end{subequations}
where $\tau_k>0$, $k \geq 0$, are the step sizes which we allow to differ from iteration to iteration.
Using the language of proximal mappings, we can rewrite this as follows
\begin{subequations}\label{eq:grad-flow-prox}
\begin{align}[left=\empheqlbrace]
    & u_{k+1} = \argmin_{u \in \Hspace} \frac12 \norm{u-u_k}^2 + \tau_k \calJ(u) \defeq \prox_{\tau_k\calJ}(u_k), \quad k \geq 0, \\
    & u^0 = f.
\end{align}
\end{subequations}
In particular, this shows that the implicit time discretization gives rise to a fixed point iteration for the proximal operator (which has unit Lipschitz constant).
Note that, just as the gradient flow, also this time discretization is mass preserving with $P_{\nullspace_\calJ} u_{k+1}=P_{\nullspace_\calJ} f$.
Furthermore, the sequence $u_k$ converges to $u_\infty:=P_{\nullspace_\calJ} f$.
This enables us to study the rescalings
\begin{align}\label{eq:rescaling_discrete}
    w_k := \frac{u_k - u_\infty}{\norm{u_k-u_\infty}}.
\end{align}

The following statement assures that the sequence $w_k$ satisfies a power iteration for the proximal operator.
\begin{proposition}
Let $\calJ \colon \Hspace \to \RI$ be a proper, convex, \lsc{,} and absolutely $p$-homogeneous functional, and assume that $f\in\Hspace$ satisfies $f\neq P_{\nullspace_\calJ} f$.
Then the rescalings \labelcref{eq:rescaling_discrete} satisfy
\begin{subequations}\label{eq:grad-flow-power}
\begin{align}[left=\empheqlbrace]
    w_{k+1} &= \frac{\prox_{\sigma_k\calJ}(w_k)}{\norm{\prox_{\sigma_k\calJ}(w_k)}}, \quad k \geq 0, \\
    w_0 &= \frac{f-u_\infty}{\norm{f-u_\infty}},
\end{align}
\end{subequations}
where $\sigma_k \defeq \tau_k\norm{u_k-u_\infty}^{p-2}$ and $\tau_k$ are are the step sizes from~\eqref{eq:grad-flow-euler}. 
\end{proposition}
\begin{proof}
Thanks to \cref{prop:subdiff_nullspace} it holds that $\prox_{\tau_k\calJ}(u_k) - u_\infty= \prox_{\tau_k\calJ}(u_k-u_\infty)$.
Hence, we can compute
\begin{align*}
    u_{k+1} - u_\infty 
    &= \prox_{\tau_k\calJ}(u_k) - u_\infty 
    = \prox_{\tau_k\calJ}(u_k-u_\infty) 
    \\
    &=
    \prox_{\tau_k\calJ}(w_k\norm{u_k-u_\infty})
    \\
    &=
    \argmin\frac{1}{2}\norm{u-w_k\norm{u_k-u_\infty}}^2 + \tau_k\calJ(u)
    \\
    &=
    \argmin\frac{1}{2}\norm{u_k-u_\infty}^2\norm{\frac{u}{\norm{u_k-u_\infty}}-w_k}^2 + \tau_k\calJ\left(u\right)
    \\
    &=
    \argmin\frac{1}{2}\norm{\frac{u}{\norm{u_k-u_\infty}}-w_k}^2 + \underbrace{\tau_k\norm{u_k-u_\infty}^{p-2}}_{=\sigma_k}\calJ\left(\frac{u}{\norm{u_k-u_\infty}}\right)
    \\
    &=
    \norm{u_k-u_\infty}\prox_{\sigma_k\calJ}(w_k).
\end{align*}
This yields
\begin{align*}
    w_{k+1} = \frac{u_{k+1}-u_\infty}{\norm{u_{k+1}-u_\infty}} = \frac{\norm{u_k-u_\infty}\prox_{\sigma_k\calJ}(w_k)}{\norm{\norm{u_k-u_\infty}\prox_{\sigma_k\calJ}(w_k)}}
    =
    \frac{\prox_{\sigma_k\calJ}(w_k)}{\norm{\prox_{\sigma_k\calJ}(w_k)}}.
\end{align*}
\end{proof}

To make sure that the proximal power method \labelcref{eq:grad-flow-power} is well defined, meaning that $\prox_{\sigma_k\calJ}(w_k)\neq 0$, we have to choose suitable step sizes $\sigma_k$.
Before proving the well-definedness result in \cref{prop:decrease-J} we need a lemma which gives a sufficient condition for the parameters $\sigma$ such that $\prox_{\sigma\calJ}(f)\neq P_{\nullspace_\calJ}(f)$.
This should be compared to the bounds for the extinction times of gradient flows obtained in \cref{ss:extinction_times}.
\begin{lemma}\label{lem:extinction_prox}
    If $f\in\H$ and $\sigma<\frac{\norm{f-P_{\nullspace_\calJ} f}^2}{\calJ(f)}$, then $\prox_{\sigma\calJ}(f)\neq P_{\nullspace_\calJ} f$.
\end{lemma}
\begin{proof}
    Assuming the contrary, the optimality conditions would imply
    \begin{align*}
        0 = P_{\nullspace_\calJ} f-f + \sigma \subgradA
    \end{align*}
    for some $\subgradA\in\partial\calJ(f)$.
    Consequently, by the definition of the subdifferential we would have
    \begin{align*}
        \sigma
        \geq 
        \sigma
        \sup_{u\in\H}
        \frac{\langle\subgradA,u-P_{\nullspace_\calJ} f\rangle}{\calJ(u)}
        =
        \sup_{u\in\H}
        \frac{\langle f-P_{\nullspace_\calJ} f,u-P_{\nullspace_\calJ} f\rangle}{\calJ(u)}
        \geq
        \frac{\norm{f-P_{\nullspace_\calJ} f}^2}{\calJ(f)}
    \end{align*}
    which is a contradiction.
\end{proof}

\begin{proposition}\label{prop:decrease-J}
Let $\calJ \colon \Hspace \to \RI$ be a proper, convex, \lsc{,} and absolutely $p$-homogeneous functional, and assume that $f\in\Hspace$ satisfies $f\neq P_{\nullspace_\calJ} f$.
Let the step sizes $\sigma_k$ be chosen as
\begin{align*}
    \sigma_k := \frac{c}{\calJ(w_0)}\qquad\text{or}\qquad\frac{c}{\calJ(w_k)},\qquad k\geq 0
\end{align*}
where $0<c<1$.
Then the proximal power method \labelcref{eq:grad-flow-power} is well-defined and satisfies
\begin{align}
    \norm{w_k} = 1,\quad\calJ(w_{k+1}) \leq \calJ(w_k),\qquad \forall k\geq 0.
\end{align}
\end{proposition}
\begin{proof}
Since by assumption $w_0$ is well-defined and satisfies $\norm{w_0}=1$ we can argue inductively. 
Assume that $w_k$ is well-defined and satisfies $\norm{w_k}=1$ as well as $\calJ(w_{k})\leq\calJ(w_{k-1})\leq\dots\leq\calJ(w_0)$.
Let us define
$$
v_k \defeq \prox_{\sigma_k\calJ}(w_k) = \argmin_{v \in \Hspace} \frac12 \norm{v-w_k}^2 + \sigma_k \calJ(v).
$$
We claim that $v_k\neq 0$.
Note that for both choices of step sizes we have
\begin{align*}
    \sigma_k \leq  
    \frac{c}{\calJ(w_k)}
    <
    \frac{\norm{w_k}^2}{\calJ(w_k)}.
\end{align*}
Taking into account that $w_k\in\nullspace_\calJ^\perp$ we can utilize \cref{lem:extinction_prox} to deduce that $v_k\neq P_{\nullspace_\calJ} w_k = 0$ which makes $w_{k+1}$ well-defined and satisfy $\norm{w_{k+1}}=1$.

The optimality of $v_k$ implies that 
\begin{align*}
    \frac{1}{2}\norm{v_k-w_k}^2+\sigma_k\calJ(v_k)\leq\frac{1}{2}\norm{v-w_k}^2+\sigma_k\calJ(v),\quad\forall v\in\H,
\end{align*}
with equality if and only if $v=v_k$, due to the strict convexity of the objective functional.
We define $v:=\frac{\norm{v_k}}{\norm{w_k}}w_k$ which satisfies
\begin{align*}
    \norm{v-w_k}^2 =\norm{v}^2-2\langle v,w_k\rangle+\norm{w_k}^2
    =\norm{v_k}^2-2\norm{v_k}\norm{w_k}+\norm{w_k}^2
    \leq\norm{v_k-w_k}^2.
\end{align*}
Plugging this into the optimality above yields
\begin{align*}
    \sigma_k\calJ(v_k) \leq \sigma_k\calJ(v)=\sigma_k\frac{\norm{v_k}^p}{\norm{w_k}^p}\calJ(w_k).
\end{align*}
Since $v_k\neq 0$, we can divide by $\norm{v_k}$, cancel $\sigma_k$, and arrive at 
$$
\calJ(w_{k+1})=\frac{\calJ(v_k)}{\norm{v_k}^p} \leq \frac{\calJ(w_k)}{\norm{w_k}^p}=\calJ(w_k).
$$
\end{proof}

From this proposition it is clear that if the sublevel sets of $\norm{\cdot}^2+\calJ$ are precompact, the power method converges to some $w_*\in\Hspace$.
Furthermore, continuity properties of the proximal operator (which we prove in \cref{lem:convergences} below) show that $w_*$ solves the nonlinear eigenvalue problem $\mu w_*=\prox_{\sigma\calJ}(w_*)$ which is equivalent to $\lambda w_*\in\dJ(w_*)$ with $\lambda=(1-\mu)/(\sigma\mu)$.

\begin{lemma}[Continuity of the proximal operator]\label{lem:convergences}
Let $(u_k)\subset\nullspace_\calJ^\perp$ be a sequence converging to $u_*$, and let $v_k:=\prox_{\sigma_k\calJ}(u^k)$ for $k\in\N$.
If the sublevel sets of $\norm{\cdot}^2+\calJ$ are precompact and the sequence of regularization parameters fulfills $\lim_{k\to\infty}\sigma_k=\sigma_*>0$ then $(v_k)_{k\in\N}$ converges to some $v_*\in\nullspace_\calJ^\perp$ and it holds $v_*=\prox_{\sigma_*\calJ}(u_*)$.
\end{lemma}
\begin{proof}
From the optimality of $v_k$ in the definition of the prox we deduce
$$\frac{1}{2}\norm{v_k-u_k}^2+\sigma_k \calJ(v_k)\leq\frac{1}{2}\norm{v-u_k}^2+\sigma_k \calJ(v),\quad\forall v\in\H.$$
Choosing $v=0$, we can infer
\begin{align}\label{ineq:optimality_vk}
    \limsup_{k\to\infty}\frac{1}{2}\norm{v_k-u_k}^2+\sigma_k\calJ(v_k)\leq\limsup_{k\to\infty}\frac{1}{2}\norm{u_k}^2<\infty,
\end{align}
since $(u_k)$ is a convergent sequence and hence bounded.
By triangle inequality it holds
\begin{align*}
    \norm{v_k}\leq\norm{v_k-u_k}+\norm{u_k}
\end{align*}
which together with \eqref{ineq:optimality_vk} shows that $\limsup_{k\to\infty}\norm{v_k}<\infty$.
Furthermore, since we have assumed that $\lim_{k\to\infty}\sigma_k=\sigma_*>0$, we also get that  $\limsup_{k\to\infty}\calJ(v_k)<\infty$.
Since the sublevel sets of $\norm{\cdot}^2+\calJ$ are precompact, a subsequence of $(v_k)$ converges to some $v_*\in\H$.
Using lower semi-continuity of $\calJ$ and the strong convergences of $(u_k)$ and $(v_k)$ it holds
\begin{align*}
    \frac{1}{2}\norm{v_*-u_*}^2+\sigma_*\calJ(v_*)&\leq\liminf_{k\to\infty}\frac{1}{2}\norm{v_k-u_k}^2+\sigma_k\calJ(u_k)\\
    &\leq\liminf_{k\to\infty}\frac{1}{2}\norm{v-u_k}+\sigma_k\calJ(v)
    =\frac{1}{2}\norm{v-u_*}^2+\sigma_*\calJ(v),\quad\forall v\in\H.
\end{align*}
This shows that $v_*=\prox_{\sigma_*\calJ}(u^*)$.
The same argument shows that in fact every convergent subsequence of $(v_k)$ converges to $v_*$ and hence the whole sequence $(v_k)$ converges to $v_*$.

Finally, note that by the optimality conditions of the proximal operator it holds that $0\in v_k-u_k+\sigma_k\partial\calJ(v_k)$ and therefore $v_k\in\nullspace_\calJ^\perp$, and as a consequence also the limit satisfies $v_*\in\nullspace_\calJ^\perp$.
\end{proof}

We are ready to prove our main theorem.

\begin{theorem}
Let $\calJ \colon \Hspace \to \RI$ be a proper, convex, \lsc{,} and absolutely $p$-homogeneous functional.
Let $(w_k)_{k\in\N}$ with $w_0\in\nullspace_\calJ^\perp$ be generated by \labelcref{eq:grad-flow-power}.
Let the step sizes $\sigma_k$ be chosen as
\begin{align*}
    \sigma_k := \frac{c}{\calJ(w_0)}\qquad\text{or}\qquad\frac{c}{\calJ(w_k)},\qquad k\geq 0
\end{align*}
where $0<c<1$.
Furthermore, assume that the sublevel sets of $\norm{\cdot}+\calJ(\cdot)$ are precompact and the following Poincar\'{e}-type inequality holds,
\begin{align}\label{eq:Poincare-num}
    \norm{u-P_{\nullspace_\calJ} u}^p \leq C_p \calJ(u),\quad\forall u\in\Hspace.
\end{align} 
\sloppy Then $(w_k)_{k\in\N}$ admits a subsequence converging to $w_\infty\in\nullspace_\calJ^\perp$ which satisfies 
$$\mu w_\infty = \prox_{\sigma(w_\infty)\calJ}(w_\infty)$$ 
with $\mu\in(0,1)$, where $\sigma(w_\infty) = \frac{c}{\calJ(w_0)}$ or $\sigma(w_\infty) = \frac{c}{\calJ(w_\infty)}$ for the constant and variable step sizes, respectively.
\end{theorem}
\begin{proof}
    Since by assumption the sublevel sets of $\norm{\cdot}+\calJ(\cdot)$ are precompact and according to \cref{prop:decrease-J} the sequence $(w_k)\subset\nullspace_\calJ^\perp$ fulfills
$$\sup_{k\in\N}\norm{w_k}+\calJ(w_k)\leq 1+\calJ(w_0)<\infty,$$
it admits a convergent subsequence (which we do not relabel).
We denote by $w_\infty\in\nullspace_\calJ^\perp$ its limit and note that it fulfills 
\begin{align*}
    \norm{w_\infty}=1,\qquad \calJ(w_\infty)\leq\liminf_{k\to\infty}\calJ(w_k),
\end{align*}
by the lower semi-continuity of $\calJ$.
For the constant parameter rule it is trivial that $\sigma_k$ converges to some positive value.
Let us therefore study the variable parameter rule $\sigma_k:=c/\calJ(w_k)$ with $0<c<1$. 
According to \cref{prop:decrease-J} it holds that $\sigma_k$ is an non-decreasing sequence; furthermore, due to~\eqref{eq:Poincare-num} it is bounded by $\sigma_k\leq c\cdot C_p$. 
Hence, $\lim_{k\to\infty}\sigma_k=\sigma_\infty$ exists and by lower semi-continuity of $\calJ$ one has
\begin{align}\label{ineq:estimate_alpha*}
    \sigma_\infty=\lim_{k\to\infty}\sigma_k=\frac{c}{\liminf_{k\to\infty}\calJ(w_k)}\leq\frac{c}{\calJ(w_\infty)}.
\end{align}

Applying \cref{lem:convergences} gives that $v_k:=\prox_{\sigma_k\calJ}(w_k)$ converges to some $v_\infty\in\nullspace_\calJ^\perp$ and it holds $v_\infty=\prox_{\sigma_\infty\calJ}(w_\infty)$.
Note that \labelcref{ineq:estimate_alpha*} together with \cref{lem:extinction_prox} implies that $v_\infty\neq 0$.

It remains to be shown that $v_\infty$ and $w_\infty$ are collinear since this implies that $w_\infty$ is an eigenvector. 
To this end, we note that by definition of $v_k$ it holds
\begin{align*}
    \frac{1}{2}\norm{v_k-w_k}^2+\sigma_k\calJ(v_k) 
    \leq\frac{1}{2}\norm{v-w_k}^2+\sigma_k\calJ(v),\quad\forall v\in\H.
\end{align*}
Choosing $v=\norm{v_k}w_k$, expanding the squared norms, and using $\norm{w_k}=1$ yields
\begin{align*}
    \frac{1}{2}\norm{v_k}^2-\langle v_k,w_k\rangle+\frac{1}{2}+\sigma_k\calJ(v_k)
    \leq \frac{1}{2}\norm{v_k}^2-\norm{v_k}+\frac{1}{2}+\sigma_k\norm{v_k}^p\calJ(w_k).
\end{align*}
This can be simplified to
\begin{align*}
    \norm{v_k}-\langle v_k,w_k\rangle \leq \sigma_k\norm{v_k}^p\left(\calJ(w_k)-\frac{1}{\norm{v_k}^p}\calJ(v_k)\right) 
    = \sigma_k\norm{v_k}^p\left(\calJ(w_k)-\calJ(w_{k+1})\right),
\end{align*}
where we recall that $w_{k+1} = \frac{v_k}{\norm{v_k}}$ by~\eqref{eq:grad-flow-power}. We note that the term $\sigma_k\norm{v_k}^p$ is uniformly bounded by some $C>0$ because both $(\sigma_k)$ and $(v_k)$ converge. By telescopic summation we obtain that
\begin{align*}
    \sum_{k=0}^\infty\norm{v_k}-\langle v_k, w_k\rangle \leq C\left(\calJ(w_0)-\calJ(w_\infty)\right).
\end{align*}
Since $\norm{w_k}=1$, clearly, $\norm{v_k}-\langle v_k, w_k\rangle > 0$ for all $k$. Hence, the sequence on the left hand side is non-negative and summable, which  yields  
\begin{align*}
    \lim_{k\to\infty}\norm{v_k}-\langle v_k,w_k\rangle = 0.
\end{align*}
Using the convergences $v_k\to v_\infty$ and $w_k\to w_\infty$, we get that $\norm{v_\infty}=\langle v_\infty,w_\infty\rangle$ which readily implies $w_\infty={v_\infty}/\norm{v_\infty}$.
Hence, we infer that
$$\prox_{\sigma_\infty\calJ}(w_\infty)=v_\infty=\norm{v_\infty}w_\infty=\tilde{\mu} w_\infty$$
with $\tilde\mu:=\norm{v_\infty}>0$.
Note that $\tilde\mu\leq 1$ since the proximal operator has unitary Lipschitz constant and thus meets
$$\tilde\mu\norm{w_\infty}=\norm{\prox_{\sigma_\infty\calJ}(w_\infty)-\prox_{\sigma_\infty\calJ}(0)}\leq\norm{w_\infty}.$$
Indeed, it even holds that $\tilde\mu<1$. To see this, recall that by definition of the proximal operator one has 
\begin{align*}
    \tilde{\mu} w_\infty = \argmin_{v} \frac12 \norm{v-w_\infty}^2 + \sigma_\infty \calJ(v),
\end{align*}
and the optimality condition is given by
\begin{align}\label{eq:OC}
    \frac{1-\tilde\mu}{\sigma_\infty}w_\infty\in\partial \calJ(w_\infty).
\end{align}
If $\tilde\mu=1$, this optimality condition reduces to $0\in\partial \calJ(w_\infty)$ which means $\calJ(w_\infty)=0$.
Due to $\norm{w_\infty}=1$ and \labelcref{eq:Poincare-num} this is impossible.

It remains to be shown that it also holds $\mu w_\infty=\prox_{\sigma(w_\infty)\calJ}(w_\infty)$ for a suitable $\mu\in(0,1)$.
For the constant parameter rule this is trivially true with $\mu=\tilde{\mu}$ since $\sigma(w_\infty)=\sigma_\infty$ in this case.

For the variable step size rule, we rewrite \labelcref{eq:OC} as $\lambda w_\infty\in\partial \calJ(w_\infty)$ where $\lambda:=(1-\tilde{\mu})/\sigma_\infty>0$.
The optimality conditions for $v = \prox_{\sigma(w_\infty)\calJ}(w_\infty)$ are given by
\begin{align*}
    0 =v-w_\infty + \sigma(w_\infty)q,\quad q\in\partial\calJ(v).
\end{align*}
Now let us make the ansatz that $v=\mu w_\infty$ with $\mu\in\R$. 
Plugging this in and using the homogeneity of $\partial\calJ$ as well as the fact that $\lambda w_\infty\in\partial\calJ(w_\infty)$ we get the equation
\begin{align*}
    0 = \mu-1+\sigma(w_\infty)\lambda\abs{\mu}^{p-1}\sign(\mu).
\end{align*}
We call the function on the right side of this equation $f_p(\mu)$ and make a case distinction.

If $p=1$ we get the equation $0=f_1(\mu)=\mu-1+\sigma(w_\infty)\lambda\sign(\mu)$ which has the unique solution $\mu = (1-\sigma(w_\infty)\lambda)_+$.
By definition of $\lambda$ and using \labelcref{ineq:estimate_alpha*} we have
\begin{align*}
    \mu = (1-\sigma(w_\infty)\lambda)_+
    =
    \left(1-\frac{\sigma(w_\infty)}{\sigma_\infty}(1-\tilde\mu)\right)_+
    \leq \tilde\mu < 1,
\end{align*}
where $(\cdot)_+ = \max\{\cdot,0\}$ denotes the positive part. On the other hand, if $\mu=0$ were true, then \cref{lem:extinction_prox} would imply that 
$$\frac{c}{\calJ(w_\infty)}=\sigma(w_\infty)\geq\frac{1}{\calJ(w_\infty)}$$ 
which is a contradiction since $c\in(0,1)$.

If $p>1$ then $\mu\mapsto f_p(\mu)$ is non-decreasing and continuous and, furthermore, satisfies $f_p(0)=-1<0$ and $f_p(1)=\sigma(w_\infty)\lambda>0$. 
By the mean value theorem there has to exist $\mu\in(0,1)$ with $f_p(\mu)=0$ and we can conclude the proof.
\end{proof}

\begingroup
\sloppy
\printbibliography[heading=bibintoc]
\endgroup
\end{document}